\ifx\shlhetal\undefinedcontrolsequence\let\shlhetal\relax\fi
\documentstyle[11pt]{article}

\input amssym.def
\input amssym.tex

\newcommand{\ens}{{\rm Ens}}
\newcommand{\suc}{{\rm succ}}
\newcommand{\posens}{{\rm PosEns}}
\newcommand{\cI}{{\cal I}}
\newcommand{\cJ}{{\cal J}}
\newcommand{\cP}{{\cal P}}
\newcommand{\cA}{{\cal A}}
\newcommand{\inc}{{\rm inc}}
\newcommand{\V}{{\bf V}}
\newcommand{\bd}{{\rm bd}}
\newcommand{\lh}{\ell g}
\newcommand{\vare}{\varepsilon}
\newcommand{\cf}{{\rm cf}}
\newcommand{\cl}{{\rm cl}}
\newcommand{\dom}{{\rm dom}}
\newcommand{\rang}{{\rm rang}}
\newcommand{\BA}{{\rm BA}}
\newcommand{\inter}{{\rm inter}}
\newcommand{\Length}{{\rm Length}}
\newcommand{\lev}{{\rm lev}}
\newcommand{\rest}{\restriction}
\newcommand{\otp}{{\rm otp}}

\newcommand{\tcf}{{\rm tcf}}
\newcommand{\pcf}{{\rm pcf}}
\newcommand{\ga}{{\frak a}}
\newcommand{\gb}{{\frak b}}
\newcommand{\gc}{{\frak c}}
\newcommand{\gd}{{\frak d}}
\newcommand{\goe}{{\frak e}}

\newcommand{\mod}{{\rm mod}}
\newcommand{\dens}{{\rm dens}}
\newcommand{\cov}{{\rm cov}}
\newcommand{\nacc}{{\rm nacc}}
\newcommand{\pp}{{\rm pp}}
\newcommand{\Reg}{{\rm Reg}}

\newcommand{\QED}{\hfill\vrule width 6pt height 6pt depth 0pt\vspace{0.1in}} 
\newcommand{\Proof}{\noindent{\sc Proof} \hspace{0.2in}} 

\newtheorem{theorem}{Theorem}[section] 

\newtheorem{claim}{Claim}[theorem]
\newtheorem{sdef}[claim]{Definition}
\newtheorem{observation}[claim]{Observation}

\newtheorem{lemma}[theorem]{Lemma} 
\newtheorem{proposition}[theorem]{Proposition} 
 
\newtheorem{definition}[theorem]{Definition}
\newtheorem{conclusion}[theorem]{Conclusion}
\newtheorem{remark}[theorem]{Remark}
\newtheorem{notation}[theorem]{Notation}

\title{$\sigma$-Entangled Linear Orders and Narrowness of Products of
Boolean Algebras} 
\author{Saharon Shelah\thanks{\ \ Partially supported by the Deutsche
Forschungsgemeinschaft, Grant No. Ko 490/7-1. \ \ \ Publication 462.}\\ 
Institute of Mathematics\\
The Hebrew University\\
91 904 Jerusalem, Israel\\
and\\
Department of Mathematics\\
Rutgers University\\
New Brunswick, NJ 08854, USA}
\date{\today}

\begin{document}
\maketitle

\begin{abstract}
We investigate $\sigma$-entangled linear orders and narrowness of Boolean
algebras. We show existence of $\sigma$-entangled linear orders in many
cardinals, and we build Boolean algebras with neither large chains nor large
pies. We study the behavior of these notions in ultraproducts.
\end{abstract}

\eject

\section*{Annotated Content}

\noindent Section 0:\quad {\bf Introduction.}
\medskip

\noindent Section 1:\quad {\bf Basic properties.}

\noindent [We define $\ens_\sigma$, $\sigma$-entangled (Definition \ref{1.1});
we give their basic properties (\ref{1.2}) and the connection between those
properties of linear orders and (the $\sigma$-completion of) the interval
Boolean algebras (Definition \ref{1.3}) which they generate (\ref{1.4}). We
recall the definition of $\inc^{(+)}(B)$ (see \ref{1.5}) and we state its
properties.  Then we formulate the properties of linear orders required to
have $\inc(B^\sigma/D)>(\inc(B))^\sigma/D$ (\ref{1.7}).]
\medskip

\noindent Section 2:\quad {\bf Constructions for $\lambda=\lambda^{<
\lambda}$.}

\noindent [In \ref{2.2}, assuming $\lambda=2^\mu=\mu^+$ (and
$\diamondsuit_\lambda$ which usually follows), we build some Boolean algebras
derived from a tree, using a construction principle (see \cite{Sh:405}). The
tree is a $\lambda^+$-Aronszajn tree, the derived linear order is locally
$\mu$-entangled (of cardinality $\lambda^+$). Next, in \ref{2.3}, we force a
subtree $T$ of ${}^{\lambda\ge}\lambda$ of cardinality $\lambda^+$, the
derived linear order is $\mu$-entangled (of cardinality $\lambda^+$). It
provides an example of Boolean algebras $B_\sigma$ (for $\sigma<\mu$) with
$\inc(B_\sigma)=\lambda$, $\inc((B_\sigma)^\mu/D)=\lambda^+$ for each uniform
ultrafilter $D$ on $\mu$.]
\medskip

\noindent Section 3:\quad {\bf Constructions Related to pcf Theory.}

\noindent [We give sufficient conditions for $\ens_\sigma(\lambda,\kappa)$
when $\lambda$ can be represented as $\tcf(\prod\limits_i \lambda_i/D)$,
$\lambda_i>\max\pcf(\{\lambda_j: j<i\})$ (see \ref{3.1}). If $2^\kappa\geq\sup
\limits_i\lambda_i$ (and more) we can get $\sigma$-entangled linear
order (\ref{3.2}). Also we can utilize $\ens_\sigma(\lambda_i,\kappa_i)$
(see \ref{3.3}, \ref{3.4}). Now relaying on a generalization of ``$\delta<
\aleph_\delta\ \Rightarrow\ \pp(\aleph_\delta)<\aleph_{|\delta|^{+
4}}$'', we prove that if $\mu=\mu^{<\sigma}$ then for many $\theta\in
[\mu,\aleph_{\mu^{+4}})$ we have $\ens_\sigma(\theta^+,\mu)$ and if
$2^\mu\geq\aleph_{\mu^{+4}}$ also $\sigma$-entangled linear orders of 
cardinality $\theta^+$ (see \ref{3.5}). Hence for each $\sigma$ for a class of
successor cardinals there is a $\sigma$-entangled linear order of cardinality
$\lambda^+$ (see \ref{3.6}).] 
\medskip

\noindent Section 4:\quad {\bf Boolean Algebras with neither pies nor 
chains.}

\noindent [Refining results in Section 3, we get Boolean algebras (again
derived from trees $\bigcup\limits_{i\leq\delta}\prod\limits_{j<i}
\lambda_j$ using $\lambda=\tcf(\prod\limits_i \lambda_i/D)$, but not as
interval Boolean algebras), which have neither large chains nor large pies.
For this we need more on how $\lambda=\tcf(\prod\limits_i \lambda_i/D)$.]
\medskip

\noindent Section 5:\quad {\bf More on entangledness.}

\noindent [In \ref{5.1}, \ref{5.2} we deal with cases $2^{<\lambda}<
2^\lambda$. Then we get finer results from assumption on $\pp(\mu)$'s,
improving Section 3. We also deal with $\pcf(\ga)$, defining $\pcf^{ex}_\kappa
(\ga)=\bigcap\{\pcf(\ga\setminus\gb): \gb\subseteq \ga, |\gb|<\kappa\}$
proving for it parallel of the old theorem and connecting it to entangledness,
mainly: if each $\mu\in\ga$ is $(\lambda,\kappa,2)$-inaccessible, then
$\theta\in\pcf^{ex}_\kappa(\ga)\ \Rightarrow\ \ens(\theta, 2^\kappa)$. We
extract from the proof of \cite[\S 4]{Sh:410} on the existence of entangled
linear orders a statement more relevant to $\pcf$. We lastly prove: for a
singular fix point $\mu$ and $\mu_0<\mu$ there is $\theta^+\in
(\mu,\pp^+(\mu))$ in which there is an entangled linear order of density
$\in(\mu_0,\mu)$ (see \ref{5.7A}(2,3).]
\medskip

\noindent Section 6:\quad {\bf Variants of entangledness in ultraproducts.}

\noindent[We investigate what kinds of entangledness (and $\inc(-)\leq\mu$)
are preserved by ultraproducts (\ref{6.5}). We also find that entangledness
can be destroyed by ultrapowers with little connection to its structure, just
its cardinality, for non separative ultrafilters. So to show the possibility
of $(\inc(B))^\omega/D>\inc(B^\omega/D)$ it suffices to find $B=
\BA_{\inter}(\cI)$ such that $|B|>(\inc(B))^{\aleph_0}$.]

\section{Introduction}
In the present paper we investigate $\sigma$-entangled linear orders and
narrowness of Boolean algebras (if $B$ is the interval Boolean algebra of a
linear order $\cI$, then the algebra $B$ is narrow if and only if $\cI$ is
entangled). On entangled $=\aleph_0$--entangled (= narrow interval Boolean
algebra) linear orders (Definition \ref{1.1}(4)) see Bonnet \cite{Bo},
Abraham, Shelah \cite{AbSh:106}, Abraham, Rubin, Shelah \cite{ARSh:153},
Bonnet, Shelah \cite{BoSh:210}, Todorcevic \cite{To} and \cite[\S4]{Sh:345},
\cite{Sh:345b}, \cite[4.9--4.14]{Sh:355}, \cite[\S4]{Sh:410}. 

We prove that for many cardinals $\lambda$ there is a $\sigma$--entangled
linear order of cardinality $\lambda$ (see \ref{3.6}). For example, if
$\lambda$ is a limit cardinal, $\lambda=\lambda^{<\sigma}$, $2^\lambda>
\lambda^{+\lambda^{+4}}$ then for some singular cardinal $\mu\in [\lambda,
\aleph_{\lambda^{+4}})$ there is one in $\mu^+$. We also prove that for a
class of cardinals $\lambda$, there is a Boolean algebra $B$ of cardinality
$\lambda^+$ with neither a chain of cardinality $\lambda^+$ nor a pie (= set
of pairwise incomparable elements) of cardinality $\lambda^+$, see \ref{4.3}. 

Another focus is a problem of Monk \cite{M1}: for a Boolean algebra $B$, let
$\inc(B)$ be $\sup\{|X|: X\subseteq B$ is a pie (see above)$\}$. He asked: are
there a Boolean algebra $B$, a cardinal $\sigma$ and an ultrafilter $D$ on
$\sigma$ such that $\inc(B^\sigma/D)>(\inc(B))^\sigma/D$, and we may ask
whether this holds for $\sigma$ but for no smaller $\sigma'<\sigma$. Now, if
$\cI$ is a $\sigma$--entangled linear order of cardinality $\lambda^+$,
$\lambda^\sigma=\lambda$ then we get examples: the interval Boolean algebra
$B$ of $\cI$ satisfies $\inc(B)=\lambda$ (hence $(\inc B)^\sigma/D=\lambda$)
but in the cases we construct $\cI$, we get that for any uniform ultrafilter
$D$ on $\sigma$, $\inc(B^\sigma/D)=\lambda^+$ (on sufficiency see \ref{1.7};
on existence see \ref{2.2}(3), \ref{3.2}, \ref{3.5}(3)). Similarly for the 
entangledness of a linear order. Unfortunately, though we know that there are
$\sigma$--entangled linear orders of cardinality $\lambda^+$ for many
cardinals $\lambda$ (as needed), we do not know this for cardinals $\lambda$
satisfying $\lambda=\lambda^\sigma$ (even $\lambda^{\aleph_0}=\lambda$), and
$\lambda<\lambda^\sigma$ implies usually $(\inc (B))^\sigma/D\geq\lambda^+$. 
Still, the unresolved case requires quite peculiar cardinal arithmetic
(everywhere): ``usually'' $2^\lambda$ is not so large in the aleph sequence,
and there are additional strong restrictions on the power structure in
$\V$. For instance, for every $\mu$ 
\[\mu^\sigma=\mu\quad \Rightarrow\quad 2^\mu <\aleph_{\mu^{+4}}\]
and 
\[\mu\mbox{ is strong limit of cofinality} >\sigma\quad\Rightarrow\quad
2^\mu<\mu^{+\mu}\ \&\ (\exists\chi)(\chi<\chi^\sigma=2^\mu)\]
and 
\[\mu\geq\beth_\omega\quad\Rightarrow\quad 2^\mu>\mu^+.\]
To make the paper more self-contained we give fully the straight
generalizations of \cite{Sh:345}, \cite{Sh:355} and \cite{Sh:410}. The
research is continued in Magidor Shelah \cite{MgSh:433}, Shafir Shelah
\cite{SaSh:553}, Ros{\l}anowski Shelah \cite{RoSh:534}, \cite{RoSh:599},
and lately \cite{Sh:620}.
\bigskip

We thank Andrzej Ros{\l}anowski and Opher Shafir for reading, correcting,
pointing out various flaws and writing down significant expansions.
\medskip

\noindent{\bf Notation}\ \ \ \ Our notation is rather standard. We will keep
the following rules for our notation: 
\begin{enumerate}
\item $\alpha,\beta,\gamma,\delta,\xi,\zeta, i,j\ldots$ will denote ordinals,
\item $\kappa,\lambda,\mu,\sigma,\ldots$ will stand for cardinal numbers,
\item a bar above a name indicates that the object is a sequence, usually
$\bar{X}$ will be $\langle X_i: i<\lh(\bar{X})\rangle$, where $\lh(\bar{X})$
denotes the length of $\bar{X}$,
\item for two sequences $\eta,\nu$ we write $\nu\vartriangleleft\eta$ whenever
$\nu$ is a proper initial segment of $\eta$, and $\nu\trianglelefteq\eta$ when
either $\nu\vartriangleleft\eta$ or $\nu=\eta$. The length of a sequence
$\eta$ is denoted by $\lh(\eta)$.
\end{enumerate}
For a set $A$ of ordinals with no last element, $J^{\bd}_A$ is the ideal of
bounded subsets of $A$. 

\section{Basic Properties}
In this section we formulate basic definitions and prove fundamental
dependencies between the notions we introduce. 

\begin{definition}
\label{1.1}
Let $\lambda,\mu,\kappa,\sigma$ be cardinal numbers.
\begin{enumerate}
\item A sequence $\bar{\cI}=\langle\cI_\vare:\vare<\kappa\rangle$ of linear
orders is {\em $(\mu,\sigma)$--entangled} if
\begin{description}
\item[($\circledast$)] {\em for every} disjoint subsets $u,v$ of $\kappa$ such
that $|u\cup v|<1+\sigma$ and sequences $\langle t^\vare_\alpha:\alpha<\mu
\rangle$ of pairwise distinct elements of $\cI_\vare$ (for $\vare\in u\cup
v$), {\em there are} $\alpha<\beta<\mu$ such that
\[\vare\in u\quad \Rightarrow\quad t^\vare_\alpha<_{\cI_\vare} t^\vare_\beta
\qquad\mbox{ and }\qquad\vare\in v\quad\Rightarrow\quad t^\vare_\alpha>_{
\cI_\vare} t^\vare_\beta.\]
\end{description}
$\ens(\lambda,\mu,\kappa,\sigma)=\ens_\sigma(\lambda,\mu,\kappa)$ means:
\begin{quotation}
\noindent there is a $(\mu,\sigma)$--entangled sequence
$\bar{\cI}=\langle\cI_\vare: \vare<\kappa\rangle$ of linear orders, each of
cardinality $\lambda$. 
\end{quotation}
\item If we omit $\mu$, this means $\lambda=\mu$ (i.e. $|\cI_\vare|=\mu$), if
we omit $\sigma$ it means $\sigma=\aleph_0$.
\item A linear order $\cI$ is {\em $(\mu,\sigma)$--entangled} if ($\cI$ has
cardinality $\geq\mu$ and) for every $\vare(*)<\sigma$ and a partition $(u,v)$ of $\vare(*)$ and pairwise distinct $t^\vare_\alpha\in\cI$ (for $\vare\in u\cup
v$, $\alpha<\mu$), {\em there are} $\alpha<\beta<\mu$ such that
\begin{description}
\item[($\oplus$)]  for each $\vare<\vare(*)$ we have 
\[\vare\in u\quad\Rightarrow\quad t^\vare_\alpha<_{\cI}t^\vare_\beta\qquad
\mbox{ and }\qquad\varepsilon\in v\quad \Rightarrow\quad t^\vare_\alpha>_{\cI}
t^\vare_\beta.\]
\end{description}
\item We omit $\mu$ if $|\cI|=\mu$ (and so we write ``$\cI$ is
$\sigma$--entangled'' instead ``$\cI$ is $(|\cI|,\sigma)$--entangled''); we
omit also $\sigma$ if it is $\aleph_0$.
\item A sequence $\langle\cI_{\zeta}: \zeta<\gamma\rangle$ of linear orders
is {\em strongly $(\mu,\sigma,\sigma')$--entangled} if
\begin{description}
\item[(a)] each of $\cI_\zeta$ is of cardinality $\geq\mu$,
\item[(b)] {\em if} $u,v$ are disjoint subsets of $\gamma$, $|u\cup v|< 1+
\sigma$, $\xi(\vare)<\sigma'$ for $\vare\in u\cup v$ and $t^\alpha_{\vare,\xi}
\in \cI_\vare$ (for $\alpha<\mu$, $\vare\in u\cup v$, $\xi<\xi(\vare)$) are
such that 
\[(\forall\vare\in u\cup v)(\forall\xi,\zeta<\xi(\vare))(\forall\alpha<\beta<
\mu)(t^\alpha_{\vare,\xi}\neq t^\beta_{\vare,\zeta})\]
{\em then}\ \ for some $\alpha<\beta<\mu$ we have:
\[\begin{array}{l}
\vare\in u\quad\Rightarrow\quad (\forall\xi<\xi(\vare))(t^\alpha_{\vare,\xi}
<t^\beta_{\vare,\xi})\qquad\mbox{ and}\\
\vare\in v\quad\Rightarrow\quad (\forall\xi<\xi(\vare))(t^\beta_{\vare,\xi}<
t^\alpha_{\vare,\xi}).
  \end{array}\]
\end{description}
\end{enumerate}
\end{definition}

\begin{proposition}
\label{1.2}
\begin{enumerate}
\item Assume $\lambda\geq\lambda_1\geq\mu_1\geq\mu$, $\kappa_1\leq\kappa$ and
$\sigma_1\leq\sigma$. Then $\ens_\sigma(\lambda,\mu,\kappa)$ implies
$\ens_{\sigma_1}(\lambda_1,\mu_1,\kappa_1)$.
\item If $\cI$ is a $(\mu,\sigma)$--entangled linear order, $\cJ\subseteq\cI$,
and $|\cI|\geq |\cJ|\geq\mu_1\geq\mu$, $\sigma_1\leq\sigma$ then $\cJ$ is
$(\mu_1,\sigma_1)$--entangled. 
\item If a linear order $\cI$ has density $\chi$, $\chi^{<\sigma}<\mu$,
$\mu=\cf(\mu)$ and $\sigma\geq 2$ 

then in Definition \ref{1.1}(3) of ``$\cI$ is $(\mu,\sigma)$--entangled'' we
can add to the assumptions
\begin{description}
\item[($\circledcirc$)] there is a sequence $\langle [a_\vare,b_\vare]:\vare<
\vare(*)\rangle$ of pairwise disjoint intervals of $\cI$ such that
$t^\vare_\alpha\in (a_\vare,b_\vare)$. 
\end{description}
\item Moreover, if a linear order $\cI$ has density $\chi$, $\chi^{<\sigma}<
\mu=\cf(\mu)$ 

then for each $\vare(*)<\sigma$ and sequences $\bar{t}_\alpha=\langle
t^\vare_\alpha: \vare<\vare(*)\rangle\subseteq\cI$ (for $\alpha<\mu$) such
that $\vare\neq \zeta\ \Rightarrow\ t^\vare_\alpha\neq t^\zeta_\alpha$, there
are $A\subseteq\mu$, $|A|=\mu$, and a sequence $\langle [a_\vare,b_\vare]:
\vare<\vare(*)\rangle$ of pairwise disjoint intervals of $\cI$ such that for
each $\vare<\vare(*)$
\[\mbox{either }\ (\forall\alpha\in A)(t^\vare_\alpha\in (a_\vare,b_\vare))\ 
\mbox{ or }\ (\forall\alpha\in A)(t^\vare_\alpha= a_\vare).\]
\item If $\sigma\geq 2$ and a linear order $\cI$ is $(\mu,\sigma)$--entangled
then $\cI$ has density $<\mu$.
\item If there exists a $(\mu,\sigma)$--entangled linear order of size
$\lambda$ then $\ens_\sigma(\lambda,\mu,\lambda)$.  
\item In Definition \ref{1.1}(3), if $\sigma$ is infinite, we can weaken
``$\alpha<\beta<\mu$'' to ``$\alpha\neq\beta$, $\alpha<\mu$, $\beta<\mu$''.
\item If there is a $(\mu,\sigma)$--entangled linear order of size $\lambda$
and $(*)_\kappa$ below holds then $\ens_\sigma(\lambda,\mu,\kappa)$, where:
\begin{description}
\item[$(*)_\kappa$] one of the following holds true:
\begin{description}
\item[($\alpha$)]\ \ $\kappa=\mu^+$ and if $\lambda=\mu$ then $\cf(\mu)\geq
\sigma$, 
\item[($\beta$)]\ \  there are $A_i\subseteq\lambda$ for $i<\kappa$, $|A_i|=
\lambda$ such that $i\neq j\ \Rightarrow\ |A_i\cap A_j|<\mu$ and $\cf(\mu)
\geq\sigma$, 
\item[($\gamma$)]\ \ there are $A_i\subseteq\lambda$ for $i<\kappa$, $|A_i|=
\lambda$ such that $\sup\{|A_i\cap A_j|: i<j<\kappa\}<\mu$.
\end{description}
\end{description}
\end{enumerate}
\end{proposition}

\Proof 1), 2) are left to the reader.
\medskip

\noindent 3) Clearly the new definition is weaker, so we shall prove that the
one from \ref{1.1}(3) holds assuming the one from \ref{1.2}(3). Let
$\cJ\subseteq\cI$ be dense in $\cI$ and $|\cJ|\leq\chi$. Thus for each $a,b\in
\cI$, $a<_{\cI}b$, there exists $s\in\cJ$ such that $a\leq_{\cI} s\leq_{\cI}
b$. 

Suppose that $\vare(*)$, $u,v$, $\langle t^\vare_\alpha: \vare<\vare(*),\alpha
<\mu\rangle$ are as in \ref{1.1}(3). For each $\vare,\zeta<\vare(*)$ and
$\alpha<\mu$ such that $t^\vare_\alpha<t^\zeta_\alpha$ there exists $s^{\vare,
\zeta}_\alpha\in\cJ$ such that $t^\vare_\alpha\leq s^{\vare,\zeta}_\alpha\leq
t^\zeta_\alpha$ (and at least one inequality is strict). Define functions
$h_0, h_1, h_2, h_3$ on $\mu$ by:
\[\begin{array}{lll}
h_0(\alpha)&=&\{\langle\vare,\zeta\rangle:\vare,\zeta<\vare(*)\mbox{ and }
t^\vare_\alpha<t^\zeta_\alpha\}\\ 
h_1(\alpha)&=&\langle s^{\vare,\zeta}_\alpha:\langle\vare,\zeta\rangle\in
h_0(\alpha)\rangle\\ 
h_2(\alpha)&=&\langle\langle\vare,\zeta,\xi,{\rm TV}(t^\xi_\alpha=s^{\vare,
\zeta}_\alpha)\rangle:\langle\vare,\zeta\rangle\in h_0(\alpha),\ \xi<\vare(*)
\rangle\\
h_3(\alpha)&=&\langle\langle\vare,\zeta,\xi,{\rm TV}(t^\xi_\alpha<s^{\vare, 
\zeta}_\alpha)\rangle:\langle\vare,\zeta\rangle\in h_0(\alpha),\ \xi<\vare(*)
\rangle,
  \end{array}\]
(where ${\rm TV}(-)$ is the truth value of $-$). Now, for each $\ell<4$,
$\dom(h_\ell)=\mu$ and $|\rang(h_\ell)|\leq |\cJ|^{|\vare(*)|^3}\leq
\chi^{<\sigma}<\mu$. Since $\cf(\mu)=\mu$, there exists $A\in
[\mu]^{\textstyle \mu}$ such that the restrictions $h_\ell\restriction A$ are
constant for $\ell = 0,1,2,3$. 

\noindent So let $s^{\vare,\zeta}_\alpha=s^{\vare,\zeta}$ for $\alpha\in A$. 
As the $t^{\vare}_\alpha$'s were pairwise distinct (for each $\vare$) we
conclude 
\[(\alpha\in A\ \& \ \vare<\vare(*))\quad \Rightarrow\quad t^\vare_\alpha
\notin\{s^{\xi,\zeta}:\langle\xi,\zeta\rangle\in h_0(\alpha)\}.\]
Define for $\vare<\vare(*)$:
\[\begin{array}{ll}
\cI_\vare=\{t\in\cI:&\mbox{for every }\zeta,\xi<\vare(*)\mbox{ such that }
s^{\xi,\zeta}\mbox{ is well defined and}\\
\ &\mbox{for every ($\equiv$ some) }\alpha\in A\mbox{ we have }\\
\ &t\leq s^{\xi,\zeta}\Leftrightarrow t^\vare_\alpha\leq s^{\xi,\zeta}\quad
\mbox{ and }\quad t\geq s^{\xi,\zeta}\Leftrightarrow t^\vare_\alpha\geq
s^{\xi,\zeta}\}.\\ 
  \end{array}\]
Note that the value of $\alpha$ is immaterial.\\
Now, easily $\cI_\vare$ does not have cofinality $>\chi$ (as $\cI$ has no
monotonic sequence of length $\geq\chi^+$, remember $\cI$ has density $\leq
\chi$). Hence we find an unbounded well ordered subset $\cJ^+_\vare\subseteq
\cI_\vare$, $|\cJ^+_\vare|\leq\chi$. Similarly there is an anti-well ordered
$\cJ^-_\vare\subseteq\cI_\vare$, $|\cJ^-_\vare|\leq\chi$, which is unbounded
from below (in $\cI_\vare$). Let $\cJ^*=\bigcup\limits_{\vare<\vare(*)}(
\cJ^+_\vare\cup \cJ^-_\vare)$. Again, for some set $A'\subseteq A$ of size
$\mu$, the Dedekind cut which $t^\vare_\alpha$ realizes in $\cJ^*$ does not
depend on $\alpha$ for $\alpha\in A'$, and $t^\vare_\alpha\notin\cJ^*$. Now we
can easily choose $(a_\vare,b_\vare)$: $a_\vare$ is any member of
$\cJ^-_\vare$ which is $<t^\vare_\alpha$ for $\alpha\in A'$ and $b_\vare$ is
any member of $\cJ^+_\vare$ which is $>t^\vare_\alpha$ for $\alpha\in A'$.
\medskip

\noindent 4) Included in the proof of \ref{1.2}(3).
\medskip

\noindent 5) By \ref{1.2}(2), without loss of generality $\sigma=2$. Suppose
that $\cI$ has density at least $\mu$. By induction on $\alpha<\mu$ we try to
choose $t^0_\alpha,t^1_\alpha$ such that
\begin{description}
\item[(i)]\ \ \ $t^0_\alpha<t^1_\alpha$,
\item[(ii)]\ \  $t^0_\alpha,t^1_\alpha\notin\{t^0_\beta,t^1_\beta:\beta<\alpha
\}$,
\item[(iii)]\   $(\forall\beta<\alpha)(\forall\ell\in\{0,1\})(t^0_\alpha<
t^\ell_\beta\Leftrightarrow t^1_\alpha<t^\ell_\beta)$.
\end{description}
Continue to define for as long as possible. There are two possible outcomes.

\noindent{\sc Outcome {\bf (a)}}: \ One gets stuck at some $\alpha<\mu$.\\
Let $\cJ=\{t^0_\beta, t^1_\beta:\beta<\alpha\}$. Then
\[(\forall t^0<t^1\in\cI\setminus\cJ)(\exists s\in\cJ)(t^0<s\Leftrightarrow
\neg[t^1<s]).\]
Since $t^0,t^1\notin\cJ$, it follows that $t^0<s<t^1$. So $\cJ$ is dense in
$\cI$ and is of cardinality $2|\alpha|<\mu$ -- a contradiction.

\noindent {\sc Outcome {\bf (b)}}: \ One can define $t^0_\alpha, t^1_\alpha$
for every $\alpha<\mu$.\\
Then $\langle t^0_\alpha, t^1_\alpha:\alpha<\mu\rangle$, $u=\{1\}$, $v=\{0\}$
constitute an easy counterexample to the $(\mu,2)$--entangledness of $\cI$.
\medskip
 
\noindent 6) Suppose $\cI$ is $(\mu,\sigma)$--entangled and $|\cI|=\lambda$. 
Take a sequence $\langle\cI_\vare:\vare<\lambda\rangle$ of pairwise disjoint
subsets of $\cI$, each of power $\lambda$. This sequence witnesses
$\ens_\sigma(\lambda,\mu,\lambda)$: suppose $u,v$ are disjoint subsets of
$\lambda$, $|u\cup v|<\sigma$ and let $t^\vare_\alpha\in\cI_\vare$ for $\alpha
<\mu$, $\vare\in u\cup v$ be pairwise distinct. Now apply ``$\cI$ is $(\mu,
\sigma)$--entangled". 
\medskip

\noindent 7) Let $u,v,t^\vare_\alpha$ (for $\vare\in u\cup v$, $\alpha<\mu$)
be as in Definition \ref{1.1}(3). Put
\[\begin{array}{c}
u'=\{2\vare:\vare\in u\}\cup\{2\vare+1:\vare\in v\},\quad v'=\{2\vare:\vare
\in v\}\cup\{2\vare+1:\vare\in u\},\\
s^{2\vare}_\alpha=t^\vare_{2\alpha},\quad s^{2\vare+1}_\alpha=t^\vare_{
2\alpha+1}.
  \end{array}\]
Now we apply Definition \ref{1.1}(3) -- the \ref{1.2}(7) version to $u'$,
$v'$, $\langle s^\vare_\alpha: \vare\in u'\cup v', \alpha<\mu\rangle$, and we
get $\alpha'\neq\beta'$ as required there. If $\alpha'<\beta'$ then $\alpha=
2\alpha'$, $\beta=2\beta'$ are as required in \ref{1.1}(3). Otherwise $\alpha'
>\beta'$ and then $\alpha=2\beta'+1$, $\beta=2\alpha'+1$ are as required in
\ref{1.1}(3).    
\medskip

\noindent 8) ($\alpha$)\ \ \ Suppose $\lambda=\mu$ (and so $\cf(\mu)\geq
\sigma$) and let $\cI$ be a $(\mu,\sigma)$--entangled linear order of size
$\lambda$. Choose a family $\{A_\vare:\vare<\mu^+\}\subseteq [\cI]^{\textstyle
\mu}$ such that $(\forall\vare<\zeta<\mu^+)(|A_\vare\cap A_\zeta|<\mu)$, and
let $\cI_\vare=\cI\rest A_\vare$ (for $\vare<\mu^+$). We claim that the
sequence $\langle \cJ_\vare:\vare<\mu^+\rangle$ witnesses $\ens_\sigma(
\lambda, \mu,\mu^+)$. Why? Clearly $|J_\vare|=\lambda=\mu$. Suppose that $u,v
\subseteq\mu^+$ are disjoint, $|u\cup v|<1+\sigma$ and for $\vare\in u\cup v$
let $\langle t^\vare_\alpha:\alpha<\mu\rangle\subseteq \cJ_\vare$ be pairwise
distinct. Since $\sigma\leq\cf(\mu)$ we find $\alpha(*)<\mu$ such that
\[(\forall \vare_0,\vare_1\in u\cap v)(\forall\alpha_0,\alpha_1<\mu)(\vare_0
\neq \vare_1\ \&\ \alpha_0,\alpha_1>\alpha(*)\quad\Rightarrow\quad t^{\vare_0
}_{\alpha_0}\neq t^{\vare_1}_{\alpha_1})\]
(remember the choice of the $A_\vare$'s). Now apply the assumption that $\cI$
is $(\mu,\sigma)$--entangled to the sequence
\[\langle t^\vare_\alpha: \vare\in u\cup v,\ \alpha\in (\alpha(*),\mu)\rangle
\subseteq\cI.\]
If $\lambda>\mu$ then we can choose a family $\{A_\vare:\vare<\mu^+\}$ of
pairwise disjoint sets from $[\cI]^{\textstyle \lambda}$ and proceed as above.

\noindent ($\beta$), ($\gamma$)\ \ \ Similarly. \QED$_{\ref{1.2}}$

\begin{definition}
\label{1.3}
Let $\cI$ be a linear order.
\begin{enumerate}
\item The interval Boolean algebra $\BA_{\inter}(\cI)$ determined by $\cI$ is
the algebra of finite unions of closed--open intervals of $\cI$ (including
$[-\infty,x)$, $[x,\infty)$, $[-\infty,\infty)$).
\item For a regular cardinal $\sigma$, $\BA^\sigma_{\inter}(\cI)$ is the
closure of the family of subsets of $\cI$ of the form $[-\infty,s)$ (for
$s\in\cI$), by complementation and unions and intersections of $<\sigma$
members\footnote{Equivalently, the Boolean algebra $\sigma$--generated by $\{
x_t:t\in\cI\}$ freely except $x_s\leq x_t$ when $\cI\models s<t$.}.
\end{enumerate}
\end{definition}

\begin{definition}
\label{1.5}
Let $B$ be an infinite Boolean algebra.
\begin{enumerate}
\item A set $Y\subseteq B$ is {\em a pie} if any two members of $Y$ are
incomparable (in $B$; ``pie'' comes from ``a set of {\bf p}airwise {\bf
i}ncomparable {\bf e}lements'').
\item $\inc(B)=\sup\{|Y|: Y\subseteq B$ is a pie $\}$.
\item $\inc^+(B)=\sup\{|Y|^+: Y\subseteq B$ is a pie $\}$. 
\item The algebra $B$ is $\mu$--narrow if there is no pie of cardinality
$\ge\mu$.
\item $\Length(B)=\sup\{|Y|: Y\subseteq B$ is a chain $\}$,

$\Length^+(B)=\sup\{|Y|^+: Y\subseteq B$ is a chain $\}$.
\end{enumerate}
\end{definition}

\begin{proposition}
\label{1.4}
Suppose that $\cI$ is a linear order and that the regular cardinals $\aleph_0
\leq\sigma<\mu$ satisfy $(\forall\theta<\mu)[\theta^{<\sigma}<\mu]$. Then the
following conditions are equivalent: 
\begin{description}
\item[(a)] The order $\cI$ is $(\mu,\sigma)$--entangled.
\item[(b)] If $\vare(*)<\sigma$, and $u,v\subseteq\vare(*)$ are disjoint and
$t^\vare_\alpha\in\cI$ (for $\vare<\vare(*)$, $\alpha<\mu$)

then for some $\alpha,\beta<\mu$ we have 
\[\vare\in u\quad\Rightarrow\quad t^\vare_\alpha\leq_{\cI}t^\vare_\beta
\qquad\mbox{ and }\qquad\vare\in v\quad \Rightarrow\quad t^\vare_\alpha
\geq_{\cI} t^\vare_\beta.\]
(Note: if the $t^\vare_\alpha$ are pairwise distinct then the inequalities are
in fact strict; as in the proof of \ref{1.2}(7) changing the demand $\alpha<
\beta$ to $\alpha\neq\beta$ does not matter.)
\item[(c)] The algebra $\BA^\sigma_{\inter}(\cI)$ is $\mu$--narrow.
\end{description}
\end{proposition}

\Proof {\bf (a)}\quad$\Rightarrow$\quad{\bf (c)}.\qquad By \ref{1.2}(5) the
order $\cI$ has density $<\mu$.\\
Let $\langle A_\alpha: \alpha<\mu\rangle$ be a sequence of distinct elements
of the algebra $\BA^\sigma_{\inter}(\cI)$. We know that for each $\alpha$
there are: an ordinal $\vare_\alpha<\sigma$, a Boolean term $\tau_\alpha$
(with all unions and intersections of size $<\sigma$ and $\vare_\alpha$ free
variables) and a sequence $\langle t^\vare_\alpha: \vare<\vare_\alpha\rangle
\subseteq\cI$ such that $A_\alpha=\tau_\alpha(\ldots,t^\vare_\alpha,\ldots)_{
\vare<\vare_\alpha}$. By \ref{1.2}(4), without loss of generality for some
$\vare(*)$ and pairwise disjoint intervals $[a_\vare, b_\vare]$ we have:
$\vare_\alpha=\vare(*)$ and for each $\vare<\vare(*)$ either $(\forall\alpha<
\mu)(a_\vare<t^\vare_\alpha<b_\vare)$ or $(\forall\alpha<\mu)(t^\vare_\alpha
=a_\vare)$. Since $\mu=\cf(\mu)>\aleph_0$ and $(\forall\theta<\mu)
(\theta^{|\vare(*)|}<\mu)$ we may apply the $\Delta$-lemma to the family
$\{x_\alpha: \alpha<\mu\}$, where $x_\alpha=:\{t^\vare_\alpha:\vare<\vare(*)
\}$. Consequently, we may assume that $\{x_\alpha:\alpha<\mu\}$ forms a
$\Delta$-system with the kernel $x$ (i.e. $\alpha<\beta<\mu\ \Rightarrow\
x_\alpha\cap x_\beta=x)$. Note that if for some $\alpha<\mu$, $t^\vare_\alpha
\in x$ then $(\forall\alpha<\beta<\mu)(t^\vare_\alpha=t^\vare_\beta)$ and if
$t^\vare_\alpha\notin x$ (for some $\alpha<\mu)$ then $(\forall\alpha<\beta<
\mu)(t^\vare_\alpha\neq t^\vare_\beta)$. Thus for each $\vare<\vare(*)$ either
$t^\vare_\alpha$  (for $\alpha<\mu$) are pairwise distinct or they are
pairwise equal. Since $\mu =\cf(\mu)>\sigma$ and $\theta^{<\sigma}<\mu$ for  
$\theta<\mu$, without loss of generality $\tau_\alpha=\tau$. Let 
\[w=\{\vare<\vare(*): \langle t^\vare_\alpha:\alpha<\mu\rangle\mbox{ are
pairwise distinct }\}.\]
Then for some disjoint sets $v,u\subseteq w$ and a set $A \subseteq\cI
\setminus\bigcup\limits_{\vare\in u\cup v}[a_\vare, b_\vare]$ we have
$A_\alpha=A\cup\bigcup\limits_{\vare\in u\cup v}\tau^\vare(a_\vare,
t^\vare_\alpha, b_\vare)$, where we let
\[\tau^\vare(x,y,z)=\left\{
\begin{array}{ll}
[x,y)&\mbox{if }\vare\in u,\\
\relax
[y,z)&\mbox{if }\vare\in v.
\end{array}\right.\]
Since $\cI$ is $(\mu,\sigma)$--entangled, we can find $\alpha<\beta$ such that 
\[(\forall\vare\in u\cup v)(t^\vare_\alpha<t^\vare_\beta\ \Leftrightarrow\
\vare\in u).\]
Clearly this implies that $A_\alpha\subseteq A_\beta$, so we are done.
\medskip

\noindent{\bf (c)}\quad$\Rightarrow$\quad{\bf (a)}.\qquad First we note that
the linear order $\cI$ has density $<\mu$.

\noindent [Why? Easily $\cI$ has no well ordered subset of power $\mu$ nor
an inverse well ordered subset of power $\mu$. Assume $\cI$ has density $\geq
\mu$. First we show that there are disjoint closed--open intervals
$\cI_0,\cI_1$ of $\cI$ with density $\geq\mu$. To prove the existence of
$\cI_0,\cI_1$ define the relation $E$ on $\cI$ by: 
\begin{quotation}
\noindent $a\;E\; b$\qquad if and only if 

\noindent $a=b$\quad or\quad [$a<b$ and ${\rm density}([a,b))<\mu$]\quad or

\noindent [$a>b$ and ${\rm density}([b,a))<\mu$].
\end{quotation}
Clearly $E$ is an equivalence relation and its equivalence classes are
convex. Moreover, the density of each $E$--equivalence class is less than
$\mu$ (as there is no monotonic sequence of length $\mu$ of members of
$\cI$). Consequently we find $a,b\in\cI$ such that $a<b$, $\neg a\;E\;b$. 
Next we can find $c,d\in(a,b)$, $c<d$ such that neither $a\;E\;c$ nor $d\;E\;
b$. Thus we may put $\cI_0=[a,c)$, $\cI_1=[d,b)$. Now for each $\cI_m $ we
choose by induction on $\beta<\mu$ elements $a^m_\beta<b^m_\beta$ from $\cI_m$
such that $[a^m_\beta,b^m_\beta]$ is disjoint from $\{a^m_\alpha,b^m_\alpha:
\alpha<\beta\}$. So $\alpha<\beta\ \Rightarrow\ [a^m_\alpha,b^m_\alpha)\not
\subseteq [a^m_\beta,b^m_\beta)$. Now, $\langle [a^0_\beta,b^0_\beta)\cup(
\cI_1\setminus [a^1_\beta,b^1_\beta)):\beta<\mu\rangle$ shows that the algebra
$\BA^\sigma_{\inter}(\cI)$ is not $\mu$-narrow, a contradiction].
\medskip

\noindent By \ref{1.2}(7) it is enough to prove that if $\vare(*)<\sigma$ and
$t^\vare_\alpha\in\cI$ are distinct for $\alpha<\mu$, $\vare<\vare(*)$ and
$u,v$ are disjoint subsets of $\vare(*)$ then we can find $\alpha\neq\beta$
such that: 
\[\vare\in u\quad\Rightarrow\quad t^\vare_\alpha<t^\vare_\beta,\quad\mbox{ and
}\quad \vare\in v\quad\Rightarrow\quad t^\vare_\alpha > t^\vare_\beta.\]
By \ref{1.2}(4), without loss of generality for some pairwise disjoint
intervals $[a_\vare,b_\vare]$ of $\cI$, we have $t^\vare_\alpha\in (a_\vare,
b_\vare)$. Let $x_\alpha=: x^1_\alpha\cup x^2_\alpha$, where
\[x^1_\alpha=:\bigcup\{[a_\vare,t^\vare_\alpha): \vare\in u\},\qquad
x^2_\alpha=:\bigcup\{[t^\vare_\alpha,b_\vare): \vare\in v\}.\]
So for $\alpha<\mu$, $x_\alpha\in\BA^\sigma_{\inter}(\cI)$. The algebra
$\BA^\sigma_{\inter}(\cI)$ is $\mu$--narrow, so for some $\alpha\neq\beta$
($<\mu$) we have $x_\alpha\subseteq x_\beta$. Then for each $\vare$:
\[\vare\in u\quad\Rightarrow\quad\cI\models t^\vare_\alpha\leq t^\vare_\beta,
\quad\mbox{ and }\quad\vare\in v\quad\Rightarrow\quad\cI\models t^\vare_\alpha
\geq t^\vare_\beta.\]
This is as required.
\medskip

\noindent{\bf (a)}\quad$\Rightarrow$\quad{\bf (b)}.\qquad It is included in
the proof of {\bf (a)}\ $\Rightarrow$\ {\bf (c)}.
\medskip

\noindent{\bf (b)}\quad$\Rightarrow$\quad{\bf (a)}.\qquad Trivial. 
\QED$_{\ref{1.5}}$

\begin{proposition}
\label{1.6}
Let $B$ be a Boolean algebra.
\begin{enumerate}
\item If $\inc^+(B)$ is a successor cardinal then $[\inc(B)]^+=\inc^+(B)$.
\item If $\inc^+(B)$ is a limit cardinal then $\inc(B)=\inc^+(B)$.
\item $B$ is $\mu$--narrow\quad if and only\quad if $\inc^+(B)\leq\mu$.
\item If $B$ is $\mu$--narrow then so is every homomorphic image of $B$. 
\item If $D$ is a filter on $\sigma$ and the product algebra $B^\sigma$ is
$\mu$--narrow then the algebra $B^\sigma/D$ is $\mu$--narrow. 
\QED$_{\ref{1.6}}$ 
\end{enumerate}
\end{proposition}

\begin{conclusion}
\label{1.7}
Assume $\lambda^*\geq\mu$, $\lambda\geq\mu=\cf(\mu)>\kappa\geq\sigma=\cf(
\sigma)\geq\aleph_0$ and $(\forall\theta<\mu)[\theta^{<\sigma}<\mu]$.
\begin{enumerate}
\item Then $(A)_{\lambda^*,\lambda,\mu,\sigma,\kappa}\ \Rightarrow\
(B)_{\lambda^*,\lambda,\mu,\sigma,\kappa}$, using $B_j=\BA_{\inter}(\cI_j+
\cJ_j)$, where
\begin{description}
\item[$(A)_{\lambda^*,\lambda,\mu,\sigma,\kappa}$] there are linear orders
$\cI_j$, $\cJ_j$ (for $j\leq\kappa$) of cardinality $\lambda$, such that each
$\cI_j+\cJ_j$ is $(\mu,\sigma)$--entangled and for any uniform ultrafilter $D$
on $\kappa$, the linear orders $\prod\limits_{j<\kappa}\cI_j/D$ and
$\prod\limits_{j<\kappa}\cJ_j/D$ have isomorphic subsets of cardinality
$\lambda^*$; 
\item[$(B)_{\lambda^*,\lambda,\mu,\sigma,\kappa}$] there are interval Boolean
algebras $B_j$ (for $j<\kappa$) which are $\mu$--narrow and of cardinality
$\lambda$ such that for any uniform ultrafilter $D$ on $\kappa$ the algebra
$B=\prod\limits_{i<\kappa}B_i/D$ is not $\lambda^*$--narrow. 
\end{description}
\item  Also $(A)^+_{\lambda^*,\lambda,\mu,\sigma,\kappa}\ \Rightarrow\ (B)^+_{\lambda^*,\lambda,\mu,\sigma,\kappa}$ (using $B=\BA_{\inter}(\cI+\cJ)$), where
\begin{description}
\item[$(A)^+_{\lambda^*,\lambda,\mu,\sigma,\kappa}$] there are linear orders
$\cI,\cJ$ of cardinality $\lambda$, such that $\cI+\cJ$ is
$(\mu,\sigma)$--entangled and for any uniform ultrafilter $D$ on $\kappa$ the
linear orders $\cI^\kappa /D$, $\cJ^\kappa/D$ have isomorphic subsets of
cardinality $\lambda^*$. 
\item[$(B)^+_{\lambda^*,\lambda,\mu,\sigma,\kappa}$] there is a $\mu$-narrow
interval Boolean algebra $B$ of cardinality $\lambda$ such that $\lambda^*<
\inc^+[B^\kappa/D]$ for any uniform ultrafilter $D$ on $\kappa$ (i.e. the
algebra is not $\lambda^*$--narrow). 
\end{description}
\item We can replace ``uniform ultrafilter $D$'' by ``regular ultrafilter
$D$'' or fix a filter $D$ on $\kappa$. 
\end{enumerate}
\end{conclusion}

\Proof Just note that
\begin{quotation}
\noindent if $B$ is a Boolean algebra, $\cI,\cJ$ are linear orders, $a_t\in B$
for $t\in\cI+\cJ$ are such that $t<s\ \Rightarrow\ a_t<_B a_s$ and $f$ is an (order) isomorphism from $\cI$ to $\cJ$

\noindent then $\{a_{f(t)}-a_t:t\in\cI\}$ is a pie of $B$.
\end{quotation}
\QED$_{\ref{1.7}}$

\begin{conclusion}
\label{1.8}
Assume that $\sigma<\mu^\kappa=\mu<\lambda$ and there is a
$(\mu,\sigma)$-entangled linear order $\cI+\cJ$ such that for a uniform
ultrafilter $D$ on $\kappa$ the linear orderings $\cI^\kappa/D$,
$\cJ^\kappa/D$ contain isomorphic subsets of cardinality $\lambda>\mu$. Then 
\[\inc^+(\BA_{\inter}(\cI+\cJ))\leq\mu\quad\mbox{ and }\quad\inc((\BA_{\inter}
(\cI+\cJ))^\kappa/D)\geq\lambda\]
and even 
\[\inc^+((\BA_{\inter}(\cI+\cJ))^\kappa/D)>\lambda\]
(so $\inc((\BA_{\inter}(\cI+\cJ))^\kappa/D)>\inc(\BA_{\inter}(\cI+\cJ
))^\kappa/D$). \QED
\end{conclusion}

\noindent {\bf Remark:}\qquad See an example in \ref{3.2}(3).

\begin{definition}
\label{1.9}
We say that a linear order $\cI$ {\em has exact
$(\lambda,\mu,\kappa)$--density} if for every $\cJ\subseteq\cI$ of cardinality
$\ge\lambda$ we have ${\rm density}(J)\in[\kappa,\mu)$.\\
If $\mu=\kappa^+$ we may omit $\mu$; if $\lambda=|\cI|$ we may omit it. We may
also say $\cI$ {\em has exact density} $(\lambda,\mu,\kappa)$ or $(\lambda,
\mu,\kappa)$ {\em is an exact density of $\cI$} (and replace $(\lambda,\mu,
\kappa)$ by $(\lambda,\mu)$ or $(\mu,\kappa)$ or $\kappa$).
\end{definition}

\begin{definition}
\label{6.4}
\begin{enumerate}
\item  A linear order $\cI$ is {\em positively $\sigma$-entangled} if for each
$\vare(*)<1+\sigma$, $u\in \{\emptyset,\vare(*)\}$ and an indexed set
$\{t^\alpha_\vare: \alpha<|\cI|,\; \vare<\vare(*)\}\subseteq\cI$ such that
\[(\forall \alpha<\beta<|\cI|)(\forall\vare<\vare(*))(t^\alpha_\vare\neq
t^\beta_\vare)\]
there exist $\alpha<\beta<|\cI|$ such that $(\forall\vare<\vare(*))(\vare\in
u\ \Leftrightarrow t^\alpha_\vare<t^\beta_\vare)$.
\item Similarly we define when $\bar{\cI}=\langle\cI_i: i<i^*\rangle$
is positively $\sigma$-entangled and $\posens_\sigma(\lambda,\mu,\kappa)$,
$\posens_\sigma(\lambda,\kappa)$.
\end{enumerate}
\end{definition}

For more on entangledness in ultraproducts see section 6.

\section{Constructions for $\lambda=\lambda^{<\lambda}$}
In this section we will build entangled linear orders from instances of GCH.
Our main tool here is the construction principle presented in \cite{HLSh:162}
and developed in \cite{Sh:405}. The main point of the principle is that for
standard $\lambda^+$--semiuniform partial orders (see \ref{semiuni} below)
there are ``sufficiently generic'' filters $G$, provided
$\diamondsuit_\lambda$ holds (actually, a weaker assumption suffices). For the
precise definition of ``sufficiently generic'' we refer the reader to
\cite[Appendix]{Sh:405} (compare \cite[\S 1]{HLSh:162} too). Here we recall
the definition of standard $\lambda^+$--semiuniform partial orders, as it
lists the conditions we will have to check later.

\begin{definition}
\label{semiuni}
Let $\lambda$ be a regular cardinal.
\begin{enumerate}
\item A set $u\subseteq\lambda^+$ is {\em closed} if $0\in u$ and\ \ $\delta=
\sup(\delta\cap u)\quad\Rightarrow\quad\delta\in u$.
\item Let $(\cP,\leq)$ be a partial order such that
\[\cP\subseteq\lambda\times\{u\subseteq\lambda^+: |u|<\lambda^+\ \&\ u\mbox{
is closed}\}.\]
If $p=(\alpha,u)\in\cP$ then we write $\dom(p)=u$.\\
For an ordinal $\beta<\lambda^+$ we let $\cP_\beta=\{p\in\cP:\ \dom(p)
\subseteq\beta\}$.\\
We say that $(\cP,\leq)$ is {\em a standard $\lambda^+$--semiuniform partial
order} if the following conditions are satisfied:
\begin{description}
\item[(a)] If $p\leq q$ then $\dom(p)\subseteq\dom(q)$.
\item[(b)] If $p\in\cP$, $\alpha<\lambda^+$ is either a successor ordinal or
$\cf(\alpha)=\lambda$ then there exists $q\in\cP$ such that $q\leq p$ and
$\dom(q)=\dom(p)\cap \alpha$; moreover there is a unique maximal such $q$
which will be denoted $p\rest\alpha$.
\item[(c)] If $p=(\alpha,u)\in\cP$, $h:u\stackrel{\rm 1 - 1}{\longrightarrow}
v\subseteq\lambda^+$ is an order isomorphism onto $v$ such that $(\forall
\alpha\in u)(\cf(\alpha)=\lambda\ \Leftrightarrow\ \cf(h(\alpha))=\lambda)$
and $v$ is closed then $h[p]\stackrel{\rm def}{=} (\alpha,v)\in\cP$; moreover,
$q\leq p$ implies $h[q]\leq h[p]$.
\item[(d)] If $p,q\in\cP$, $\alpha<\lambda^+$ is either a successor ordinal or
$\cf(\alpha)=\lambda$ and $p\rest \alpha\leq q\in\cP_\alpha$ then there is
$r\in\cP$ such that $p,q\leq r$.
\item[(e)] If $\langle p_i:i<\delta\rangle\subseteq\cP$ is an increasing
sequence, $\delta<\lambda$ then there is $q\in\cP$ such that 
\[\dom(q)=\cl(\bigcup\limits_{i<\delta}\dom(p_i))\quad\mbox{and}\quad(\forall 
i<\delta)(p_i\leq q).\]
\item[(f)] Suppose that $\langle p_i:i<\delta\rangle\subseteq \cP_{\beta+1}$
is increasing, $\delta<\lambda$ and $\beta<\lambda^+$ has cofinality
$\lambda$.  Assume that $q\in\cP_\beta$ is such that $(\forall i<\delta)(
p_i\rest\beta\leq q)$. Then the family $\{p_i:i<\delta\}\cup\{q\}$ has an
upper bound $r$ such that $q\leq r\rest\beta$.
\item[(g)] Assume that $\langle\beta_i: i<\delta\rangle\subseteq\lambda^+$ is
strictly increasing, each $\beta_i$ is either a successor or has cofinality
$\lambda$, $\delta<\lambda$ is a limit ordinal. Suppose that $q\in\cP$ and
$(\forall i<\delta)(q\rest\beta_i\leq p_i\in\cP_{\beta_i})$ and $\langle
p_i:i<\delta\rangle\subseteq\cP$ is increasing. Then the family $\{p_i:i<
\delta\}\cup\{q\}$ has an upper bound $r\in\cP$ such that $(\forall i<\delta)
(p_i\leq r\rest\beta_i)$.
\item[(h)] Suppose that $\delta_1,\delta_2<\lambda$ are limit ordinals and
$\langle\beta_j: j<\delta_2\rangle\subseteq\lambda^+$ is a strictly increasing
sequence of ordinals, each $\beta_j$ either a successor or of cofinality
$\lambda$. Let 
\[\langle p_{i,j}: (i,j)\in (\delta_1+1)\times (\delta_2+1)\setminus\{(
\delta_1,\delta_2)\}\rangle\subseteq\cP\]
be such that
\[p_{i,j}\in\cP_{\beta_j},\quad i\leq i'\ \Rightarrow\ p_{i,j}\leq p_{i',j},
\quad j\leq j'\ \Rightarrow\ p_{i,j}\leq p_{i,j'}\rest \beta_j.\]
Then the family $\{p_{i,j}:(i,j)\in (\delta_1+1)\times (\delta_2+1)\setminus
\{(\delta_1,\delta_2)\}\}$ has an upper bound $r\in\cP$ such that
$(\forall j<\delta_2)(r\rest\beta_j=p_{\delta_1,j})$.
\end{description}
\end{enumerate}
\end{definition}

\begin{notation}
\label{2.1}
Let $\lambda,\mu$ be cardinals and $T$ be a tree.
\begin{enumerate}
\item  For an ordinal $\alpha$, the $\alpha$-th level of the tree $T$ is
denoted by $T_\alpha$;\quad for $x\in T$, $\lev(x)$ is the unique
$\alpha$ such that $x\in T_\alpha$.
\item  We say that the tree $T$ is normal if for each $y,z\in T$ we have:

if $(\forall x\in T)(x<_T y\equiv x <_T z)$ and $\lev(y)=\lev(z)$ is a
limit ordinal

then $y = z$.

\noindent Usually we assume that $T$ is normal.
\item We say that the tree $T$ is $\lambda^+$--Aronszajn if it has
$\lambda^+$ levels, each level is of size $\leq\lambda$, there is no
$\lambda^+$-branch in $T$, $T$ is normal, and
\[y\in T,\ \ \lev(y)<\beta<\lambda^+\quad\Rightarrow\quad (\exists z\in T)
[y<_T z\ \&\ \lev(z)=\beta].\]
\item  For ordinals $\zeta$ and $\alpha$ let $T_\alpha^{[\zeta]}$ be the
set of all sequences of length $\zeta$ with no repetition from
$T_\alpha$. We let $T^{[\zeta]}=\bigcup\limits_{\alpha}
T_\alpha^{[\zeta]}$, but we may identify $T^{[1]}$ and $T$ (and similarly for
$T^{\langle 1\rangle}$ below);
\item For a sequence $\bar{x}\in T^{[\zeta]}$, let $\lev(\bar{x})$ be
the unique $\alpha$ such that $\bar{x}\in T_\alpha^{[\zeta]}$.
\item For $\bar{x},\bar{y}\in T^{[\zeta]}$, let $\bar{x}<\bar{y}$ mean
$(\forall\vare<\zeta)(x_\vare<_T y_\vare)$; similarly for
$\bar{x}\leq\bar{y}$.
\item Let $\bar{x}\in T^{[\zeta]}_\alpha$. We define $T^{[\zeta]}_{\bar
x} \stackrel{\rm def}{=}\{\bar{y}\in T^{[\zeta]}: \bar{x}<_T\bar{y}\}$.
\item $T^{\langle\zeta\rangle}=\bigcup\{T^{[\zeta]}_{\alpha}:$ either
$\alpha$ is a successor ordinal or $\cf(\alpha)=\lambda\}$,\\
$T^{\langle\zeta\rangle}_{\bar x}=T^{[\zeta]}_{\bar x}\cap T^{\langle\zeta
\rangle}$.
\item For $x,y\in T$ let $x\wedge y$ be their maximal lower bound (in $T$,
exists when $T$ is normal).
\item For $x\in T$ and an ordinal $\alpha\leq\lev(x)$ let $x\rest\alpha$ is the
unique $y\leq_T x$ such that $\lev(y)=\alpha$.
\item For $\bar{x}=\langle x_\vare:\vare<\zeta\rangle\in T^{[\zeta]}$ and an
ordinal $\alpha\leq\lev(\bar{x})$ let $\bar{x}\rest\alpha=\langle x_\vare\rest
\alpha: \vare<\zeta\rangle$.
\item Let $H^1_\mu$ be the family of all functions $h$ with domains included
in $\bigcup\{{}^\zeta(\mu^+\times\mu^+): \zeta<\mu\}$ and such that for 
$\zeta<\mu$, $\bar{x}\in{}^\zeta(\mu^+\times\mu^+)$ we have:\quad
$h(\bar{x})\subseteq {}^\zeta(\mu^+)$ (if defined, and then) there are $\mu^+$
members of $h(\bar{x})$ with pairwise disjoint ranges.

If $h\in H^1_\mu$ is a partial function, $\zeta<\mu$, $\bar{x}\in {}^\zeta(
\mu^+)$ and $h(\bar{x})$ is not defined then $h(\bar{x})$ will mean
${}^\zeta(\mu^+)$.

We use mainly $h\in H_{\mu,*}^1$ where
\[H^1_{\mu,*}=\bigcup_{\zeta<\mu} H^1_{\mu,\zeta},\quad\quad
H^1_{\mu,\zeta}=\{h\in H^1_\mu: \dom(h)={}^\zeta(\mu^+\times\mu^+)\}.\]
\item Let $H^0_\mu$ be the set of all $h$ from $H^1_\mu$ such that
the value of $h(\langle(\alpha^0_\vare,\alpha^1_\vare): \vare<\zeta\rangle)$
does not depend on $\langle\alpha^0_\vare: \vare<\zeta\rangle$ (so we may
write $h(\langle\alpha^1_\vare:\vare<\zeta\rangle)$).
\item Let $H^3_\mu$ be the family of all functions $h$ with domain $\mu$
such that $h(\zeta)$ is a subset of ${}^\zeta((\mu^+)^3)$ with the following 
property
\begin{description}
\item[$(\boxtimes)$] for each $\langle\alpha^0_\vare:\vare<\zeta\rangle
\subseteq\mu^+$ and every $\beta<\mu^+$ there is $\langle\alpha^1_\vare:\vare
<\zeta\rangle\subseteq (\beta,\mu^+)$ with
\[(\forall\beta'<\mu^+)(\exists\langle\alpha^2_\vare:\vare<\zeta\rangle
\subseteq(\beta',\mu^+))(\langle(\alpha^0_\vare,\alpha^1_\vare,\alpha^2_\vare):
\vare<\zeta\rangle\in h(\zeta)).\]
\end{description}
\item $H^2_\mu$ is the collection of those $h\in H^3_\mu$ that the
truth value of ``$\langle (\alpha^0_\vare,\alpha^1_\vare,\alpha^2_\vare):
\vare<\zeta\rangle\in h(\zeta)$''does not depend on $\langle\alpha^0_\vare:
\vare\in\zeta\rangle$ (so we may write $\langle(\alpha^1_\vare,
\alpha^2_\vare): \vare<\mu\rangle\in h(\zeta)$).
\end{enumerate}
\end{notation}

\begin{theorem}
\label{2.2}
Suppose $\lambda=\mu^+=2^\mu$ and $\diamondsuit_\lambda$ (the second follows
e.g.~if $\mu\geq\beth_\omega$; see \cite[3.5(1)]{Sh:460}). Then:
\begin{enumerate}
\item There exists a dense linear order $\cI$ of cardinality $\lambda^+$ and
density $\lambda^+$ (really exact density $\lambda^+$, see \ref{1.9}) such
that: 
\begin{description}
\item[$(\ast)_1$] $\cI$ is hereditarily of the cellularity $\lambda^+$, i.e.
every interval in $\cI$ contains $\lambda^+$ pairwise disjoint subintervals,
and
\item[$(\ast)_2$] $\cI$ is $\mu$--locally entangled, i.e.~ if $\kappa<\mu$,
$(a_i,b_i)_{\cI}$ (for $i<\kappa$) are pairwise disjoint intervals then
the sequence $\langle\cI\rest(a_i,b_i): i<\kappa)\rangle$ is
$\kappa^+$--entangled\footnote{Note: for $\mu\in (2,\lambda)$, $\lambda=\cf(
\lambda)$ such that $(\forall\alpha<\lambda)(|\alpha|^{<\mu}<\lambda)$ and a
linear order $\cI$ of cardinality $\lambda$ we have:\quad $\cI$ is
$\mu$--entangled\ \ if and only if\ \ $\cI$ is $\mu$--locally entangled of
density $<\lambda$}.
\end{description}
\item Let $H^{1,*}_\mu\subseteq H^1_\mu$ and $H^{3,*}_\mu\subseteq
H^3_\mu$ have cardinality $\leq\lambda$. There is a $\lambda^+$--Aronszajn
tree $T\subseteq {}^{\lambda^+>}\lambda$ in which each node has
$\lambda$ immediate successors and there are two functions $c,\bar{d}$ such
that: 
\begin{description}
\item[(a)] $c$ is a function from $T$ to $\lambda$,
\item[(b)] for every $\bar{x}\in T^{[\mu]}$ and a function $h\in H^{1,*}_\mu
\cup H^{3,*}_\mu$ we have $d_{\bar{x},h}:T^{[\mu]}_{\bar{x}}\longrightarrow
\lambda$ such that\quad if $\bar{y},\bar{z}\in T^{\langle\mu\rangle}_{
\bar{x}}$ are distinct, $d_{\bar{x},h}(\bar{y})=d_{\bar{x},h}(\bar{z})$ then
for some $\bar{t}\in T^{[\mu]}_{\bar x}$ we have 
\begin{description}
\item[($\alpha$)] $t_i=y_i\wedge z_i$, and the values of $\lev(t_i)$ do not
depend on $i$,
\item[($\beta$)] $\lev(\bar{t})<\lev(\bar{z})$, $\lev(\bar{t})<
\lev(\bar{y})$,
\item[($\gamma$)] if $d_{\bar{x},h}(\bar{y})<d_{\bar{x},h}(\bar{t})$ then
$(\forall\vare<\mu)(\exists^{\mu}i<\mu)(c(t_i)=\vare)$,
\item[($\delta$)] if $\zeta<\mu$ and $h\in H^{1,*}_\mu$ then for $\mu$
ordinals $i<\mu$ divisible by $\zeta$ we have\\
{\rm (i)}\qquad either $\langle z_{i+\vare}(\lev(\bar{t})):\vare<\zeta
\rangle$ belongs to 
\[h(\langle c(t_{i+\vare}),y_{i+\vare}(\lev(\bar{t})):\vare<\zeta
\rangle)\]
{\rm (ii)}\qquad or $\langle y_{i+\vare}(\lev(\bar{t})):\vare<\zeta
\rangle$ belongs to 
\[h(\langle c(t_{i+\vare}),z_{i+\vare}(\lev(\bar{t})):
\vare<\zeta\rangle),\]
\item[($\vare$)] if $\zeta<\mu$ and $h\in H^{3,*}_\mu$ and $d_{\bar{x},h}(
\bar{y})<d_{\bar{x},h}(\bar{t})$ then for $\mu$ ordinals $i<\mu$ divisible
by $\zeta$ we have\\
{\rm (i)}\qquad either $\langle(c(t_{i+\vare}),y_{i+\vare}(\lev(\bar{t})),
z_{i+\vare}(\lev(\bar{t}))): \vare<\zeta\rangle$ belongs to $h(\zeta)$\\
{\rm (ii)}\qquad or $\langle(c(t_{i+\vare}),z_{i+\vare}(\lev(\bar{t})),
y_{i+\vare}(\lev(\bar{t}))): \vare<\zeta\rangle$ belongs to $h(\zeta)$.
\end{description}
\end{description}
\end{enumerate}
\end{theorem}

\noindent{\bf Explanation:}\quad Some points in \ref{2.2}(2) may look
unnatural. 
\begin{enumerate} 
\item Why $\bar{y},\bar{z}\in T^{\langle\mu\rangle}_{\bar{x}}$ and not $\dom(
d_{\bar{x},h})$? As in proving amalgamation we should compare $\langle x^{
\beta}_i:i<\mu\rangle$, $\langle x^\alpha_i:i<\mu\rangle$; necessary when
$\alpha=\sup(w^q)$. However, working a little bit harder we may wave this.
\item Why e.g~in clause (b)($\vare$) we demand $d_{\bar{x},h}(\bar{y})<
d_{\bar{x},h}(\bar{t})$? Otherwise we will not be able to prove the 
density of 
\[{\cal D}_{\bar{x},h,\bar{y}}=:\{p\in\cA\cP: \bar{y}\in\dom(d_{\bar{x},
h})\quad (\mbox{ or }\quad \neg\bar{x}<\bar{y})\}.\]
\item Why do we have clauses (b)($\delta$) and (b)($\vare$)? For the
application here (b)($\delta$) suffices, if this is enough for the reader
then clause (J) in the definition of $\cA\cP$ may be omitted. But they
both look ``local maximal''. 
\end{enumerate}

\noindent{\sc Proof of \ref{2.2}(1)}\hspace{0.2in} We will use 2).\\
Let $T \subseteq {}^{\lambda^+>}\lambda$ and $c,\bar{d}$ be as there and let 
all the functions $h_{\kappa,u}\in H^0_\mu$ defined in the continuation 
of the proof of \ref{2.2}(1) below be in $H^{1,*}_\mu$. We may assume 
that if $x\in T$ and $\alpha<\lambda$ then $x^\frown\!\langle\alpha\rangle
\in T$. Let $<^\otimes$ be a linear order on $\lambda$ such that 
$(\lambda,<^\otimes)$ has neither first nor last element and is
$\lambda$-dense  (i.e.~if $\alpha_i<^\otimes\beta_j$ for $i<i_0<\lambda$,
$j<j_0<\lambda$ then for some $\gamma$, $\alpha_i\leq^\otimes\gamma
\leq^\otimes\beta_j$). We define the order $<_{\cI}$ on $\cI=T^{\langle 1
\rangle}=\{x\in T: \lev(x)$ is a successor or of cofinality $\lambda\}$:
\begin{quotation}
\noindent $y <_{\cI} z$\qquad  if and only if\\
either $z = (y \wedge z) \vartriangleleft y$ or $y(\alpha))<^\otimes
z(\alpha)$, where $\alpha=\lev(y \wedge z)$.
\end{quotation}
Clearly $(\cI,<_{\cI})$ is a dense linear order of the density $\lambda^+$
and the size $\lambda^+$. To show that $\cI$ has exact density $\lambda^+$
(i.e.~its exact density is $(\lambda^+,\lambda^{++},\lambda^+)$) assume
that $\cJ\subseteq\cI$, $|\cJ|=\lambda^+$. We want to show that $\cJ$ has
density $\lambda^+$. Suppose that $\cJ_0\subseteq\cJ$, $|\cJ_0|\leq\lambda$.
Then for some $\alpha(*)<\lambda^+$, $\cJ_0\subseteq\bigcup\limits_{\alpha<
\alpha(*)} T_\alpha$ and we may find distinct $x,y\in\cJ\setminus
\bigcup\limits_{\alpha\leq\alpha(*)}T_\alpha$ such that $x\rest\alpha(*)=
y\rest\alpha(*)$. Then $x,y$ show that $\cJ_0$ is not dense in $\cJ$.
\medskip

Now we are proving that $\cI$ satisfies $(*)_2$.
\smallskip

\noindent Suppose that $\kappa<\mu$ and $(a_i,b_i)$ are disjoint intervals in
$\cI$ (for $i<\kappa)$. Suppose that $\bar{y}_{\xi}=\langle y^\xi_i:
i<\kappa\rangle$ (for $\xi<\lambda^+)$ are such that $a_i<_{\cI} y^\xi_i
<_{\cI} b_i$ and $y^\xi_i$'s are pairwise distinct for $\xi<\lambda^+$.
Let $u\subseteq\kappa$. Take $\alpha(*)<\lambda^+$ such that $(\forall i<
\kappa)(\lev(a_i),\lev(b_i)<\alpha(*))$. As $y^\xi_i$'s are pairwise distinct
we may assume that
\[(\forall \xi<\lambda^+)(\forall i<\kappa)(\alpha(*)<\lev(y^\xi_i)\;\;\mbox{
and }\;\;\xi\leq\lev(y^\xi_i)).\]
Note that if $i<j<\kappa$, and $\xi,\zeta<\lambda^+$ then $y^\xi_i\rest
\alpha(*)\neq y^\zeta_j\rest \alpha(*)$. Now the following claim (of self
interest) is applicable to $\langle\langle y^\xi_i:i<\kappa\rangle: \xi\in
[\alpha(*),\lambda)\rangle$ and as we shall see later this finishes the proof
of \ref{2.2}(1) shortly.
\begin{claim}
\label{2.2B}
Assume (for the objects constructed in \ref{2.2}(2)):
\begin{description}
\item[(a)] $\kappa<\mu$, 
\item[(b)] for each $\xi<\lambda^+$ we have a sequence $\bar{y}^\xi=\langle
y^\xi_i:i<\kappa\rangle$ such that $y^\xi_i\in T$ and either
\begin{description}
\item[($\alpha$)] $(\xi^1,i_1)\neq (\xi^2,i_2)\quad \Rightarrow\quad
y^{\xi^1}_{i_1}\neq y^{\xi^2}_{i_2}$\qquad or
\item[($\beta$)] $\lev(y^\xi_i)\geq\xi$,
\end{description}
\item[(c)] $h\in H^{1,*}_\mu\cup H^{3,*}_\mu$,
\item[(d)] for some $\alpha(*)<\lambda^+$,
\[\xi,\zeta<\lambda^+,\ i<j<\kappa\quad \Rightarrow\quad y^\xi_i\rest
\alpha(*)\neq y^\zeta_j\rest \alpha(*).\]
\end{description}
Then we can find $\xi_1<\xi_2<\lambda^+$ such that clause (b)($\delta$)(i) (or
(b)($\vare$)(i), respectively) of \ref{2.2}(2) holds with $(\bar{y}^{\xi_1},
\bar{y}^{\xi_2})$ standing for $(\bar{y},\bar{z})$, i.e.:
\begin{description}
\item[($\delta$)] if $h\in H^{1,*}_\mu$ then\\
(i)\quad $\langle y^{\xi_2}_{\vare}(\lev(\bar{t})):\vare<\kappa\rangle\in h(
\langle c(t_{\vare}),y^{\xi_1}_{\vare}(\lev(\bar{t})):\vare<\kappa\rangle$,
\item[($\vare$)]  if $h\in H^{3,*}_\mu$ then\\
(i)\quad $\langle(c(t_{\vare}), y^{\xi_1}_{\vare}(\lev(\bar{t})),
y^{\xi_2}_{\vare}(\lev(\bar{t}))): \vare<\kappa\rangle$ belongs to
$h(\kappa)$,
\end{description}
where $\bar{t}=\langle t_j:j<\kappa\rangle$ and $t_j=y^{\xi_1}_j\wedge
y^{\xi_2}_j$ etc.
\end{claim}

\noindent{\em Proof of the claim:}\qquad First note that, by easy thining (as
either (b)($\alpha$) or (b)($\beta$) holds; remember clause (e)) we can assume
(b)($\alpha$)\ \&\ (b)($\beta$). As $\lambda=\lambda^\kappa$ we may assume
that $\langle y^\delta_i\rest \alpha(*): i<\kappa\rangle$ is the same for all
$\delta\in \lambda^+$. Let 
\[Z=\{\alpha\in [\alpha(*),\lambda^+): (\exists Y\in [T_\alpha]^{<\lambda})
(|\{\xi<\lambda^+:(\forall i<\kappa)(y^\xi_i\rest \alpha\notin Y)\}|\leq
\lambda)\}.\]
First we are going to show that $Z\neq[\alpha(*),\lambda^+)$. If not then for
each $\alpha\in[\alpha(*),\lambda^+)$ we find a set $Y_\alpha\in
[T_\alpha]^{<\lambda}$ and an ordinal $\xi(\alpha)<\lambda^+$ such that
\[\{\xi<\lambda^+: (\forall i<\kappa)(y^\xi_i\rest\alpha\notin Y)\}\subseteq
\xi(\alpha).\]
For $\alpha\in [\alpha(*),\lambda^+)$, choose $\bar x^\alpha\in
T^{[\mu]}_\alpha$ such that
\begin{description}
\item[(i)]  $(\forall i<\kappa)(y^\alpha_i\rest\alpha=x^\alpha_i)$,
\item[(ii)] $Y_\alpha\subseteq\{x^\alpha_i: i<\mu\}$.
\end{description}
For each $\delta\in [\alpha(*),\lambda^+)$ with $\cf(\delta)=\lambda$ we can
find $\gamma_\delta<\delta$ such that $\bar{x}^*_\delta=\langle x^\delta_i
\rest\gamma_\delta: i<\mu\rangle$ is with no repetition (recall that
$x^\delta_i\in T_\delta$ are pairwise distinct, $i<\mu<\lambda$ and the tree
$T$ is normal). By Fodor's lemma, for some $\gamma^*$ the set
\[S_0\stackrel{\rm def}{=}\{\delta\in [\alpha(*),\lambda^+):
\cf(\delta)=\lambda\ \&\ \gamma_\delta=\gamma^*\}\]
is stationary. For $\delta\in S_0$ there are at most $\lambda^\mu=\lambda$
possibilities for $\langle x^\delta_i\rest\gamma^*:i<\mu\rangle$ and hence
for some $\bar{x}^*=\langle x^*_i:i<\mu\rangle\in T^{[\mu]}_{\gamma^*}$ the set
\[S_1\stackrel{\rm def}{=}\{\delta\in S_0:\bar{x}^*=\langle x^\delta_i\rest
\gamma^*: i<\mu\rangle\}\]
is stationary. Hence the set
\[S_2\stackrel{\rm def}{=}\{\delta\in [\alpha(*),\lambda^+)\cap S_1:
(\forall i<\mu)(x^\delta_i\rest\gamma^*=x^*_i)\ \&\ (\forall\alpha<\delta)
(\xi(\alpha)<\delta)\}\]
is stationary. Look at $d_{\bar{x}^*,h}$ (really any $h'\in H^{\ell,*}_\mu$
will do here). Note that $\bar{x}^\delta\in\dom(d_{\bar{x}^*,h})$ for $\delta
\in S_2$ (remember $\cf(\delta)=\lambda$), and therefore we find $\delta_1,
\delta_2\in S_2$, $\delta_1<\delta_2$ such that $d_{\bar{x}^*,h}(\bar{x}^{
\delta_1})=d_{\bar{x}^*,h}(\bar{x}^{\delta_2})$. As $\xi(\delta_1)<\delta_2$
we can find $i<\kappa$ such that $y_i^{\delta_2}\rest\delta_1\in Y_{
\delta_1}$. Thus for some $j<\mu$ necessarily $y_i^{\delta_2}\rest\delta_1
=x^{\delta_1}_j$ and hence $x^{\delta_2}_i\rest\delta_1=x_j^{\delta_1}$ and
consequently $x^{\delta_2}_i\rest\gamma^*=x^{\delta_1}_j\rest\gamma^*$. This
implies $i=j$ (as $\delta_1,\delta_2\in S_2\subseteq S_0$ and hence
$\gamma_{\delta_1}=\gamma^*=\gamma_{\delta_2}$ and now apply the definition of
$\gamma_{\delta_1},\gamma_{\delta_2}, S_2$) and thus $x^{\delta_1}_i
\vartriangleleft x^{\delta_2}_i$, what contradicts clause $(\beta)$ of
\ref{2.2}(2b). So we have proved $Z\neq [\alpha(*),\lambda^+)$, but by its
definition $Z$ is an initial segment of $[\alpha(*),\lambda^+)$. Hence for
some $\beta(*)\in [\alpha(*),\lambda^+)$ we have $Z=[\alpha(*),\beta(*))$.
\smallskip

\noindent Let $\beta\in [\beta(*),\lambda^+)$ be a successor ordinal.

\noindent By induction on $\vare<\mu$ choose\footnote{We could have done it for
$\vare<\lambda$, but no need here.} pairwise disjoint $\bar{x}^\vare=\langle
x^\vare_i: i<\kappa\rangle\in T^{[\kappa]}_\beta$ such that the sets
\[Z_\vare\stackrel{\rm def}{=}\{\xi\in [\beta,\lambda^+): (\forall
i<\kappa)(y^\xi_i\rest\beta=x^\vare_i)\}\]
are of the size $\lambda^+$. Suppose we have defined $\bar{x}^\vare$ for
$\vare<\vare'$. Let $Y=\{x^\vare_i: i<\kappa,\; \vare<\vare'\}$. Then $Y\in
[T_\beta]^{<\lambda}$ and by the choice of $\beta$ the set $\{\xi<\lambda^+:
(\forall i<\kappa)(y^\xi_i\rest\beta\notin Y)\}$ is of the size $\lambda^+$.
As $\lambda^\kappa=\lambda$ we can find $\bar{x}^{\vare'}\in
T^{[\kappa]}_\beta$ such that
\[|\{\xi<\lambda^+:(\forall i<\kappa)(y^\xi_i\rest\beta=x^{\vare'}_i\notin
Y)\}|=\lambda^+.\]
Now for $\vare<\mu$ and $i<\kappa$ let $x_{\kappa{\cdot}\vare+i}=x^\vare_i$
and let $\bar{x}=\langle x_j: j<\mu\rangle$. Thus $\bar{x}\in
T^{[\mu]}_\beta$. For each $\alpha\in [\beta,\lambda^+)$ choose $\langle
\zeta_{\alpha,\vare}:\vare<\mu\rangle$ such that $\zeta_{\alpha,\vare}
\in Z_\vare$ and $\alpha_1<\alpha_2$ implies $\alpha_1<\zeta_{\alpha_1,
\vare}<\zeta_{\alpha_2,\vare}$. For $\vare<\mu$, $i<\kappa$ and $\alpha\in
[\beta,\lambda^+)$ put $z^\alpha_{\kappa{\cdot}\vare+i}=y^{\zeta_{\alpha,
\vare}}_i\rest\alpha$. Then $\bar{z}^\alpha=\langle z^\alpha_j:j<\mu\rangle
\in T^{[\mu]}_\alpha$, $\bar{x}<\bar{z}^\alpha$. Consider the function
$d_{\bar{x},h}$. Let $S_3\stackrel{\rm def}{=}\{\delta\in (\beta,\lambda^+): 
\cf(\delta)=\lambda\}$ and note that for each $\alpha\in S_3$ the restriction
$d_{\bar{x},h}\rest\{\bar{z}^\alpha\rest\gamma:\beta<\gamma\leq\alpha\}$ is a
one-to-one function (since $\cf(\alpha)=\lambda$ and we have clause ($\beta$)
of \ref{2.2}(2b)). Consequently, for $\alpha\in S_3$, we find $\delta(\alpha)$
such that 
\[\delta(\alpha)\leq\delta<\alpha\quad \Rightarrow\quad d_{\bar{x},h}(
\bar{z}^{\alpha}\rest\delta)>d_{\bar{x},h}(\bar{z}^\alpha);\qquad\qquad
\delta(\alpha)>\beta.\]
Next, applying Fodor Lemma, we find a stationary set $S_4'\subseteq S_3$
and $\delta^*>\beta$ such that $d_{\bar{x},h}(\bar{z}^{\alpha}\rest\delta)
>d_{\bar{x},h}(\bar{z}^{\alpha})$ for $\alpha\in S_4'$, $\delta\in
(\delta^*,\alpha)$ and $\bar{z}^{\alpha_1}\rest\delta^*=\bar{z}^{\alpha_2}
\rest\delta^*$ for all $\alpha_1,\alpha_2\in S_4'$. Let
\[S_4\stackrel{\rm def}{=}\{\delta\in S_4': (\forall\alpha<\delta)(|\{\beta
\in S_4':\bar{z}^\beta\rest\alpha=\bar{z}^\delta\rest\alpha\}|=
\lambda^+)\}.\]
If $S_4$ is not stationary then for $\delta\in S_4'\setminus S_4$ choose
$\alpha_\delta<\delta$ contradicting the demand in the definition of
$S_4$. For some stationary $S^*\subseteq S_4'\setminus S_4$ we have
$\alpha_\delta=\alpha^*$ for $\delta\in S^*$ and we get an easy contradiction.

\noindent As $\rang(d_{\bar{x},h})\subseteq\lambda$, for some stationary
$S_5\subseteq S_4$, $d_{\bar{x},h}\rest\{\bar{z}^\alpha:\alpha\in S_5\}$
is constant. Choose $\alpha_1\neq\alpha_2$ from $S_5$. By clauses
($\alpha$), ($\beta$) of \ref{2.2}(2b) we get that if $t_j=z_j^{\alpha_1}
\wedge z^{\alpha_2}_j$ (for $j<\mu$) then $\bar{t}\in T^{[\mu]}_{\bar{x}}$
and $\lev(\bar{x})\leq\lev(\bar{t})<\lev(\bar{z}^{\alpha_1})$, $\lev(\bar{t})
<\lev(\bar{z}^{\alpha_2})$. Moreover by the definition of $S_4'$ we have
$\lev(\bar{t})\geq\delta^*$ and $d_{\bar{x},h}(\bar{t})>d_{\bar{x},h}
(\bar{z}^{\alpha_1})$. Hence by the clause $(\delta)$ of \ref{2.2}(2b) we
find $i<\mu$ divisible by $\kappa$ such that either (i) or (ii) of clause
($\delta$) of \ref{2.2}(2b) holds with $(\bar{y},\bar{z})$ there standing
for $(\bar{z}^{\alpha_1},\bar{z}^{\alpha_2})$ here. By the symmetry wlog (i)
of \ref{2.2}(2b)($\delta$) holds. Let $\lev(\bar{t})<\beta_\ell<\alpha_\ell$
(for $\ell=1,2$), hence by the definition of $S_4$ (and as $\alpha_2\in
S_5\subseteq S_4$) and by the character of the requirement on $\alpha_1,
\alpha_2$ wlog $\alpha_1<\alpha_2$, so we are done. \hfill
$\square_{\ref{2.2B}}$
\medskip

\noindent{\sc Continuation of the proof of \ref{2.2}(1):}\hspace{0.2in}
Remember we have $\kappa<\mu$, $(a_i,b_i)$ (for $i<\mu$), $u\subseteq\kappa$
and we have $\langle\bar{y}^\xi: \xi\in [\alpha(*),\lambda^+)\rangle$,
$\bar{y}^\xi=\langle y^\xi_i:i<\kappa\rangle$. Define $h=h_{\kappa,u}\in
H^0_\mu$ (and assume that for each $\kappa<\mu$ and $u\subseteq\kappa$ we have
$h_{\kappa,u}\in H^{1,*}_\mu$):
\[h(\langle\beta^1_i:i<\kappa\rangle)=\{\langle\beta^2_i:i<\kappa\rangle\in
{}^\kappa\lambda:\ (\forall i<\kappa)(\beta^1_i<^\otimes\beta^2_i\
\Leftrightarrow\ i\in u)\}.\]
By the choice of $<^\otimes$ it is easy to check that $h\in H^0_\mu$. So by
Claim \ref{2.2B} for $\kappa$, $h$, $\bar{y}^\xi$, there are $\alpha^1<
\alpha^2$ as there and we are done.
\medskip

\noindent We still have to prove $(*)_1$.\\
Suppose that $a,b\in\cI$, $a<_{\cI}b$. By the definition of the order
there is $t\in T$ such that $t\trianglelefteq s\in T\ \Rightarrow\ s\in
(a,b)_{\cI}$. As the tree $T$ is $\lambda^+$--Aronszajn we find $\bar{x}\in
T^{[\mu]}_\alpha$ (for some $\alpha\in (\lev(t),\lambda^+)$) such that
$(\forall j<\mu)(t\trianglelefteq x_j)$. Next for every $\beta\in
(\alpha,\lambda^+)$ we can choose $\bar{y}_\beta\in T^{[\mu]}_\beta$ such
that $\bar{x}<\bar{y}_\beta$. Take any $h\in H^{1,*}_\mu\cup H^{3,*}_\mu$.
For some unbounded $S\subseteq(\alpha,\lambda^+)$ the sequence $\langle
d_{\bar{x},h}(\bar{y}_\beta):\beta\in S\rangle$ is constant. Consequently,
elements $y^\beta_0$ for $\beta\in S$ are pairwise
$\vartriangleleft$--incomparable (in the tree $T$). Hence $\{\{y\in T:
y^\beta_0\trianglelefteq y\}: \beta\in S\}$ is a family of pairwise
disjoint convex subsets of $(a,b)_{\cI}$ (each with $\lambda^+$
elements), so we have finished. \hfill$\square_{2.2(1)}$
\medskip

\noindent{\sc Proof of \ref{2.2}(2):}\hspace{0.2in} We want to apply
\cite[Appendix]{Sh:405}. For this we have to define a set $\cA\cP$ of 
approximations and check the conditions of \ref{semiuni}.
\begin{sdef}
\label{2.2C}
The set $\cA\cP$ of approximations consists of all tuples $p=\langle t,w,
\leq,D,\bar{d},\bar{e},f,c\rangle$ (we may write $t^p$, $w^p$ etc) such
that
\begin{description}
\item[(A)] $t$ is a subset of $\lambda^+$ of cardinality $<\lambda$, $t\cap
[0,\lambda)=\{0\}$,
\item[(B)] $w=\{\alpha<\lambda^+:[\lambda{\cdot}\alpha,\lambda{\cdot}
(\alpha+1))\cap t\neq\emptyset\}$ is such that 
\[(\forall \alpha<\lambda^+)(\alpha\in w\quad\Leftrightarrow\quad\alpha+1
\in w),\qquad\mbox{ and}\]
the set $\{\alpha<\lambda^+:\omega\cdot\alpha\in w\}$ is closed, 
\item[(C)] $\leq=\leq_t$ is a partial order on $t$ such that $t=t^p=(t,\le)$
is a normal tree (so $x\wedge y$ is well defined) and for each
$\alpha\in w$ the set $t\cap [\lambda{\cdot}\alpha,\lambda{\cdot}(\alpha+1))$
is the $\otp(w\cap\alpha)$-th level of $t$ (but possibly $\alpha<\beta$ are
in $w$, and for some $x$ in the $\alpha$-th level of $t$ there is no $y$
in the $\beta$-th level, $x\leq_t y$)

\noindent{\em  [obviously, the intention is: $t$ approximates $T$; we may use
$t$ for $(t,\le_t)$]},
\item[(D)] $D$ is a set of $<\lambda$ pairs $(\bar{x},h)$ such that
$\bar{x}\in t^{[\mu]}$ and $h\in H^{1,*}_\mu\cup H^{3,*}_\mu$,
\item[(E)] $\bar{d}=\langle d_{\bar x,h}: (\bar{x},h)\in D\rangle$, each
$d_{\bar{x},h}$ is a partial function from $t^{[\mu]}_{\bar{x}}$ to $\lambda$
with domain of cardinality $<\lambda$ such that: 
\[[\bar{y}\in t^{[\mu]}_{\bar{x}}\ \ \&\ \ \bar{x}<_t\bar{y}<_t\bar{z}\ \ \&\
\ \bar{z}\in\dom(d_{\bar{x},h})]\quad\Rightarrow\quad\bar{y}\in\dom(d_{
\bar{x},h}),\]
\item[(F)] $\bar{e}=\langle e_{\bar{x},h}: (\bar{x},h)\in D\rangle$,
each $e_{\bar{x},h}$ is a partial function from $t^{[\mu]}_{\bar{x}}\times
\lambda$ to $\{0,1,2\}$ of size $<\lambda$ and such that:
\begin{description}
\item[(i)] $(\bar{y},\gamma)\in\dom(e_{\bar{x},h})\quad\Rightarrow\quad \bar{y}
\in\dom(d_{\bar{x},h})$,\qquad and 
\[\hspace{-0.5cm}\dom(e_{\bar{x},h})\supseteq\{(\bar{y},\gamma)\!:\bar{y}\!
\in\!\dom(d_{\bar{x},h})\ \&\ (\exists\bar{z}\!\in\!\dom(d_{\bar{x},h}))(
\gamma\leq d_{\bar{x},h}(\bar{z}))\},\]
\item[(ii)] if $\bar{y}\in t^{[\mu]}_{\bar{x}}$, $\bar{x}<_t\bar{y}<_t\bar{z}
\in\dom(d_{\bar{x},h})$ and $e_{\bar{x},h}(\bar{z},\beta)$ is defined then
$e_{\bar{x},h}(\bar{y},\beta)$ is defined and $e_{\bar{x},h}(\bar{y},\beta)
\leq e_{\bar{x},h}(\bar{z},\beta)$ and if $\bar{y}\in t^{\langle\mu\rangle}_{
\bar{x}}$ then at most one of them is $1$,

[here, we interpret $t^{\langle\mu\rangle}_{\bar{x}}$ as the set of those
$\bar{y}\in t^{[\mu]}_{\bar{x}}$ that $\lev(\bar{y})$ is either a successor
ordinal or is $\otp(w\cap\alpha)$ for some $\alpha\in w$ such that
$\cf(\alpha)=\lambda$],
\item[(iii)] $d_{\bar{x},h}(\bar{y})=\alpha\quad\Leftrightarrow\quad
e_{\bar{x},h}(\bar{y},\alpha)=1$,
\end{description}
\noindent{\em [the intention: $e_{\bar{x},h}$ is not explicitly present in
\ref{2.2}(2b), but $e_{\bar{x},h}(\bar{y},\gamma)=\ell$ will mean that:\quad
if $\ell=0$ then for some $\bar{z}$ we have $\bar{y}<\bar{z}$ and
$d_{\bar{x},h}(\bar{z})=\gamma$;\ if $\ell=1$ then $d_{\bar{x},h}(\bar{y})=
\gamma$ \ and if $\ell=2$ then none of these]},
\item[(G)] $f$ is a function from $t^+=\{\alpha\in t:\ \alpha$ is of a
successor level in $t\}$ to $\lambda$ such that:
\begin{quotation}
\noindent if $\gamma\neq\beta$ are immediate successors (in $t$) of some 
$\alpha$\\
then $f(\beta)\neq f(\gamma)$,
\end{quotation}
\noindent{\em  [the intention is that if $\alpha$ represents $\eta\in
T_{i+1}$, then $f(\alpha)=\eta(i)$]},
\item[(H)] $c$ is a function from $t$ to $\lambda$,
\item[(I)] if $(\bar{x},h)\in D$, $e_{\bar{x},h}(\bar{y},\alpha),e_{\bar{x},
h}(\bar{z},\alpha)\le 1$, $\neg[\bar{y}\leq_t\bar{z}]$, $\neg [\bar{z}\leq_t
\bar{y}]$

then clauses ($\alpha$)---($\vare$) of \ref{2.2}(2b) hold (with $\langle
f(y_\vare\rest (i+1)),f(z_\vare\rest (i+1))\rangle$ replacing $\langle
y_\vare(i),z_\vare(i)\rangle)$ and $\alpha$ in place of $d_{\bar{x},h}
(\bar{y})=d_{\bar{x},h}(\bar{z})$),
\item[(J)] if $\alpha\in w^p$, $\bar{x}<\bar{y}^0<\bar{y}^1$, all in
$(t^p)^{[\mu]}$, $(\bar{x},h)\in D$, $h\in H^{3,*}_\mu$ and $\bar{y}^0,
\bar{y}^1\in\dom(d_{\bar{x},h})$, $\lev(\bar{y}^0)=\alpha$, $\lev(\bar{y}^1)
=\alpha+1$, $e_{\bar{x},h}(\bar{y}^1,\gamma)\leq 1$ and $\gamma<d_{\bar{x},h}
(\bar{y}^0)$

then in looking at $\bar{y}^0,\bar{y}^1$ as candidates for
$\bar{t},\bar{y}$ (or $\bar{t},\bar{z})$ in clause (I) (i.e.~in clauses
($\alpha$)---($\vare$) from \ref{2.2}(2b)) in ($\vare$) there, for each
$\zeta<\mu$, for $\mu$ ordinals $i<\mu$ divisible by $\zeta$, the values we
have i.e.  $c^p(y^0_{i+\vare})$, $f^p(y^1_{i+\vare})$ for $\vare<\zeta$,
are compatible with the demand (i) there, i.e. for every $\beta<\lambda$
there are $\alpha^2_\vare\in (\beta,\lambda)$ for $\vare<\zeta$ such that
\[\langle(c^p(y_{i+\vare}^0),f^p(y^1_{i+\vare}),\alpha^2_\vare):\vare<
\zeta\rangle\in h(\zeta),\]
\item[(K)] if $\bar{x}<\bar{y}$ are in $t^{[\mu]}$, $(\bar{x},h)\in D$,
$\bar{y}\in\dom(d_{\bar{x},h})$, $\alpha<d_{\bar{x},h}(\bar{y})$ and
$e_{\bar{x},h}(\bar{y},\alpha)=0$ then $(\forall\vare<\mu)(\exists^\mu
i<\mu)(c(y_i)=\vare)$ (i.e.~looking at $\bar{y}$ as a candidate for $\bar{t}$
in \ref{2.2}(2b)($\gamma$) the values we have are compatible with the demand
there). 
\end{description}
The set $\cA\cP$ of approximations is equipped with the natural partial
order.
\end{sdef}

We will want to apply the machinery of \cite[Appendix]{Sh:405} to the partial
order $(\cA\cP,\leq)$. For this we have to represent it as a standard
$\lambda^+$--semiuniform partial order. In representing it as a partial order
on $\lambda\times[\lambda^+]^{\textstyle{<}\lambda}$ we define the set of terms
such that:
\begin{description}
\item[(a)] $\{\tau(u): \tau$ a term with $\otp(u)$ places$\}=\{p\in\cA\cP: 
\{\alpha<\lambda^+: \omega{\cdot}\alpha\in w^p\}=u\}$,\qquad for a closed set
$u\in [\lambda^+]^{<\lambda}$,
\item[(b)] if $p_\ell=\tau(u_\ell)$ for $\ell=1,2$ then $\otp(t^{p_1},
\leq)=\otp(t^{p_2},\leq)$ and the one-to-one order preserving mapping $g$ 
from $t^{p_1}$ onto $t^{p_2}$ maps $p_1$ to $p_2$ (i.e.~$\alpha\leq^{p_1}
\beta\ \Leftrightarrow\ g(\alpha)\leq^{p_2} g(\beta)$, etc). 
\end{description}
Note that for $p\in\cA\cP$, its domain $\dom(p)$ (in the sense of
\ref{semiuni}) is $\{\alpha:\omega\cdot\alpha\in w^p\}$. Hence, $\cA\cP_\alpha
=\{p\in\cA\cP:w^p\subseteq\omega\cdot\alpha\}$. 

Now we have to check that (with this representation) $\cA\cP$ satisfies the
demands \ref{semiuni}(2)(a)--(h). Clauses (a) and (c) there should be clear.  
 
To deal with the clause (b) of \ref{semiuni}(2), for an approximation
$p\in\cA\cP$ and $\alpha<\lambda^+$ such that either $\alpha$ is a successor
ordinal or $\cf(\alpha)=\lambda$, we define $q=p\rest \omega\cdot\alpha$ by:
\begin{itemize}
\item $t^q=t^p\cap\lambda{\cdot}(\omega\cdot\alpha)$, $w^q=w^p\cap \omega\cdot
\alpha$, $\leq^q=\leq^p\rest t^q$, 
\item $D^q=\{(\bar{x},h)\in D^p: \bar{x}\subseteq t^q\}$,
\item if $(\bar{x},h)\in D^q$ then $d^q_{\bar{x},h}=d^p_{\bar{x},h}\rest
(t^q)_{\bar{x}}^{[\mu]}$ and $e^q_{\bar{x},h}=e^p_{\bar{x},h}\rest ((t^q)_{
\bar{x}}^{[\mu]}\times\lambda)$,
\item $f^q=f^p\rest(t^q\cap\dom(f^q))$, $c^q=c^p\rest t^q$.
\end{itemize}

\begin{observation}
\label{2.2E}
If $p\in\cA\cP$ is an approximation and $\alpha<\lambda^+$ is either a
successor or of cofinality $\lambda$ then $p\rest\omega\alpha\in\cA\cP$
is a unique maximal approximation such that $p\rest\omega\cdot\alpha\leq p$
and $w^{p\rest\omega\cdot\alpha}=w^p\cap\omega\cdot\alpha$. \hfill$\square$ 
\end{observation}
Thus $p\rest\omega\cdot\alpha$ corresponds to $p\rest\alpha$ as required in 
\ref{semiuni}(2c). The main difficulty of the proof is checking the
amalgamation property \ref{semiuni}(2d). Before we deal with this demand we
will check that some sets are dense in $\cA\cP$ (which will allow us to
simplify some arguments and will be of importance in drawing conclusions) and
we will deal with existing of some upper bounds.

\begin{claim}
\label{2.2D}
In $\cA\cP$, if $\langle p_i:i<\delta\rangle$ is increasing, $\delta<\lambda$
and for $i_0<i_1<\delta$ we have $w^{p_{i_0}}=w^{p_{i_1}}$ then its union
(defined naturally) is its least upper bound in $\cA\cP$. 
\hfill$\square$
\end{claim}

\begin{claim}[Density Observation]
\label{2.2F}
Assume $p\in\cA\cP$.
\begin{enumerate}
\item Suppose that $\alpha\in w^p$, $u\subseteq [\lambda{\cdot}\alpha,
\lambda{\cdot}\alpha+\lambda)\setminus t_p$, $|u|<\lambda$ and for $i\in u$
we are given a full branch $A_i$ of $(t^p\cap\lambda{\cdot}\alpha,\leq^p)$
(i.e. $A_i$ is linearly ordered by $\leq^p$ and $\beta\in w^p\cap\alpha
\Rightarrow [\lambda{\cdot}\beta,\lambda{\cdot}\beta+\lambda)\cap A_i\neq
\emptyset$), $i\neq j\Rightarrow A_i\neq A_j$. Furthermore, assume that
if $\alpha$ is limit then
\[\gamma\in t^p\cap[\lambda{\cdot}\alpha,\lambda{\cdot}\alpha+\lambda)\ \&\
i\in u\quad\Rightarrow\quad A_i\neq\{y\in t^p: y<^p\gamma\}.\]
Then there is $q\in\cA\cP$, $p\leq q$ such that $t^q=t^p\cup u$, $<^q=<^p\cup
\{(y,i): y\in A_i, i\in u\}$, and the rest is equal (i.e. $w^q=w^p$, $D^q=
D^p$, $\bar{d}^q=\bar{d}^p$, $\bar{e}^{q}=\bar{e}^{p}$, $f^q\supseteq f^p$,
$c^q\supseteq c^p$ naturally).
\item If $\alpha\in w^p$, $i\in [\lambda{\cdot}\alpha,\lambda{\cdot}\alpha+
\lambda)$ then there is $q\in\cA\cP$, $p\leq q$ such that $i\in t^q$,
$w^p=w^q$, $D^q=D^p$, $\bar{d}^q =\bar{d}^p$, $\bar{e}^q=e^p$, and naturally
$f^q\supseteq f^p$, $c^q\supseteq c^p$.
\item If $\bar{x}\in(t^p,\leq^p)^{[\mu]}$ and $h\in H^{1,*}_\mu\cup
H^{3,*}_\mu$ then there is $q\in\cA\cP$ such that $p\le q$ and $D^q=
D^p\cup\{(\bar{x},h)\}$.
\item If $(\bar{x},h)\in D^p$, $\bar{x}<\bar{y}\in (t^p,\le^p)^{[\mu]}$ then
for some $q\in\cA\cP$, $p\le q$ we have:\quad $\bar{y}\in\dom(d^q_{\bar{x},
h})=\dom(d^p_{\bar{x},h})\cup\{\bar{y}':\bar{x}<\bar{y}'\leq\bar{y},\
\bar{y}'\in(t^p,\leq^p)^{[\mu]}\}$, $t^p=t^q$, $w^p=w^q$, $\leq^p=\leq^q$,
$D^p=D^q$, $d^q_{\bar{x}',h'}=d^p_{\bar{x}',h'}$, $e^q_{\bar{x}',h'}=e^p_{
\bar{x}',h'}$ for $(\bar{x}',h')\in D^p\setminus\{(\bar{x},h)\}$, and
$f^q=f^p$, $c^q=c^p$. 
\end{enumerate}
\end{claim}

\noindent{\em Proof of the claim:}\hspace{0.2in} 1)\ \ \ Check.
\medskip

\noindent 2)\ \ \ Iterate (1) (in the $j$-th time -- on the $j$-th level
of $t$) using \ref{2.2D} for limit stage. More elaborately, for each $\beta\in
w^p\setminus\{0\}$ choose $i_\beta\in [\lambda\cdot\beta,\lambda\cdot\beta+
\lambda)\setminus t^p$  such that $i_\alpha=i$. Next by induction on $\beta\in
w^p\cap (\alpha+1)$ choose an increasing sequence $\langle p_\beta:\beta\in
w^p\cap (\alpha+1)\rangle\subseteq\cA\cP$ of approximations such that $p_0=p$
and $t^{p_\beta}=t^p\cup\{i_\gamma:0<\gamma\in w^p\cap(\beta+1)\}$ (so
$w^{p_\beta}=w^p$) and the sequence $\langle i_\gamma:0<\gamma\in w^p\cap
(\beta+1)\rangle$ is $\leq^{p_\beta}$--increasing. 
\medskip

\noindent 3)\ \ \ We just put: $t^q=t^p$, $w^q=w^p$, $f^q=f^p$, $c^q=c^p$,
$d^q_{\bar{x}',h'}=d^p_{\bar{x}',h'}$ if $(\bar{x}',h')\in D^p$,
$d^q_{\bar{x},h}$ is empty if $(\bar{x},h)\notin D^p$, and similarly for
$e^q_{\bar{x}',h'}$.
\medskip

\noindent 4)\ \ \ Let $\{\bar{y}_\vare:\vare<\zeta\}$ list $\dom(d^q_{\bar{x},
h})\setminus\dom(d^p_{\bar{x},h})$ in the $<^p$--increasing way and let
$\alpha^*=\sup(\rang(d^p_{\bar{x},h})\cup\{\gamma: (\exists\bar{y}')((
\bar{y}',\gamma)\in\dom(e^p_{\bar{x},h}))\})$. Now put
\[\dom(e^q_{\bar{x},h})=\{(\bar{y},\alpha):\bar{y}\in\dom(d^q_{\bar{x},h})\
\&\ \alpha<\alpha^*+1+\zeta\},\]
declare that $e^p_{\bar{x},h}\subseteq e^q_{\bar{x},h}$ and 
\begin{itemize}
\item if $\bar{y}\in\dom(d^q_{\bar{x},h})$, $\bar{y}<\bar{y}_\vare$, $\vare<
\zeta$ then $e^q_{\bar{x},h}(\bar{y},\alpha^*+1+\vare)=0$,
\item $e^q_{\bar{x},h}(\bar{y}_\vare,\alpha^*+1+\vare)=1$,
\item $e^q_{\bar{x},h}(\bar{y},\gamma)=2$ in all other instances.
\end{itemize}
It should be clear that this defines correctly an approximation $q\in\cA\cP$
and that it is as required. [Note that clauses ($\gamma$), ($\vare$) of
\ref{2.2}(2b) relevant to clauses (J), (K) in the definition of $\cA\cP$
hold by the requirement ``$d_{\bar{x},h}(\bar{t})<d_{\bar{x},h}(\bar{y})=
d_{\bar{x},h}(\bar{z})$''.] \hfill$\square_{\ref{2.2F}}$

\begin{claim}
\label{2.2G}
If $p\in\cA\cP$ and an ordinal $\alpha\in\lambda^+\setminus w^p$ is divisible
by $\omega$ then for some $q\in\cA\cP$, $p\leq q$ we have $w^q=w^p\cup
[\alpha,\alpha+\omega)$, $D^p=D^q$ and
\[\beta\in w^p\quad\Rightarrow\quad t^q\cap [\lambda{\cdot}\beta,\lambda
{\cdot}\beta+\lambda)=t^p\cap [\lambda{\cdot}\beta,\lambda{\cdot}\beta+
\lambda).\]
\end{claim}

\noindent{\em Proof of the claim:}\hspace{0.2in} Let $\beta=\min(w^p\setminus
\alpha)$ (if $\beta$ is undefined then it is much easier; of course $\beta>
\alpha$ as $\alpha\notin w^p$ and therefore $\beta\geq\alpha+\omega)$. Let
$t^p\cap[\lambda{\cdot}\beta,\lambda{\cdot}\beta+\lambda)=\{y^\beta_i:i<i^*\}$
be an enumeration with no repetitions and for $n<\omega$ let
\[\{y^{\alpha+n}_i: i<i^*\}\subseteq [\lambda{\cdot}(\alpha+n),\lambda{\cdot}
(\alpha+n)+\lambda)\]
be with no repetition. Let $t^q=t^p\cup\{y^{\alpha+n}_i:n<\omega, i<i^*\},$ and
\[\begin{array}{ll}
\leq^q=\leq^p\cup&\{(y^{\alpha+n}_i,y^{\alpha+m}_i): n\leq m<\omega,i<i^*\}
\cup\\
\ &\{(y^{\alpha+n}_i,x): n<\omega, i<i^*, y^\beta_i\leq^p x\}\cup\\
\ &\{(x,y^{\alpha+n}_i): n<\omega, i<i^*, x<^p y^\beta_i\}.\\
\end{array}\]
By \ref{2.2F} we may assume that $i^*=\mu$. For $(\bar{x},h)\in D^p=D^q$ we
let 
\[\begin{array}{ll}
\dom(d^q_{\bar{x},h})=&\dom(d^p_{\bar{x},h})\cup\{\bar{y}\in (t^q)^{[\mu]}: 
(\exists\bar{y}'\in\dom(d^p_{\bar{x},h}))(\bar{x}<\bar{y}<\bar{y}')\ \mbox{
and}\\ 
\ &\qquad(\exists n<\omega)(\forall\vare<\mu)(y_\vare\in [\lambda\cdot
(\alpha+n),\lambda\cdot(\alpha+n+1)))\},
  \end{array}\]
and let $\alpha^*_{\bar{x},h}=\sup(\rang(d^p_{\bar{x},h})\cup\{\gamma:
(\exists\bar{y}')((\bar{y}',\gamma)\in\dom(e^p_{\bar{x},h}))\})$. Fix an
enumeration $\{\bar{y}_\xi:\xi<\zeta\}$ of $\dom(d^q_{\bar{x},h})\setminus\dom
(d^p_{\bar{x},h})$ such that $\bar{y}_{\xi_0}<\bar{y}_{\xi_1}\ \ \Rightarrow\
\ \xi_0<\xi_1$. We put $d^q_{\bar{x},h}(\bar{y}_\xi)=\alpha^*_{\bar{x},h}+
1+\xi$ for $\xi<\zeta$ (and we declare $d^q_{\bar{x},h}\supseteq d^p_{\bar{x},
h}$). Next we define $e^q_{\bar{x},h}$ similarly as in \ref{2.2F}(4) putting
the value 2 whenever possible (so $\dom(e^q_{\bar{x},h})=\{(\bar{y},\gamma):
\bar{y}\in \dom(d^q_{\bar{x},h})\ \&\ \gamma<\alpha^*_{\bar{x},h}+1+\zeta\}$).
Now comes the main point: we have to define functions $f^q,c^q$ (extending
$f^p,c^p$, respectively) such that clauses (I) + (J) + (K) hold. But it should
be clear that each instance of clause (I) in $t^q$ can be reduced to an
instance of this clause in $t^p$ (just look at the definitions of $t^q,
d^q_{\bar{x},h},e^q_{\bar{x},h}$). Thus what we really have to take care of
are instances of (J) and (K). For this we define $c^q\rest\{y^{\alpha+n}_i:i<
\mu\}$ and $f^q\rest\{y^{\alpha+n+1}_i: i<\mu\}$ by induction on $n<\omega$. 
At the first stage (for $n=0$) we let 
\[\begin{array}{ll}
P=\big\{(\zeta,\bar{x},h,\bar{y},\bar{z}):&(\bar{x},h)\in D^q\mbox{ and }\zeta<
\mu\mbox{ and }\bar{x}<\bar{y}<\bar{z}\in\dom(d^q_{\bar{x},h}),\\
\ &\mbox{and }\bar{z}\subseteq \{y^{\alpha+1}_i:i<\mu\}\mbox{ and }\bar{y}
\subseteq\{y^\alpha_i:i<\mu\}\big\}.
  \end{array}\]
Take a list $\langle(X^\Upsilon,\zeta^\Upsilon,\bar{x}^\Upsilon,h^\Upsilon,
\bar{y}^\Upsilon,\bar{z}^\Upsilon):\Upsilon<\mu\rangle$ of 
\[\{(X,\zeta,\bar{x},h,\bar{y},\bar{z}): X\in \{J,K\}\ \&\ (\zeta,\bar{x},h,
\bar{y},\bar{z})\in P\}\]
in which each 6-tuple appears $\mu$ times and $\zeta^\Upsilon\leq
1+\Upsilon$. Next by induction on $\Upsilon<\mu$ choose a sequence $\langle
c_\Upsilon,f_\Upsilon:\Upsilon\leq\mu\rangle$ such that
\begin{description}
\item[$(\alpha)$] $c_\Upsilon:\dom(c_\Upsilon)\longrightarrow\mu$,
$\dom(c_\Upsilon)\subseteq\{y^\alpha_i:i<\mu\}$, $|\dom(c_\Upsilon)|<\aleph_0
+|\Upsilon|^+$,
\item[$(\beta)$] $f_\Upsilon:\dom(f_\Upsilon)\longrightarrow\lambda$,
$\dom(f_\Upsilon)\subseteq\{y^{\alpha+1}_i:i<\mu\}$, $|\dom(f_\Upsilon)|<
\aleph_0+|\Upsilon|^+$,
\item[$(\gamma)$] $\langle c_\Upsilon:\Upsilon\leq\mu\rangle$, $\langle
f_\Upsilon:\Upsilon\leq\mu\rangle$ are increasing continuous,
\item[$(\delta)$] for each $\Upsilon<\mu$ there is $i^\Upsilon<\mu$ divisible
by $\zeta^\Upsilon$ such that 
\[\rang(\bar{z}^\Upsilon\rest [i^\Upsilon,i^\Upsilon+\zeta^\Upsilon))\subseteq
\dom(f_{\Upsilon+1})\setminus\dom(f_{\Upsilon})\qquad\mbox{and}\]
\[\rang(\bar{y}^\Upsilon\rest [i^\Upsilon,i^\Upsilon+\zeta^\Upsilon))\subseteq
\dom(c_{\Upsilon+1})\setminus\dom(c_{\Upsilon}),\]
\item[$(\vare)$] if $X^\Upsilon=J$ and $h^\Upsilon\in H^{3,*}_\mu$ then
condition \ref{2.2C}(J) holds for $\bar{x}^\Upsilon$, $h^\Upsilon$, 
$\bar{y}^\Upsilon$, $\bar{z}^\Upsilon$ with $i=i^\Upsilon$,
\item[$(\zeta)$] if $X^\Upsilon=K$ then $c_{\Upsilon+1}((\bar{y}^\Upsilon)_{
i^\Upsilon})=\zeta^\Upsilon,$ 
\item[$(\iota)$] $y^\alpha_\Upsilon\in\dom(c_{\Upsilon+1})$,
$y^{\alpha+1}_\Upsilon\in\dom(f_{\Upsilon+1})$.
\end{description}
There are no difficulties with carrying out the construction: the only
possible troubles could come from demand $(\vare)$ above. But look at the
definition \ref{2.1}(14) of $H^3_\mu$. Taking sufficiently large $\beta<\mu^+
=\lambda$, the respective sequences $\langle\alpha^0_\vare:\vare<
\zeta^\Upsilon\rangle$, $\langle\alpha^1_\vare:\vare<\zeta^\Upsilon\rangle$
will be good candidates for $c_{\Upsilon+1}\rest (\bar{y}^\Upsilon\rest
[i^\Upsilon,i^\Upsilon+\zeta^\Upsilon))$ and $f_{\Upsilon+1}\rest (
\bar{z}^\Upsilon\rest [i^\Upsilon,i^\Upsilon+\zeta^\Upsilon))$ in clause
$(\vare)$.\\ 
The functions $c_\mu,f_\mu$ will be the respective restrictions 
$c^q\rest \{y^{\alpha}_i:i<\mu\}$ and $f^q\rest\{y^{\alpha+1}_i:
i<\mu\}$. Next, arriving to a stage $n+1$ of the definition we repeat the
above procedure with no changes. Note that at this stage we know $c^q\rest 
\{y^{\alpha+n}_i:i<\mu\}$, $f^q\rest\{y^{\alpha+n+1}_i: i<\mu\}$ but they have
no influence on defining $c^q,f^q$ at levels $\alpha+n+1$, $\alpha+n+2$. 
\hfill$\square_{\ref{2.2G}}$

\begin{claim}[The Amalgamation Property]
\label{2.2H}
Assume that $\alpha<\lambda^+$ is either a successor ordinal or $\cf(\alpha)=
\lambda$, $p,q\in\cA\cP$ and $p\rest\omega\cdot\alpha\leq q\in\cA\cP_\alpha$. 
Then there is $r\in\cA\cP$ such that $p,q\leq r$.
\end{claim}

\noindent{\em Proof of the claim:}\hspace{0.2in} First try just the $r$
defined by:
\[w^r=w^p\cup w^q,\quad t^r=t^p\cup t^q,\quad\leq^r=\leq^p\cup\leq^q,\quad
D^r=D^p\cup D^q,\]
\[\begin{array}{ll}
d^r_{\bar{x},h}\ \ \mbox{ is:}\; & d^p_{\bar{x},h}\quad\mbox{ if }(\bar{x},h)\in
D^p\setminus D^q\\
\ & d^q_{\bar{x},h}\quad\mbox{ if }(\bar{x},h)\in D^q\setminus D^p\\
\ & d^p_{\bar{x},h}\cup d^q_{\bar{x},h}\quad\mbox{ if }(\bar{x},h)\in D^q
\cap D^p
\end{array}\]
and $e^r_{\bar{x},h}\supseteq e^p_{\bar{x},h}\cup e^q_{\bar{x},h}$ (defined
naturally, i.e.~with $\dom(e^r_{\bar{x},h})$ minimal possible to satisfy
demand (F) and value 2 whenever possible), $f^r=f^p\cup f^q$, $c^r=c^p\cup
c^q$. Clearly $p\le r$, $q\le r$, $w^r=w^p\cup w^q$, but does $r$ belong
to $\cA\cP$? The things that might go wrong are:
\begin{itemize}
\item there is $y\in t^p\setminus t^q$ which has nothing below it in
some levels, or
\item $y\wedge z$ not defined for some $y,z$, or
\item the relevant cases of clauses (I)--(K) fail.
\end{itemize}
Let $\beta_0=\bigcup\{\gamma:\omega\cdot\gamma\in w^p\cap\omega\cdot\alpha\}=
\sup(\dom(p)\cap\alpha)$. Note that $\beta_0\in\dom(p)\cap\alpha$ (as $\alpha$
is either successor or of cofinality $\lambda$) and
\begin{description}
\item[($\circledast$)] \quad if $w^q\subseteq\omega\cdot\beta_0+\omega$ then 
$r\in\cA\cP$.
\end{description}
So we assume from now on that $w^q\not\subseteq\omega\cdot\beta_0+\omega$ (by
($\circledast$) above). Then necessarily $\beta_0+1<\alpha$ (as $w^q\subseteq
\omega\cdot\alpha$). Without loss of generality $\dom(p)\setminus\alpha\neq
\emptyset$ (as if $w^p\subseteq\omega\cdot\alpha$ we can let $r=q$) and
$q\rest\omega\cdot(\beta_0+1)=p\rest\omega\cdot\alpha$. Let $\alpha^*
\stackrel{\rm def}{=}\min(\dom(p)\setminus\alpha)$ and $\beta^*\stackrel{\rm
def}{=}\min(\dom(q)\setminus(\beta_0+1))$. By \ref{2.2G} we may assume that  
$\beta^*=\beta_0+1$ (i.e.~$\omega\cdot\beta_0+\omega\in w^q$). Let $\{x^{
\omega\cdot\alpha^*}_i: i<i^*\}$ list $t^p\cap [\lambda\cdot(\omega\cdot
\alpha^*),\lambda\cdot(\omega\cdot\alpha^*)+\lambda)$. By \ref{2.2F}(1) 
(i.e.~increasing $q$ only by increasing $t^q\cap [\lambda\cdot(\omega\cdot
\beta^*),\lambda\cdot((\omega\cdot\beta^*)+1))$ we may assume that there is a
list $\{x^{\omega\cdot\beta^*}_i: i<i^*\}$ of distinct members of $t^q\cap
[\lambda\cdot(\omega\cdot\beta^*),\lambda\cdot(\omega\cdot\beta^*)+\lambda)$
such that $(\forall z\in t^p\cap t^q)[z<^p x^{\omega\cdot\alpha^*}_i\equiv 
z<^q x^{\omega\cdot\beta^*}_i]$. Let $x^\beta_i\in [\lambda{\cdot}\beta,
\lambda{\cdot}\beta+\lambda)$ (for $\beta\in w^q\setminus (\omega\cdot\beta^*+
1)$ and $i<i^*$) be pairwise distinct and not in $t^q$. Now we shall
``correct'' $r $ to $r^*$: 
\[t^{r^*}=t^r\cup\{x^\beta_i:\beta\in w^q\setminus(\omega\cdot\beta^*+1),i<
i^*\},\quad\quad w^{r^*}=w^r,\]
\[\begin{array}{ll}
\leq^{r^*}=\leq^r\cup &\{(x^\beta_i,x):i<i^*,x\in t^p,x^{\omega\cdot
\alpha^*}_i\leq^p x,\beta\in w^q\setminus\omega\cdot\beta^*\}\cup\\
\ &\{(x,x^\beta_i):i<i^*,x\in t^q,x\leq^q x^{\omega\cdot\beta^*}_i,\beta\in 
w^q\setminus\omega\cdot\beta^*\}\cup\\
\ &\{(x^{\beta_0}_i,x^{\beta_1}_i): \beta_0,\beta_1\in w^q\setminus\omega\cdot
\beta^*,\beta_0\leq\beta_1, i<i^*\}.\\
\end{array}\]
Put $D^{r^*}=D^r$. If $(\bar{x},h)\in D^r\setminus D^{p\rest\omega\cdot
\alpha}$ then we can let $d^{r^*}_{\bar{x},h}=d^r_{\bar{x},h}$, but if
$(\bar{x},h)\in D^{p\rest\omega\cdot\alpha}$ then we first let
\[\gamma^*_{\bar{x},h}=\sup(\rang(d^r_{\bar{x},h})\cup\{\gamma:(\exists
\bar{y}')((\bar{y},\gamma)\in\dom(e^p_{\bar{x},h})\cup\dom(e^q_{\bar{x},h}))\}
)\qquad\mbox{ and}\]
\[\begin{array}{ll}
\dom(d^{r^*}_{\bar{x},h})=&\{\bar{y}\in (t^{r^*})^{[\mu]}:\bar{y}\in\dom(
d^r_{\bar{x},h})\ \ \mbox{ or \ \ for some}\\
\ &\ \bar{z}=\langle x^{\omega\cdot\alpha^*}_{\vare(j)}\!:j\!<\!\mu\rangle\in
(t^p)^{[\mu]}_{\otp(w^p\cap\omega{\cdot}\alpha^*)}\mbox{ and }\beta\in w^q
\setminus (\omega{\cdot}\beta^*+1)\\
\ &\ \mbox{we have}\ \ \ \bar{z}\in\dom(d^p_{\bar{x},h})\ \mbox{ and }\ 
\bar{y}=\langle x^\beta_{\vare(j)}:j<\mu\rangle\}.
\end{array}\]
Choose $d^{r^*}_{\bar{x},h}$ in such a manner that $d^{r^*}_{\bar{x},h}
\supseteq d^r_{\bar{x},h}$ and the values $d^{r^*}_{\bar{x},h}(\bar{y})$, if
not defined before, are distinct ordinals from $(\gamma^*_{\bar{x},h},
\lambda)$. Thus, in particular,
\[d^{r^*}_{\bar{x},h}(\bar{y})=d^{r^*}_{\bar{x},h}(\bar{z})\ \&\ \bar{y}\neq
\bar{z}\quad\Rightarrow\quad\{\bar{y},\bar{z}\}\subseteq\dom(d^r_{\bar{x},
h}).\]
Next we define $e^{r^*}_{\bar{x},h}$ extending $e^{r}_{\bar{x},h}$ to 
satisfy clause (F) -- we put the value 2 whenever it is possible. [Note that
this is the place in which the assumption that $\bar{y}\in t^{\langle\mu
\rangle}_{\bar{x}}$ in clause (F)(ii), and so the respective assumption in
\ref{2.2}(2b), play role: the values of $e^q_{\bar{x},h}$ at the level
$\omega\cdot\beta^*=\omega\cdot\beta_0+\omega$ do not interfere with the
values of $e^p_{\bar{x},h}$ at the level $\omega\cdot\alpha^*$ since $\beta^*$
is a successor, not of cofinality $\lambda$.] Now we have to define
$c^{r^*}\supseteq c^r$, $f^{r^*}\supseteq f^r$, i.e.~to define 
\[\begin{array}{l}
c^{r^*}\rest\{x^\beta_i:i<i^*,\beta\in w^q\setminus (\omega\cdot\beta^*+1)\},
\quad\quad\mbox{and}\\
f^{r^*}\rest\{x^\beta_i:i<i^*,\beta\in w^q\setminus\omega\cdot\beta^*\mbox{ is
a successor}\},
\end{array}\]
such that clauses \ref{2.2C}(G)--(K) hold. This is done like in \ref{2.2G},
but now the clause (I) is ``active'' too. Of course, the point is that we have
$\mu$ commitments, each has ``$\mu$ disjoint chances", so we list them in a
list of length $\mu$ and inductively we can easily do it (for $\mu$ singular
-- at place $i$ of the list there may appear only a commitment of ``size
$\leq|i|+\aleph_0$"). More fully, let 
\[\begin{array}{ll}
P^1=&\big\{(\zeta,\bar{x},h,\bar{y},\bar{z}):(\bar{x},h)\in D^p\cap D^q\mbox{
and }\zeta<\mu\mbox{ and for some }\vare\\
\ &(\bar{z},\vare)\in\dom(e^{r^*}_{\bar{x},h})\setminus\dom(e^r_{\bar{x},h})
\mbox{ and } (\bar{y},\vare)\in\dom(e^q_{\bar{x},h})\mbox{ and}\\
\ &\ e^{r^*}_{\bar{x},h}(\bar{z},\vare)\leq 1\mbox{ and }e^q_{\bar{x},h}(
\bar{y},\vare)\leq 1\mbox{ and}\\
\ &\bar{y}\subseteq [\lambda(\omega\cdot\beta^*+1),\lambda(\omega\cdot\beta^*+
1)+\lambda)\mbox{ and }\bar{z}\subseteq\{x^{\omega\cdot\beta^*+1}_i:i<i^*\}
\big\}.\\
\end{array}\]
Defining $f$ we have to take care of condition (I) for all $(\zeta,\bar{x},
h,\bar{y},\bar{z})\in P^1$. We also have to take care of conditions (J), (K)
for 
\[\begin{array}{ll}
P^2=&\big\{(\zeta,\bar{x},h,\bar{y},\bar{z}):(\bar{x},h)\in D^p\cap D^q\mbox{
and }\zeta<\mu\mbox{ and}\\
\ &\bar{z}\in\dom(d^{r^*}_{\bar{x},h}),\ \bar{z}\subseteq\{x^{\gamma+1}_i:
i<i^*\},\ \bar{y}=\bar{z}\rest\gamma,\ \gamma\in w^q\setminus\omega\cdot
\beta^*\}.\\
\end{array}\]
So we use a list $\langle(X^\Upsilon,\zeta^\Upsilon,\bar{x}^\Upsilon,
h^\Upsilon,\bar{y}^\Upsilon,\bar{z}^\Upsilon): \Upsilon<\mu\rangle$ of 
$\{(X,\zeta,\bar{x},h,\bar{y},\bar{z}):X\in\{1,2\}$ and $(\zeta,
\bar{x},h,\bar{y},\bar{z})\in P^1\cup P^2\}$ in which each 6-tuple appears
$\mu$ times and $\zeta^\Upsilon\leq 1+\Upsilon$. Let $\{x^\Upsilon:
\Upsilon<\mu\}$ list $t^{r^*}\setminus t^r$. Now we define by induction on
$\Upsilon\leq\mu$ functions $c_\Upsilon$, $f_\Upsilon$ such that
\begin{description}
\item[$(\alpha)$] $c_\Upsilon $ is a function extending $c^r$,
$\rang(c_\Upsilon)\subseteq\lambda$, 
\item[$(\beta)$]  $\dom(c_\Upsilon)\setminus\dom(c^r)$ is a subset of 
$\{x^\beta_i:i<i^*,\beta\in w^q\setminus(\omega\cdot\beta^*+1)\}$ of
cardinality $<\aleph_0+|\Upsilon|^+$, 
\item[$(\gamma)$] $f_\Upsilon$ is a function extending $f^r$,
$\rang(f_\Upsilon)\subseteq\lambda$, 
\item[$(\delta)$]  $\dom(f_\Upsilon)\setminus\dom(f^r)$ is a subset of 
\[\{x^\beta_i:i<i^*,\beta\mbox{ is a successor ordinal and }\beta\in w^q
\setminus (\omega\cdot\beta^*+1)\}\]
of cardinality $<\aleph_0+|\Upsilon|^+$,
\item[$(\vare)$]  the sequences $\langle c_\Upsilon:\Upsilon\leq\mu\rangle$,
$\langle f_\Upsilon:\Upsilon\leq\mu\rangle$ are increasing continuous,
\item[$(\zeta)$]  for each $\Upsilon$ there is $i^\Upsilon<\mu$ divisible by
$\zeta^\Upsilon$ such that:

if $X^\Upsilon=I$ then 
\[\hspace{-0.3cm}\rang\big(\bar{z}^\Upsilon\rest[i^\Upsilon,i^\Upsilon+
\zeta^\Upsilon)\big)\subseteq [\dom(c_{\Upsilon+1})\setminus\dom(c_\Upsilon)]
\cap[\dom(f_{\Upsilon+1})\setminus\dom(f_\Upsilon)],\]
and if $X^\Upsilon\in\{J,K\}$, $(\zeta^\Upsilon,\bar{x}^\Upsilon,h^\Upsilon,
\bar{y}^\Upsilon,\bar{z}^\Upsilon)\in P^2$ then 
\[\rang(\bar{z}^\Upsilon\rest [i^\Upsilon,i^\Upsilon+\zeta^\Upsilon))\subseteq
\dom(f_{\Upsilon+1})\setminus\dom(f_\Upsilon),\qquad\mbox{and}\]
\[\rang(\bar{y}^\Upsilon\rest [i^\Upsilon,i^\Upsilon+\zeta^\Upsilon))\subseteq
(\dom(c_{\Upsilon+1})\setminus\dom(c_\Upsilon))\cup c^r,\]
\item[$(\iota)$] if $X^\Upsilon=I$ and $(\zeta^\Upsilon,\bar{x}^\Upsilon,
h^\Upsilon,\bar{y}^\Upsilon,\bar{z}^\Upsilon)\in P^1$ then condition
\ref{2.2C}(I) holds for $(\zeta^\Upsilon,\bar{x}^\Upsilon,h^\Upsilon,
\bar{y}^\Upsilon,\bar{z}^\Upsilon,i^\Upsilon)$,
\item[$(\kappa)$] if $X^\Upsilon=J$, $(\zeta^\Upsilon,\bar{x}^\Upsilon,
h^\Upsilon,\bar{y}^\Upsilon,\bar{z}^\Upsilon)\in P^2$  and $h^\Upsilon\in
H^{3,*}_\mu$ then condition \ref{2.2C}(J) holds for $(\bar{x}^\Upsilon,
h^\Upsilon,\bar{y}^\Upsilon,\bar{z}^\Upsilon,i^\Upsilon)$,
\item[$(\lambda)$] if $X^\Upsilon=K$, $\bar{y}^\Upsilon\not\subseteq
\dom(c^r)$ then $c_{\Upsilon+1}((\bar{y}^\Upsilon)_{i^\Upsilon})=
\zeta^\Upsilon$, 
\item[$(\mu)$] $x^\Upsilon\in\dom(c_{\Upsilon+1})$ and if $x^\Upsilon$ is 
from a successor level of $t^{r^*}$ then $x^\Upsilon\in\dom(f_{\Upsilon+1})$,
\item[$(\nu)$] $\rang(f^\Upsilon\rest\{x^{\omega\cdot\beta^*+1}_i:i<i^*\})\cap
\rang(f^q)=\emptyset$.
\end{description}
The definition is carried out as in \ref{2.2G}. The new points are clause
$(\iota)$ and instances of clause $(\kappa)$ for $\Upsilon$ such that
$\bar{y}^\Upsilon\subseteq\{x^{\omega\cdot\beta^*}_i:i<i^*\}$. In the second
case a potential trouble could be caused by the fact that the function
$c_\Upsilon$ is defined on $\bar{y}^\Upsilon$ already. But the definition
\ref{2.1}(14) of $H^3_\mu$ was exactly what we need to handle this: we may
find suitable values for $f_{\Upsilon+1}\rest (\bar{z}\rest [i^\Upsilon,
i^\Upsilon,\zeta^\Upsilon))$. To deal with clause $(\iota)$ note that if
$h^\Upsilon\in H^{3,*}_\mu$ then demand \ref{2.2C}(J) for $q$ provides the
needed candidates for values of $f_{\Upsilon+1}$; if $h^\Upsilon\in
H^{1,*}_\mu$ then the definition \ref{2.1}(12) of $H^1_\mu$ works. 

The functions $c_\mu$, $f_\mu$ are as required. \hfill$\square_{\ref{2.2H}}$
\medskip

The demands (e)--(h) of \ref{semiuni}(2) are easy now:
\begin{claim}
\label{cl6}
\begin{enumerate}
\item If a sequence $\langle p_i:i<\delta\rangle\subseteq\cA\cP$ is increasing,
$\delta<\lambda$ then it has an upper bound $q\in\cA\cP$ such that
$\dom(q)=\cl(\bigcup\limits_{i<\delta}\dom(p_i)))$.
\item Assume $\beta<\lambda^+$, $\cf(\beta)=\lambda$, $\delta<\lambda$. Let 
$\langle p_i:i<\delta\rangle\subseteq \cA\cP_{\beta+1}$ be an increasing
sequence and let $q\in\cA\cP_\beta$ be an upper bound to $\langle p_i\rest
\omega\cdot\beta:i<\delta\rangle$. Then the family $\{p_i:i<\delta\}\cup\{q\}$
has an upper bound $r$ such that $r\rest\omega\cdot\beta\geq q$.
\item Assume that $\langle\beta_i: i<\delta\rangle\subseteq\lambda^+$ is
strictly increasing, each $\beta_i$ is either a successor or has cofinality
$\lambda$, $\delta<\lambda$ is limit. Suppose that $q\in\cA\cP$ and $\langle
p_i:i<\delta\rangle\subseteq\cA\cP$ is an increasing sequence such that 
\[(\forall i<\delta)(q\rest\omega\cdot\beta_i\leq p_i\in\cA\cP_{\beta_i}).\]
Then the family $\{p_i:i<\delta\}\cup\{q\}$ has an upper bound $r\in\cA\cP$
such that $(\forall i<\delta)(p_i\leq r\rest\omega\cdot\beta_i)$.
\item Suppose that $\delta_1,\delta_2<\lambda$ are limit ordinals and
$\langle\beta_j: j<\delta_2\rangle\subseteq\lambda^+$ is a strictly increasing
sequence of ordinals, each $\beta_j$ either a successor or of cofinality
$\lambda$. Let 
\[\langle p_{i,j}: (i,j)\in (\delta_1+1)\times (\delta_2+1)\setminus\{(
\delta_1,\delta_2)\}\rangle\subseteq\cA\cP\]
be such that 
\[p_{i,j}\in\cA\cP_{\beta_j},\quad i\leq i'\ \Rightarrow\ p_{i,j}\leq
p_{i',j},\quad j\leq j'\ \Rightarrow\ p_{i,j}\leq p_{i,j'}\rest \omega\cdot
\beta_j.\]
Then the family $\{p_{i,j}:(i,j)\in (\delta_1+1)\times (\delta_2+1)\setminus
\{(\delta_1,\delta_2)\}\}$ has an upper bound $r\in\cA\cP$ such that
$(\forall j<\delta_2)(r\rest\omega\cdot\beta_j=p_{\delta_1,j})$.
\end{enumerate}
\end{claim}

\noindent{\em Proof of the claim:}\hspace{0.2in} 1)\ \ \ The first try may be
to take the natural union of the sequence $\langle p_i:i<\delta\rangle$. 
However, it may happen that we will not get a legal approximation, as
$\bigcup\limits_{i<\delta}\dom(p_i)$ does not have to be closed. But we may
take its closure $\cl(\bigcup\limits_{i<\delta}\dom(p_i))$ and apply a
procedure similar to the one described in \ref{2.2G} (successively at each
element of $\cl(\bigcup\limits_{i<\delta}\dom(p_i))\setminus\bigcup\limits_{
i<\delta}\dom(p_i)$) and construct the required $q$.
\medskip

\noindent 2)--4)\ \ \ Similarly as 1) above plus the proof of \ref{2.2H}. 
\hfill$\square_{\ref{cl6}}$
\medskip

Now we apply \cite[Appendix]{Sh:405}: we find a ``sufficiently generic''
$G\subseteq \cA\cP$ which gives the $T, c, d$ we need (remember \ref{2.2F}): 
\[T=\{\eta_\vare:\vare\in t^p\ \mbox{ for some }p\in G\}\]
where for $\vare\in [\lambda\alpha,\lambda\alpha+\lambda)$ we define 
$\eta_\vare\in {}^\alpha\lambda$ by:

$\gamma=\eta_\vare(\beta)$\qquad if and only if 
\[(\exists p\in G)(\exists\vare'\in [\lambda{\cdot}(\beta+1),\lambda{\cdot}
(\beta+1)+\lambda))(t^p\models\mbox{``}\vare'<\vare\mbox{''}\ \&\ 
f^p(\vare')=\gamma).\]
This finishes the proof. \QED$_{\ref{2.2}}$

\begin{remark}
\label{2.2K}
{\em
\begin{enumerate}
\item Theorem \ref{2.2} is close to \cite{Sh:50}, which is a strengthening of
the construction of special Aronszajn trees. There essentially we replace
$(\gamma)+(\delta)$ by
\begin{description}
\item[$(\gamma)'$] $y_i(\lev(t_i)), z_i (\lev(t_i))$ do not depend on $i$.
\end{description}
\item By the proof of \ref{2.2}(2), $T$ is $\lambda$-complete.
\item We may add to \ref{2.2}(2):
\begin{description}
\item[(c)] $T$ is special, i.e.~there is a function $d:\bigcup\limits_\alpha
T_\alpha\longrightarrow\lambda$ such that
\[d(\eta)=d(\rho)\quad\Rightarrow\quad\neg[\eta<_T \rho].\]
\end{description}
[Just in the definition of $p\in{\cA\cP}$ (see \ref{2.2C}) add such $d^p:
t^p\longrightarrow\lambda$.]
\item The reader may wonder why we need the condition ``$h(\bar{x})$ has
$\mu$ pairwise disjoint members''. The point is that when we amalgamate $p$
and $q$ when $p\rest\alpha\leq q$, it may happen that $p$ gives information on
levels $\alpha_n<\alpha_{n+1}$ (for $n<\omega$), $\beta=\bigcup\limits_{n<
\omega}\alpha_n<\alpha$, $q$ gives information on the level $\beta$, and when
amalgamating the function $f^q$ gives information on $f$ on this level and $f$
is supposed to be one-to-one on every $\suc_T(\eta)$. So considering
$\bar{x}\in(t^q)_{\beta+1}$, $\bar{y}\in(t^p)_\alpha$, when we try to define
$f\rest\rang(\bar{y}\rest(\beta+1))$, some values are excluded.
\end{enumerate}
}
\end{remark}

\begin{theorem}
\label{2.3}
Assume that $\lambda=\mu^+=2^\mu$. Then there is a forcing notion ${\Bbb P}$
which is $(<\lambda)$--complete of size $\lambda^+$ and satisfies the
$\lambda^+$--cc (so it preserves cardinalities, cofinalities and cardinal
arithmetic) and such that in $\V^{\Bbb P}$:
\begin{enumerate}
\item There is a $\mu$-entangled linear order $\cI$ of cardinality $\lambda^+$
and density $\lambda$. 
\item Let $\sigma\leq\mu$ be a regular cardinal. There exist linear orders
$\cI'$, $\cI''$ of the cardinality $\lambda^+$ such that for any uniform
ultrafilter $D$ on $\sigma$ the linear orders $(\cI')^\sigma/D$,
$(\cI'')^\sigma/D$ have isomorphic subsets of cardinality $\lambda^+$, but
$\cI'+\cI''$ is $\mu$-entangled.\\
Hence there is a Boolean algebra $B$ which is $\lambda^+$-narrow but
$B^\sigma/D$ is not $\lambda^+$-narrow for any uniform ultrafilter $D$ on
$\sigma$.
\item There are a set $R\subseteq {}^\lambda\lambda$, $|R|=\lambda^+$ and
functions $c,\bar{d}$ such that, letting $T^+=({}^{\lambda>}\lambda\cup R,
\vartriangleleft)$ ($\vartriangleleft$ is being initial segment), we have:
\begin{description}
\item[(a)] $c$ is a function from ${}^{\lambda>}\lambda$ to $\lambda$,
\item[(b)] $R=\{\eta_\alpha:\alpha<\lambda^+\}$ (with no repetition),
$<_R=\{(\eta_\alpha,\eta_\beta):\alpha<\beta\}$; define
\[R^*=\{\langle\eta_{\alpha_i}:i<\mu\rangle: \alpha_i<\lambda^+,\; \langle
\alpha_i: i<\lambda^+\rangle\ \mbox{ is increasing}\},\]
\item[(c)] for every $\bar{x}\in T^{[\mu]}_{<\lambda}$, $\zeta<\mu$, and $h\in
{}^\zeta(\lambda\times\lambda\times\lambda)$, $d_{\bar{x},h}$ is a function
from $\{\bar{y}\in R^*:\bar{x}<\bar{y}\}$ to $\lambda$ such that:

\quad [$d_{\bar{x},h}(\bar{y})=d_{\bar{x},h}(\bar{z})\ \&\ \bar{y},\bar{z}\in
T^{[\mu]}_{\bar{x}}$ are distinct]\qquad implies

\quad $\sup\{\alpha:\eta_\alpha\mbox{ appears in }\bar{y}\}\neq\sup\{\alpha:
\eta_\alpha\mbox{ appears in }\bar{z}\}$\\
and for some $\bar{t}\in T^{[\mu]}_{\bar{x}}\cap T^{[\mu]}_{<\lambda}$ and
$i^*<\mu$ we have:
\begin{description}
\item[($\alpha$)]  $t_i=y_i\wedge z_i$, $\lev(t_i)=\lev(t_{i^*})$ for
$i\ge i^*$, 
\item[($\beta$)]   $(\forall\vare<\mu)(\exists^\mu i<\mu)(c(t_i)=\vare)$,
\item[($\gamma$)]  for $\mu$ ordinals $i<\mu$ divisible by $\zeta$ we have
\begin{description}
\item[(i)] either there are $\xi_0<\xi_1<\lambda$ such that
\[(\forall\vare<\zeta)(\zeta\cdot\xi_0\leq y_{i+\vare}(\lev(\bar{t}))<
\zeta\cdot\xi_1\leq z_{i+\vare}(\lev(\bar{t})))\quad\mbox{ and}\]
\[h=\langle(c(t_{i+\vare}),y_{i+\vare}(\lev(\bar{t}))-\zeta\cdot\xi_0, z_{i+
\vare}(\lev(\bar{t}))-\zeta\cdot\xi_1): \vare<\zeta\rangle,\]
\item[(ii)] or a symmetrical requirement interchanging $\bar{y}$ and
$\bar{z}$. 
\end{description}
\end{description}
\end{description}
\end{enumerate}
\end{theorem}

\Proof 1)\ \ \ We apply \ref{2.3}(3): let $R,c,\bar{d}$ be as there (of course
we are in the universe $\V^{\Bbb P}$ all time). We define the order $<_{\cI}$
on $\cI=R$ by 
\begin{quotation}
\noindent $y<_{\cI}z$\quad if and only if

either $c(y\wedge z)=0$ and $y(\alpha)<z(\alpha)$

or $c(y\wedge z)\neq 0$ and $y(\alpha)>z(\alpha)$,

\noindent where $\alpha=\lev(y \wedge z)$.
\end{quotation}
Clearly  $<_{\cI}$ is a linear order of the density $\lambda$,
$|\cI|=\lambda^+$. To show that it is $\mu$-entangled suppose that
$y^\alpha_\vare\in R$ (for $\alpha<\lambda^+$, $\vare<\vare(*)<\mu$) are
pairwise distinct, $u\subseteq\vare(*)$. Let $y^\alpha_\vare=\eta_{\beta(
\alpha,\vare)}$ (for $\alpha<\lambda^+$, $\vare<\vare(*)$). We may assume that
the truth value of ``$\beta(\alpha,\vare_1)<\beta(\alpha,\vare_2)$'' does not
depend on $\alpha<\lambda^+$. For simplicity, we may assume that for each
$\alpha<\lambda^+$, $\vare<\vare'<\vare(*)$ implies $y^\alpha_\vare<_R
y^\alpha_{\vare'}$. Finally, without loss of generality we may assume that
if $\alpha<\alpha'<\lambda^+$ and $\vare,\vare'<\vare(*)$ then $y^\alpha_\vare
<_R y^{\alpha'}_{\vare'}$ (i.e.~$\beta(\alpha,\vare)<\beta(\alpha',\vare')$).
For $\vare<\vare(*)$, $i<\mu$ and $\alpha<\lambda^+$ let $z^\alpha_{i{\cdot}
\vare(*)+\vare}=y^{\alpha{\cdot}\mu+i{\cdot}\vare(*)}_{\vare}$ and
$\bar{z}^\alpha=\langle z^\alpha_i: i<\mu\rangle$. Clearly each
$\bar{z}^\alpha$ is in $R^*$. Now for $\alpha<\lambda^+$ choose
$\xi(\alpha)<\lambda$ such that $z^\alpha_i\rest\xi(\alpha)$ (for $i<\mu$)
are pairwise distinct. Without loss of generality we may assume that
$\xi(\alpha)=\xi$ for $\alpha<\lambda^+$. There are $\lambda^\mu=\lambda$
possibilities for $\langle z^\alpha_i\rest\xi: i<\mu\rangle$, so we may assume
that for all $\alpha<\lambda^+$
\[\langle z^\alpha_i\rest\xi: i<\mu\rangle=\langle x_i: i<\mu\rangle=
\bar{x}\in T^{[\mu]}_\xi.\]
Let $h\stackrel{\rm def}{=}\langle(0,\ell^u_\vare,1-\ell^u_\vare): \vare<
\vare(*)\rangle\in{}^{\vare(*)}(\lambda\times\lambda\times\lambda)$, where
$\ell^u_\vare$ is 0 if $\vare\in u$ and 1 otherwise. For some distinct
$\alpha_1,\alpha_2<\lambda^+$ we have $d_{\bar{x},h}(\bar{z}^{\alpha_1})=
d_{\bar{x},h}(\bar{z}^{\alpha_2})$. By the properties of $d_{\bar{x},h}$,
possibly interchanging $\alpha_1,\alpha_2$, we find $i<\mu$, ordinals $\xi_0<
\xi_1<\lambda$ and $\bar{t}\in T^{[\mu]}_{\bar{x}}$ such that
\[\begin{array}{l}
(\forall\vare<\vare(*))(z^{\alpha_1}_{i{\cdot}\vare(*)+\vare}\wedge
z^{\alpha_2}_{i{\cdot}\vare(*)+\vare}=t_{i{\cdot}\vare(*)+\vare})\qquad
\mbox{ and}\\
(\forall\vare<\vare(*))(\vare(*)\cdot\xi_0\leq z^{\alpha_1}_{i\cdot\vare(*)+
\vare}(\beta)<\vare(*)\cdot\xi_1\leq z^{\alpha_2}_{\i\cdot\vare(*)+\vare}(
\beta)\quad\mbox{ and}\\
h=\langle(c(t_{i{\cdot}\vare(*)+\vare}),z^{\alpha_1}_{i{\cdot}\vare(*)
+\vare}(\beta)-\vare(*)\cdot\xi_0,z^{\alpha_2}_{i{\cdot}\vare(*)+\vare}(\beta)
-\vare(*)\cdot\xi_1): \vare<\vare(*)\rangle,
\end{array}\]
where $\beta=\lev(\bar{t})$. Then $\alpha_1{\cdot}\mu+i{\cdot}\vare(*)\neq
\alpha_2{\cdot}\mu+i{\cdot}\vare(*)$ (for $i<\mu$) and
\[(\forall\vare<\vare(*))(y^{\alpha_1{\cdot}\mu+i{\cdot}\vare(*)}_{
\vare}<_{\cI}y^{\alpha_2{\cdot}\mu+i{\cdot}\vare(*)}_{\vare}\ \ \mbox{
iff }\ \ \vare\in u)\]
(by the definition of the order $<_{\cI}$ and the choice of $h$), so using
\ref{1.2}(7) we are done. 
\medskip

\noindent 2)\ \ \ As in 1) above, we work in $\V^{\Bbb P}$ and we use
\ref{2.3}(3). Suppose $\sigma\leq\mu$. For a set $A\subseteq\lambda$ we define
the order $<_A$ on $R$ by 
\begin{quotation}
\noindent $y<_A z$\quad if and only if

either $c(y\wedge z)\in A$ and $y(\alpha)<z(\alpha)$

or $c(y\wedge z)\notin A$ and $z(\alpha)<y(\alpha)$,

\noindent where $\alpha=\lev(y\wedge z)$.
\end{quotation}
[Note that the order $<_{\cI}$ from part 1) above is just
$<_{\{0\}})$.]\\
Clearly $<_A$ is a linear order. As in the proof of \ref{2.3}(1) one can show
that it is $\mu$-entangled. As $\sigma^{<\sigma}\leq\lambda$ we may choose
sets $A_i\subseteq\lambda$ for $i\leq\sigma$ such that
\begin{description}
\item[(i)]   for each $\alpha<\lambda$ the set $\{i<\sigma:\alpha\in A_i\
\equiv\ \alpha\notin A_\sigma\}$ has cardinality $<\sigma$, and
\item[(ii)]  if $v\subseteq\sigma$, $|v|<\sigma$, $h:v\cup\{\sigma\}
\longrightarrow\{0,1\}$ then there is $\alpha\in\lambda$ such that
\[(\forall i\in v\cup\{\sigma\})(\alpha\in A_i\quad\Leftrightarrow\quad h(i)=
1).\]
\end{description}
For $i\le\sigma$ let $\cI_i=(R,<_{A_i})$. Put $\cI'=\sum\limits_{i<\sigma}
\cI_i$, $\cI''=\cI_\sigma$, $\cI=\cI'+\cI''=\sum\limits_{i\leq\sigma}\cI_i$.
So it is notationally clearer to let $\cI_i=(\{i\}\times R,<_i)$, $(i,y_1)<_i
(i,y_2)$ iff $y_1<_{A_i} y_2$, and $\cI=((\sigma+1)\times R,<^*)$,
$(i_1,y_1)<^*(i_2,y_2)$ iff $i_1<i_2$ or $(i_1=i_2\ \&\ y_1<_{A_{i_1}} y_2)$.

\begin{claim}
\label{cl1}
If $y_0, y_1\in R$, $y_0<_{A_\sigma} y_1$ then that the set $\{i<\sigma: y_0
<_{A_i} y_1\}$ is co-bounded.
\end{claim}

\noindent{\em Proof of the claim:}\hspace{0.2in} Let $t=y_0\wedge y_1$,
$y_0<_{A_\sigma}y_1$. The set
\[u=:\{i<\sigma: c(t)\in A_i\ \equiv\ c(t)\in A_\sigma\}\]
satisfies $u\subseteq\sigma\ \&\ |\sigma\setminus u|<\sigma$. If $c(t)\in
A_\sigma$ then $y_0(\lev(t))<y_1(\lev(t))$. Hence $y_0<_{A_i}y_1$ for $i\in
u$. If $c(t)\notin A_\sigma$ then $y_0(\lev(t))>y_1(\lev(t))$ and $y_1<_{A_i}
y_0$ for $i\in u$. \hfill$\square_{\ref{cl1}}$
\smallskip

Let $\pi:\cI''\longrightarrow\prod\limits_{i<\sigma}\cI'$ be such that for
$y\in\cI_\sigma$, $\pi(y)(i)$ is the element of $\cI_i$ that corresponds to
$y$; recall that all orders $\cI_i$ are defined on $R$. Now, if $D$ is a
uniform ultrafilter on $\sigma$ then $\pi/D:\cI''\to(\cI')^\sigma/D$ is an
embedding. Thus both $(\cI')^\sigma/D$ and $(\cI'')^\sigma/D$ contain a copy
of $\cI''$. Now we will finish by the following claim.

\begin{claim}
\label{cl2}
The linear order $\cI=\cI'+\cI''$ is $\mu$-entangled.
\end{claim}

\noindent{\em Proof of the claim:}\hspace{0.2in} Suppose that $\vare(*)<
\sigma$ and $(j^\alpha_\vare,y^\alpha_\vare)\in\cI$ for $\alpha<\lambda^+$,
$\vare<\vare(*)$ are pairwise distinct and $u\subseteq\vare(*)$. As
$2^\sigma\leq\lambda$ wlog $j^\alpha_\vare=j_\vare$. Let $\{z^\alpha_\zeta:
\zeta<\zeta_\alpha\}$ be an enumeration of $Y_\alpha=:\{y^\alpha_\vare:
\vare<\vare(*)\}$ and let $\{i_\xi:\xi<\xi^*\}$ enumerate $v=:\{j_\vare:
\vare<\vare(*)\}$. Wlog the sequences $\langle z^\alpha_\zeta:\zeta<
\zeta_\alpha\rangle$ are $<_R$-increasing and pairwise disjoint (for
$\alpha<\lambda^+$) as each $z$ may appear in at most $\vare(*)$ of
these sequences. Moreover we may assume that $\{(\zeta_\alpha,\vare,\zeta,\xi):
(j_\vare,y^\alpha_\vare)=(i_\xi,z^\alpha_\zeta)\}$ does not depend on $\alpha$,
so $\zeta_\alpha=\zeta^*$ and by enlarging and renaming instead $\langle(
j^\alpha_\vare,y^\alpha_\vare):\vare<\vare(*)\rangle$ we have $\langle(i_\xi,
z^\alpha_\zeta): \zeta<\zeta^*,\xi<\xi^*\rangle$ and so $u\subseteq\zeta^*
\times\xi^*$. Now, for each $\zeta<\zeta^*$ we choose $c[\zeta]<\lambda$ such
that
\[(\forall \xi<\xi^*)[c[\zeta]\in A_{i_\xi}\quad\Leftrightarrow\quad (\xi,
\zeta)\in u]\]
and we proceed as in earlier cases (considering $h=\langle(c[\zeta],0,1):\zeta
<\zeta^*\rangle\in {}^{\zeta^*}(\lambda\times\lambda\times\lambda)$)\footnote{
Alternatively, let $R$ be the disjoint union of $R_i$ ($i\leq\sigma$) and use
$\cI^*_i=\cI_i\rest(\{i\}\times R_i)$.}. \hfill$\square_{\ref{cl2}}$
\medskip

\noindent 3)\ \ \ The definition of the forcing notion $\Bbb P$ is somewhat
similar to that of the approximations is \ref{2.2}(2).
\begin{sdef}
\label{2.3A}
{\em A condition in $\Bbb P$} is a tuple $p=\langle t,\delta,w,\leq,D,\bar{d},
\bar{e},c\rangle$ (we may write $t^p$, $w^p$, etc) such that:
\begin{description}
\item[(A)] $w\subseteq\lambda^+$ is a set of cardinality $<\lambda$, $\delta$
is a limit ordinal $<\lambda$; let $w^{[\mu]}$ be the family of all increasing
sequences $\bar{y}\subseteq w$ of length $\mu$, 
\item[(B)] $t$ is a non-empty closed under initial segments subset of
${}^{\delta\geq}\lambda$ of cardinality $<\lambda$, 
\item[(C)] $\leq$ is such that $(t\cup w,\leq)$ is a normal tree, $\leq^p\rest
t$ is $\vartriangleleft$, and for $\alpha\in w$, $b^p_\alpha\stackrel{\rm
def}{=}\bigcup\{\nu\in t:\nu\leq^p\alpha\}\in t\cap{}^\delta \lambda$, and
$\alpha\neq\beta\quad\Rightarrow\quad b^p_\alpha\neq b^p_\beta$,
\item[(D)] $c$ is a function from $t$ to $\lambda$,
\item[(E)] $D$ is a set of $<\lambda$ pairs $(\bar{x},h)$ such that
$\bar{x}\in t^{[\mu]}$ and $h\in \bigcup\limits_{\zeta<\mu}{}^\zeta(\lambda
\times\lambda\times\lambda)$, 
\item[(F)] $\bar{d}=\langle d_{\bar{x},h}: (\bar{x},h)\in D\rangle$, each
$d_{\bar{x},h}$ is a partial function from $w^{[\mu]}$ to $\delta^p$ with
domain of cardinality $<\lambda$,
\item[(G)] $\bar{e}=\langle e_{\bar{x},h}: (\bar{x},h)\in D\rangle$, each
$e_{\bar{x},h}$ is a partial function from $t^{[\mu]}_{\bar{x}}\times\delta^p$
to $\{0,2\}$ such that for $(\bar{x},h)\in D$:
\begin{description}
\item[(i)]  $\dom(e_{\bar{x},h})\supseteq\{(\bar{y},\gamma):$ for some
$\bar{z}_0, \bar{z}_1\in\dom(d_{\bar{x},h})$ we have 

\qquad\qquad\qquad $\bar{x}<\bar{y}\leq\bar{z}_0$ and $\gamma=d_{\bar{x},h}(
\bar{z}_1))\}$,
\item[(ii)] if $\bar{y}\in t^{[\mu]}_{\bar{x}}$, $\bar{x}<_t\bar{y}<_t
\bar{z}$ and $e_{\bar{x},h}(\bar{z},\beta)$ is defined then $e_{\bar{x},h}
(\bar{y},\beta)$ is defined and $e_{\bar{x},h}(\bar{y},\beta)\leq
e_{\bar{x},h}(\bar{z},\beta)$,
\item[(iii)] if $\bar{z}\in\dom(d_{\bar{x},h})$ and $\bar{x}<\bar{y}<\bar{z}$
then $e_{\bar{x},h}(\bar{y},d_{\bar{x},h}(\bar{z}))=0$,
\end{description}
and letting $h=\langle (\alpha^0_\vare,\alpha^1_\vare,\alpha^2_\vare): \vare<
\zeta\rangle$ we have
\begin{description}
\item[(iv)] if $(\bar{y},\alpha)\in\dom(e_{\bar{x},h})$, $\lev(\bar{y})=
\delta^p$ and $e_{\bar{x},h}(\bar{y},\alpha)=0$ then 
\[(\forall\vare<\mu)(\exists^\mu i<\mu)(c(u_i)=\vare)\qquad\mbox{and}\]
for $\mu$ ordinals $i<\mu$ divisible by $\zeta$, for some $\xi<\lambda$
\[(\forall\vare<\zeta)(c(y_{i+\vare})=\alpha^0_\vare),\]
\item[(v)]  if $\bar{x}<\bar{y}^0<\bar{y}^1$, $\alpha\leq\lev(\bar{y}^0)+1=
\lev(\bar{y}^1)<\delta^p$ and $e_{\bar{x},h}(\bar{y}^1,\alpha)=0$, then for
$\mu$ ordinals $i<\mu$ divisible by $\zeta$, for some $\xi<\lambda$
\[(\forall\vare<\zeta)(c(y^0_{i+\vare})=\alpha^0_\vare\ \ \&\ \ \zeta\cdot
\xi+\alpha^1_\vare=y^1_{i+\vare}(\lev(\bar{y}^0))),\]
\end{description}
\item[(H)] if $e_{\bar{x},h}(\bar{y},\alpha)=e_{\bar{x},h}(\bar{z},\alpha)=0$,
$\neg[\bar{y}\leq_t\bar{z}]$, $\neg [\bar{z}\leq_t\bar{y}]$ 
then clauses ($\alpha$), $(\beta)$, $(\gamma)$ of \ref{2.3}(4c) hold.
\end{description}
$\Bbb P$ is equipped with the natural partial order:

$p\leq q$\qquad if and only if

$t^p\subseteq t^q$, $\delta^p\leq\delta^q$, $w^p\subseteq w^q$, $c^p\subseteq
c^q$, $D^p\subseteq D^q$,\ \ \ and
\[(\bar{x},h)\in D^p\quad\Rightarrow\quad d^p_{\bar{x},h}\subseteq
d^q_{\bar{x},h}\ \&\ e^p_{\bar{x},h}\subseteq
e^q_{\bar{x},h}.\]
\end{sdef}
For a condition $p\in {\Bbb P}$ and an ordinal $\alpha<\lambda^+$ we define
$q=p\rest\alpha$ by:
\begin{itemize}
\item $\delta^q=\delta^p$, $t^q=t^p$, $w^q=w^p\cap \alpha$, $\leq^q=\leq^p
\rest t^q$, $D^q=D^p$, $c^q=c^p$,
\item if $(\bar{x},h)\in D^q$ then $d^q_{\bar{x},h}=d^p_{\bar{x},h}\rest
(w^q)^{[\mu]}$ and $e^q_{\bar{x},h}=e^p_{\bar{x},h}
\rest ((t^q)_{\bar{x}}^{[\mu]}\times\delta^q)$.
\end{itemize}

\begin{observation}
\label{2.3C}
If $p\in {\Bbb P}$ and $\alpha<\lambda^+$ is either a successor or of
cofinality $\lambda$ then $p\rest\alpha\in\cA\cP$ is the unique maximal
condition such that $p\rest\alpha\leq p$ and $w^{p\rest\alpha}=w^p\cap\alpha$.
\hfill$\square$  
\end{observation}

\begin{claim}[Density Observation]
\label{2.3D}
Assume $p\in {\Bbb P}$.
\begin{enumerate}
\item Suppose $\eta\in{}^{\lambda>}\lambda$. Then there is $q\in {\Bbb P}$
such that
\[p\leq q,\quad \eta\in t^q,\quad w^p=w^q,\quad D^p=D^q,\quad \bar{d}^q=
\bar{d}^p.\]
\item For each $\beta\in\lambda^+\setminus w^p$ there is $q\in {\Bbb P}$
such that $p\leq q$ and $\beta\in w^q$.
\item If $\bar{x}\in(t^p,\leq^p)^{[\mu]}$, $\zeta<\mu$ and $h\in {}^\zeta(
\lambda\times\lambda\times\lambda)$ then there is $q\in {\Bbb P}$ such that
$p\leq q$ and $D^q=D^p\cup\{(\bar{x},h)\}$. 
\item If $(\bar{x},h)\in D^p$, $\bar{x}<\bar{y}\in (w^p)^{[\mu]}$ then
there is $q\in\cA\cP$, $p\le q$ such that $w^q=w^p$ and $\bar{y}\in
\dom(d_{\bar{x},h})$. 
\end{enumerate}
\end{claim}

\noindent{\em Proof of the claim:}\hspace{0.2in} 1)\ \ \ We may assume that
$|t^p\cap {}^{\delta^p}\lambda|=\mu$, as otherwise it is even easier. So let
$\langle\nu_\beta:\beta<\mu\rangle$ enumerate $t^p\cap{}^{\delta^p}\lambda$.
Put $\delta^q=\max\{\delta^p+\omega,\lh(\eta)+\omega\}$ and fix $\rho\in
{}^{\delta^q}\lambda$ such that $\eta\vartriangleleft\rho$. Next, by induction
on $\gamma\in (\delta^p,\delta^q]$ define sequences $\nu^\gamma_\beta$ for
$\beta<\mu$ and functions $c_\gamma$ such that
\begin{description}
\item[(a)] $c_\gamma:\{\nu^\gamma_\beta:\beta<\mu\}\longrightarrow\lambda$, 
\quad $\nu^\gamma_\beta\in {}^\gamma\lambda$.
\item[(b)] $\delta^p<\gamma_0<\gamma_1\leq\delta^q\quad\Rightarrow\quad
\nu_\beta\vartriangleleft\nu^{\gamma_0}_\beta\vartriangleleft
\nu^{\gamma_1}_\beta$.
\item[(c)] Suppose that $(\bar{x},h)\in D^p$, $h=\langle(\alpha^0_\vare,
\alpha^1_\vare,\alpha^2_\vare):\vare<\zeta\rangle\in{}^\zeta(\lambda\times
\lambda\times\lambda)$, $\bar{x}<\bar{y}=\langle\nu_{\vare(i)}:i<\mu\rangle
\subseteq t^p\cap {}^{\delta^p}\lambda$, $(\bar{y},\alpha)\in\dom(e^p_{\bar{x},
h})$ and $e_{\bar{x},h}^p(\bar{y},\alpha)=0$. Then 
\begin{description}
\item[$(\alpha)$] for $\mu$ ordinals $i<\mu$ divisible by $\zeta$, for some
$\xi<\lambda$ 
\[(\forall j<\zeta)(c^p(\nu_{\vare(i+j)})=\alpha^0_j\ \&\ \nu^{\delta^p+1}_{
\vare(i+j)}(\delta^p)=\zeta\cdot\xi+\alpha^1_j),\]
\item[$(\beta)$]  for each $\gamma\in (\delta^p,\delta^q)$ we have
\[(\forall\xi<\mu)(\exists^\mu i<\mu)(c_\gamma(\nu^\gamma_{\vare(i)})=\xi)
\quad\mbox{ and}\]
for $\mu$ ordinals $i<\mu$ divisible by $\zeta$ there is $\xi<\lambda$ such
that 
\[(\forall j<\zeta)(c_\gamma(\nu^\gamma_{\vare(i+j)})=\alpha^0_j\ \&\
\nu^{\gamma+1}_{\vare(i+j)}=\zeta\cdot\xi+\alpha^1_j),\]
\item[$(\gamma)$] $(\forall\xi<\mu)(\exists^\mu i<\mu)(c_{\delta^q}(
\nu^{\delta^q}_{\vare(i)})=\xi)$\quad and for $\mu$ ordinals $i<\mu$ divisible
by $\zeta$
\[(\forall j<\zeta)(c_{\delta^q}(\nu^{\delta^q}_{\vare(i+j)}=\alpha^0_j).\]
\end{description}
\end{description}
The construction is easy and can be done like the one in \ref{2.2G}. Next we
put
\[\begin{array}{l}
t^q=t^p\cup\{\rho\rest\gamma:\gamma\leq\delta^q\}\cup\{\nu^\gamma_\beta:
\beta<\mu, \delta^p<\gamma\leq\delta^q\},\\
w^q=w^p,\quad D^q=D^p,\quad \bar{d}^q=\bar{d}^p.
\end{array}\]
The function $c^q$ is defined as any extension of $c^p\cup\{c_\gamma: \delta^p
<\gamma\leq\delta^q\}$ (note that the possibly non-defined values are that at
some initial segments of $\rho$). For $\alpha\in w^q$ let $b^q_\alpha$ be
$\nu^{\delta^q}_\beta$ for beta such that $b^p_\alpha=\nu_\beta$. This
determines the tree ordering $\leq^q$ of $t^q\cup w^q$ (remember
\ref{2.3A}(C)). To define $e^q_{\bar{x},h}$ we let 
\[\begin{array}{ll}
e^q_{\bar{x},h}\supseteq e^p_{\bar{x},h}, &\ \\
\dom(e^q_{\bar{x},h})=\dom(e^p_{\bar{x},h})\cup\{(\bar{y},\iota):&
\bar{y}=\langle\nu^\gamma_{\beta(i)}:i<\mu\rangle>\bar{y}'=\langle\nu_{\beta
(i)}:i<\mu\rangle,\\
\ &\gamma\in (\delta^p,\delta^q],\ (\bar{y}',\iota\in \dom(e^p_{\bar{x},h})\
\},\\
\mbox{if }\bar{x}<\bar{y}'<\bar{y},\ (\bar{y}',\iota)\in\dom(e^p_{\bar{x},h}),
& \lev(\bar{y}')=\delta^p\mbox{ then }e^q_{\bar{x},h}(\bar{y},\iota)=e^p_{
\bar{x},h}(\bar{y}',\iota).
\end{array}\]
Now one easily checks that the condition $q$ defined above is as required. 
\medskip

\noindent 2)\ \ \ Choose $\nu\in {}^{\delta^p}\lambda$ such that
$\nu\rest 1\notin t^p$. Let 
\[\delta^q=\delta^p,\quad t^q=t^p\cup\{\nu\rest\gamma:\gamma<\delta^p\}\quad
\mbox{ and}\quad w^q=w^p\cup\{\beta\}.\]
Define $\leq^q$ by letting $b^q_\alpha=b^p_\alpha$ for $\alpha\in w^p$ and
$b^q_\beta=\{\nu\rest\vare:\vare<\delta^q\}$. Finally put $D^q=D^p$,
$d^q_{\bar{x},h}=d^p_{\bar{x},h}$ and $e^q_{\bar{x},h}=e^q_{\bar{x},h}$.
\medskip

\noindent 3)\ \ \ Let $t^q=t^p$, $w^q=w^p$, $c^q=c^p$, $d^q_{\bar{x}',h'}=
d^p_{\bar{x}',h'}$ if $(\bar{x}',h')\in D^p$, $d^q_{\bar{x},h}$ is empty if
$(\bar{x},h)\notin D^p$, and similarly for $e^q_{\bar{x}',h'}$.
\medskip

\noindent 4)\ \ \ Let $h=\langle(\alpha^0_\vare,\alpha^1_\vare,\alpha^2_\vare)
:\vare<\zeta\rangle$. Declare that $\delta^q=\delta^p+\omega$ and fix an
enumeration $\langle\nu_\beta:\beta<\mu\rangle$ of $t^p\cap{}^{\delta^p}
\lambda$. Let $\bar{y}'=\langle\nu_{\beta(i)}: i<\mu\rangle$ be such that
$\bar{y}'<\bar{y}$ ($\in w^{[\mu]}$). Next, like in (1) above build
$c_\gamma,\nu^\gamma_\beta$ for $\gamma\in (\delta^p,\delta^q]$, $\beta<\mu$
satisfying demands (a)--(c) there plus:
\begin{description}
\item[(c)]
\begin{description}
\item[$(\delta)$] for each $\gamma\in (\delta^p,\delta^q)$
\[(\forall\xi<\mu)(\exists^\mu i<\mu)(c_\gamma(\nu^\gamma_{\beta(i)})=
\xi)\quad\mbox{ and}\]
for $\mu$ ordinals $i<\mu$ divisible by $\zeta$ there is $\xi<\lambda$ such
that 
\[(\forall j<\zeta)(c_\gamma(\nu^\gamma_{\beta(i+j)})=\alpha^0_j\ \&\
\nu^{\gamma+1}_{\beta(i+j)}(\gamma)=\zeta\cdot\xi+\alpha^1_j),\]
\item[$(\vare)$]  $(\forall\xi<\mu)(\exists^\mu i<\mu)(c_{\delta^q}(
\nu^{\delta^q}_{\beta(i)})=\xi)$\quad and for $\mu$ ordinals $i<\mu$ divisible
by $\zeta$ 
\[(\forall j<\zeta)(c_{\delta^q}(\nu^{\delta^q}_{\beta(i+j)}=\alpha^0_j).\]
\end{description}
\end{description}
Next we put $\dom(d^q_{\bar{x},h})=\dom(d^p_{\bar{x},h})\cup\{\bar{y}\}$,
$d^q_{\bar{x},h}\supseteq d^p_{\bar{x},h}$, $d^q_{\bar{x},h}(\bar{y}=\delta^p
+2$, $d^q_{\bar{x}',h'}=d^p_{\bar{x}',h'}$ for $(\bar{x}',h')\in D^p\setminus
\{(\bar{x},h)\}$, and $D^q=D^p$. The functions $c^q$ and $e^q_{\bar{x}',h'}$
(for $(\bar{x}',h')\in D^q$) are defined as in (1), but dealing with
$(\bar{x},h)$ we take into account the new obligation: $d^q_{\bar{x},h}(
\bar{y})=\delta^p+2$ (and we put the value 2 whenever possible). There is no
problem with it as we demanded clauses (c)$(\delta,\vare)$. Now one easily
checks that $q$ is as required. \hfill$\square_{\ref{2.3D}}$
\medskip

\begin{claim}
\label{2.3B}
The forcing notion $\Bbb P$ is $(<\lambda)$--complete, i.e.~if $\bar{p}=
\langle p^i:i<i^*\rangle\subseteq {\Bbb P}$ is increasing, $i^*<\lambda$ then
$\bar{p}$ has an upper bound in ${\Bbb P}$.  
\end{claim}

\noindent{\em Proof of the claim:}\hspace{0.2in} It is easy: the first
candidate for the upper bound is the natural union of the $p^i$'s. What may
fail is that the tree $\bigcup\limits_{i<i^*} t^{p^i}$ does not have the last
level. But this is not a problem as we may use the procedure of \ref{2.3D}(1)
to add it. \hfill$\square$

\begin{claim}[The Amalgamation Property]
\label{2.3E}
If $\alpha<\lambda^+$ is either a successor ordinal or of cofinality
$\lambda$, $p,q\in {\Bbb P}$, $p\rest\alpha\leq q$ and $w^q\subseteq\alpha$ 
then there is $r\in {\Bbb P}$ such that $p\le r$, $q\leq r$ and $w^r=w^p\cup
w^q$. 
\end{claim}

\noindent{\em Proof of the claim:}\hspace{0.2in} By \ref{2.3D}(1) we may
assume that $\delta^p<\delta^q$. Moreover we may assume that $|w^p\setminus
\alpha|=\mu$ (as otherwise everything is easier). Let $\delta^r=\delta^q$ and
$w^r=w^p\cup w^q$. By induction on $\gamma\in [\delta^p,\delta^q]$ choose
sequences $\langle\nu_{\beta,\gamma}:\beta\in w^p\setminus\alpha\rangle$
and functions $c_\gamma$ such that
\begin{description}
\item[$(\alpha)$] $\nu_{\beta,\gamma}\in{}^\gamma\lambda$ are 
$\vartriangleleft$--increasing with $\gamma$,
\item[$(\beta)$]  $\nu_{\beta,\delta^p}=b^p_\beta$,\quad $\nu_{\beta,\delta^p
+1}\notin t^q$

[note that $\langle\nu_{\beta,\delta^p}:\beta\in w^p\setminus\alpha\rangle$ is
with no repetition],
\item[$(\gamma)$] $c_\gamma:\{\nu_{\beta,\xi}:\beta\in w^p\setminus\alpha,\
\xi\in [\delta^p,\gamma)\}\longrightarrow\lambda$ are continuously increasing
with $\gamma$, $c_{\delta^p+1}$ is $c^p$ restricted to $\{\nu_{\beta,\delta^p}
: \beta\in w^p\setminus\alpha\}$,
\end{description}
and for each $(\bar{x},h)\in D^p$, $h=\langle(\alpha^0_\vare,\alpha^1_\vare,
\alpha^2_\vare):\vare<\zeta\rangle$, $\bar{z}\in\dom(d^p_{\bar{x},h})\setminus
\dom(d^q_{\bar{x},h})$ and $i^*<\mu$ such that $z_i\geq\alpha$ for $i\geq i^*$
we have
\begin{description}
\item[$(\delta)$] for each $\gamma\in [\delta^p,\delta^q)$, for $\mu$ ordinals
$i\in [i^*,\mu)$ divisible by $\zeta$, for some $\xi<\lambda$
\[(\forall\vare<\zeta)(c_{\gamma+1}(\nu_{z_{i+\vare},\gamma})=\alpha^0_\vare\
\&\ \nu_{z_{i+\vare},\gamma+1}(\gamma)=\zeta\cdot\xi+\alpha^1_\vare),\]
\item[$(\vare)$] for $\mu$ ordinals $i\in [i^*,\mu)$ divisible by $\zeta$ 
\[(\forall\vare<\zeta)(c_{\delta^q+1}(\nu_{z_{i+\vare},\delta^q})=
\alpha^0_\vare),\]
\item[$(\zeta)$] if $\bar{x}<\bar{y}<\bar{z}$, $\lev(\bar{y})=\delta^p$,
$\bar{y}<\bar{y}'$, $\lev(\bar{y}')=\delta^p+1$, $(\bar{y}',d^p_{\bar{x},h}(
\bar{z}))\in\dom(e^q_{\bar{x},h})$ and $e^q_{\bar{x},h}(\bar{y}',d^p_{\bar{x},
h}(\bar{z}))=0$ then for $\mu$ ordinals $i\in [i^*,\mu)$ divisible by $\zeta$
there are $\xi_0<\xi_1<\lambda$ such that
\[(\forall\vare<\zeta)(\zeta\cdot\xi_0\leq y^\prime_{i+\vare}(\delta^p)<
\zeta\cdot\xi_1\leq y_{i+\vare}(\delta^p))\quad\mbox{ and}\]
\[(\forall\vare<\zeta)(c^p(\nu_{z_{i+\vare},\delta^p})=\alpha^0_\vare,\ \
y^\prime_{i+\vare}(\delta^p)=\zeta\cdot\xi_0+\alpha^1_\vare,\ \ y_{i+\vare}(
\delta^p)=\zeta\cdot\xi_1+\alpha^2_\vare,\]
\item[$(\iota)$]  for each $\gamma\in (\delta^p,\delta^q]$, for every
$\vare<\mu$ there are $\mu$ ordinals $i<\mu$ such that $c_{\gamma+1}(\nu_{
z_i,\gamma})=\vare$. 
\end{description}
[Our intension here is that $b^r_\beta=\nu_{\beta,\delta^q}$ and $c^r\supseteq
c_{\delta^q}$.] We have actually $\mu$ demands, each of which can be satisfied
by $\mu$ pairwise disjoint cases of size $\zeta<\mu$. So we may carry out the
procedure analogous to that from the end of the proof of \ref{2.2H}. Note that
in handling instances of clause $(\zeta)$ we use demand \ref{2.3A}(G)(v) for
$q$ (applicable as $d^p_{\bar{x},h}(\bar{z})<\delta^p$) and for clause
$(\delta)$ we use \ref{2.3A}(G)(iv). After the construction is carried out we
easily define a condition $r$ as required. \hfill$\square_{\ref{2.3E}}$ 

\begin{claim}
\label{cc}
The forcing notion $\Bbb P$ satisfies the $\lambda^+$--cc.
\end{claim}

\noindent{\em Proof of the claim:}\hspace{0.2in} Suppose that $\langle
p_\alpha:\alpha<\lambda^+\rangle$ is an antichain in $\Bbb P$. By passing to a
subsequence we may assume that $\otp(w^{p_\alpha})$ is constant and that the
order isomorphism of $w^{p_\alpha}$, $w^{p_\beta}$ carries the condition
$p_\alpha$ to $p_\beta$ (so $t^{p_\alpha}=t^{p_\beta}$,
$D^{p_\alpha}=D^{p_\beta}$ etc). Moreover, we may assume that the family $\{
w^{p_\alpha}:\alpha<\lambda^+\}$ forms a $\Delta$--system with kernel $w^*$
(remember $\lambda=2^\mu=\mu^+$). Now we may find an ordinal $\alpha^*<
\lambda^+$ of cofinality $\lambda$ and $\alpha_0<\alpha_1<\lambda^+$ such that
$w^{p_{\alpha_0}}\subseteq \alpha^*$, $w^{p_{\alpha_1}}\cap\alpha^* =w^*$ and
$w^*$ is an initial segment of both $w^{p_{\alpha_0}}$ and
$w_{p_{\alpha_1}}$. Note that then $p_{\alpha_1}\rest \alpha^*\geq
p_{\alpha_0}$. So applying \ref{2.3E} we conclude that the conditions
$p_{\alpha_0}, p_{\alpha_1}$ have a common upper bound, a contradiction. 
\hfill$\square_{\ref{cc}}$
\medskip

To finish the proof note that if $G\subseteq {\Bbb P}$ is a generic filter
over $\V$ then, in $\V[G]$ we may define the tree $T$ by:
\[T=({}^{\lambda>}\lambda)\cup\{\eta_\alpha:\alpha<\lambda^+\}\]
where for $\alpha<\lambda^+$ we define $\eta_\alpha\in {}^\lambda\lambda$ by
\[\eta_\alpha=\bigcup\{\nu\in{}^{\lambda>}\lambda: (\exists p\in
G)(\nu\in b^p_\alpha)\}\]
(and $c,d$ are defined similarly). By \ref{2.3D} and \ref{2.3B} (no new
$\mu$--sequences of ordinals are added) we easily conclude that these objects
are as needed. \QED$_{\ref{2.3}}$ 
\medskip

We may want to improve \ref{2.2}(2) so that it looks more like \ref{2.3}(4);
we can do it at some price.

\begin{proposition}
\label{2.4}
Let $\cJ^*$ be a linear order, $\cJ^*=\sum\limits_{\alpha<\lambda^+}
\cI_\alpha$, each $\cI_\alpha$ a $\lambda$-dense linear order of cardinality
$\lambda$ (as in the proof of \ref{2.2}(1)). Then $\omega\times\cJ^*$ is a
$\lambda^+$-like linear order such that every $j\in\cI$ which is neither
successor nor the first element (under $<_{\cI}$) satisfies
\[\cf(\{i:i<_{\cJ}j\},<_{\cJ})=\lambda\]
and each member of $\omega\times\cJ$ has an immediate successor.\QED
\end{proposition}

\begin{definition}
\label{2.5}
For a $\lambda^+$-like linear order $\cJ$, a $\cJ$-Aronszajn tree is
$T=(T,\leq_T,\lev)$ such that
\begin{description}
\item[(a)] $T$ is a set of cardinality $\lambda^+$,
\item[(b)] $(T,\leq_T)$ is a partial order which is a tree, i.e. for every
$y\in T$ the set $\{x: x\leq_T y\}$ is linearly ordered by $\leq_T$,
\item[(c)] $\lev$ is a function from $T$ to $\cJ$, $T_j=:\{y\in T:\lev(y)=
j\}$,
\item[(d)] for every $y\in T$, $\lev$ is a one-to-one order preserving
function from $\{x: x<_T y\}$ onto $\{j\in \cJ: j<_{\cJ}\lev(y)\}$, so
$y\rest j$ is naturally defined,
\item[(e)] for $y\in T$ and $j\in\cJ$, $\lev(y)<_{\cJ} j$ there is $z$ such
that $y<_T z\in T$, $\lev(z)=j$,
\item[(f)] $\{y: \lev(y)=j\}$ has cardinality $\lambda$,
\item[(g)] {\em normality:}\qquad if $y\neq z$, both in $T_j$, $j$ is neither
successor nor the first element of $(\cJ,<_{\cJ})$ then $\{x: x<_T y\}\neq
\{x: x<_T z\}$,
\item[(h)] for $y\neq z\in T$ there is $j\in\cJ$ such that $y\rest j= z
\rest j$ and
\[(\forall i)(j<_{\cJ} i\quad\Rightarrow\quad y\rest i\neq z\rest i)\]
[we write $y\rest j= z\rest j$ for $z\wedge y$].
\end{description}
\end{definition}

\begin{theorem}
\label{2.6}
Assume that $\lambda=\mu^+=2^\mu$ and $\diamondsuit_\lambda$ (the second
follows e.g.~if $\mu\geq \beth_\omega$ -- see \cite[3.5(1)]{Sh:460}) and
$\cJ$ is as constructed in \ref{2.4}. Then there are a $\cJ$-Aronszajn
tree $T$ and functions $f$, $c$, $d$ such that
\begin{description}
\item[(a)] $f$, $c$ are functions from $T$ to $\lambda$, if $y$ is the
successor of $j$ in $\cJ$, $y\in T_i$ then $f$ is one to one from
$\{z\in T_j: z\restriction i=y\}$ onto $\lambda$
\item[(b)] for every $\bar x\in T^{[\mu]}$ ($=\bigcup\limits_{j\in {\cJ}}
T_j^{[\mu]}$) and $h\in H^\otimes_{\mu,\zeta}$, $\zeta<\mu$ we have
$d_{\bar{x},h}:T^{[\mu]}_{\bar{x}}\longrightarrow\lambda$ such that:\\
if $d_{\bar{x},h}(\bar{x})=d_{\bar{x},h}(\bar{y})$ then for some $\bar{t}
\in T^{[\mu]}_{\bar{x}}$ we have
\begin{description}
\item[$(\alpha)$] $t_i=y_i\wedge z_i$,
\item[$(\beta)$]  $\lev(\bar{t})<\lev(\bar{y})$, $\lev(\bar{t})<\lev(
\bar{z})$,
\item[$(\gamma)$] $(\forall\vare<\mu)(\exists^\mu i<\mu)(c(t_i)=\vare)$,
\item[$(\delta)$] for $\mu$ ordinals $i<\mu$ divisible by $\zeta$ we
have
\[h(\zeta)=\langle(c({t}_i),f(y_i\rest (\alpha+_{\cJ}1)),
f(z_i\rest(\alpha+_{\cJ}1))),\]
where $\alpha=\lev(\bar{t})$.
\end{description}
\end{description}
\end{theorem}

\Proof Like \ref{2.2}. \QED

\section{Constructions Related to pcf Theory}

\begin{lemma}
\label{3.1}
\begin{enumerate}
\item Suppose that
\begin{description}
\item[(A)] $\langle\lambda_i:i<\delta\rangle$ is a strictly increasing
sequence of regular cardinals, $|\delta|<\lambda_i<\lambda=\cf(\lambda)$
for $i<\delta$ and $D$ is a $\sigma$-complete filter on $\delta$ containing
all co-bounded subsets of $\delta$ (follows by clause (D); hence
$\cf(\delta)\geq\sigma$),
\item[(B)] $\tcf(\prod\limits_{i<\delta}\lambda_i/D)=\lambda$, i.e.~there is
a sequence $\langle f_\alpha:\alpha<\lambda\rangle\subseteq\prod\limits_{i<
\delta}\lambda_i$ such that
\begin{description}
\item[(i)]  $\alpha<\beta<\lambda$ implies  $f_\alpha<_D f_\beta$,
\item[(ii)] $(\forall f\in\prod\limits_{i<\delta}\lambda_i)(\exists
\alpha<\lambda)(f<_D f_\alpha)$,
\end{description}
\item[(C)] sets $A_i\subseteq\delta$ (for $i<\kappa$) are such that the
family $\{A_i: i<\kappa\}$ is $\sigma$-independent in $\cP(\delta)/D$
(i.e.~if $u,v$ are disjoint subsets of $\kappa$, $|u\cup v|<\sigma$ then
$\bigcap\limits_{i\in u} A_i\setminus\bigcup\limits_{j\in v} A_j\neq
\emptyset\ \mod\;D$),
\item[(D)] $|\{f_\alpha\rest i:\alpha<\lambda\}|^{<\sigma}<\lambda_i$
for each $i<\delta$.
\end{description}
Then  $\ens_\sigma(\lambda,\kappa)$.
\item The linear orders in part (1) have exact density $\mu=:\sum\limits_{i<
\delta}\lambda_i$ (see Definition \ref{1.9}) and they are positively
$\mu$-entangled (see Definition \ref{6.4}). Moreover, if $\langle f_\alpha:
\alpha<\lambda\rangle$ is as gotten in \cite[\S 1]{Sh:355} (i.e.~it is
$\mu$-free) then they have exact density $(\mu^+,\mu^+,\mu)$.
\end{enumerate}
\end{lemma}

\noindent{\bf Remark:}\qquad By \cite[3.5]{Sh:355}, if $\delta<\lambda_0$ and
$\max\pcf(\{\lambda_i: i< \alpha\})<\lambda_\alpha$ (for $\alpha<\delta$) then
we can have (D); i.e.~we can find $f_\alpha$ (for $\alpha<\lambda$) satisfying
(i) + (ii) of (B) and (D). If in addition $\lambda\in\pcf_{\sigma{\rm
-complete}}\{\lambda_i: i<\delta\}$ (but $\lambda\notin\pcf(\{\lambda_i:i<
\alpha\})$ for $\alpha<\lambda$) then we can find a filter $D$ as required in
(A) + (B). So if $\mu>\cf(\mu)=\kappa$ and $\alpha<\mu\quad\Rightarrow\quad
|\alpha|^\kappa<\mu$, then we can find $\bar{\lambda}=\langle\lambda_i: i<
\kappa\rangle$ as above and $\bar{\lambda}$ strictly increasing with limit
$\mu$. 
\medskip

\Proof  1)\ \ \ Let $\cI=\{f_\alpha:\alpha<\lambda\}$. For each
$\zeta<\kappa$ we define a linear order $<^*_\zeta$ of $\cI$. It is
$<^*_{A_\zeta}$, where for $A\subseteq\delta$ we define $<^*_A$ by:

$f_\alpha<^*_A f_\beta$\qquad if and only if
\[(\exists i<\delta)(f_\alpha(i)\neq f_\beta(i)\ \&\ f_\alpha\rest i=f_\beta
\rest i\ \&\ [f_\alpha(i)<f_\beta(i)\quad\Leftrightarrow\quad i\in
A]).\]
Let $u,v$ be disjoint subsets of $\kappa$, $|u\cup v|<\sigma$ and for each
$\vare\in u\cup v$ let $t^\vare_\alpha=f_{\gamma(\vare,\alpha)}$ be pairwise
distinct (for $\alpha<\lambda$). We should find $\alpha<\beta<\lambda$
as in \ref{1.1}(1). Let
\[\begin{array}{ll}
g_\alpha(i)=:&\min\{f_{\gamma(\vare,\alpha)}(i): \vare\in u\cup v\},\\
i_\alpha=:   &\min\{i<\delta:\langle f_{\gamma(\vare,\alpha)}\rest
i:\vare\in u\cup v\rangle\ \mbox{ are pairwise distinct}\}.
\end{array}\]
Since $|u\cup v|<\sigma\leq\cf(\delta)$, $i_\alpha<\delta$. Clearly
$g_\alpha\in \prod\limits_{i<\delta}\lambda_i$. Without loss of generality
$i_\alpha=i^*$ for every $\alpha<\lambda$. Let
\[B=\{i<\delta:(\forall \xi<\lambda_i)(\exists^\lambda\alpha<\lambda)(
g_\alpha(i)>\xi)\}.\]

\begin{claim}
\label{3.1A}
$B\in D$.
\end{claim}

\noindent{\em Proof of the claim:}\hspace{0.2in} Assume not, so $\delta
\setminus B\neq\emptyset\ \mod\,D$. For $i\in\delta\setminus B$ let $\xi_i<
\lambda_i$, $\beta_i<\lambda$ exemplify $i\notin B$, i.e. $\alpha\in
[\beta_i,\lambda)\Rightarrow g_\alpha(i)\leq\xi_i$. Define $h\in\prod\limits_{
i<\delta}\lambda_i$ by:
\[h(i)=:\left\{\begin{array}{ll}
\xi_i+1 &\mbox{if }\ i\in\delta\setminus B;\\
0       &\mbox{if }\ i\in B.
\end{array}\right.\]
Now $\langle f_\alpha/D:\alpha<\lambda\rangle$ is cofinal in $\prod\limits_{i
<\delta} \lambda_i/D$ (i.e. clause (ii) of (B)), so there exists $\beta<
\lambda$ such that $h<f_\beta\ \mod\,D$. Without loss of generality $\sup
\limits_{i\in\delta\setminus B}\beta_i<\beta$ (remember that $\delta\setminus
B\subseteq\delta$, $|\delta|<\lambda=\cf(\lambda)$ and $(\forall i\in\delta
\setminus B)(\beta_i<\lambda)$). Since, for each $\vare\in u\cup v$,
$\gamma(\vare,\alpha)$ (for $\alpha<\lambda$) are pairwise distinct and
$\beta<\lambda$, there exists $\alpha<\lambda$ such that $(\forall\vare\in
u\cup v)(\gamma(\vare,\alpha)>\beta)$. Without loss of generality $\beta<
\alpha$ and hence $\sup\limits_{i\in\delta\setminus B}\beta_i<\alpha$. Now by
the choice of $\alpha$ we have $(\forall\vare\in u\cup v)(f_\beta< f_{\gamma(
\vare,\alpha)}\ \mod\,D)$ and for every $i\in\delta\setminus B$, $g_\alpha(i)
\leq\xi_i$. Hence $E_\vare=:\{i<\delta: f_\beta(i)<f_{\gamma(\vare,\alpha)}
(i)\}\in D$ and as $D$ is $\sigma$-complete and $\sigma>|u\cup v|$ we get
$\bigcap\limits_{\vare\in u\cup v} E_\vare\in D$. By $g_\alpha$'s definition
and the choice of $\beta$, it now follows that $\{i<\delta:h(i)<g_\alpha(i)\}
\in D$ and thus 
\[(\delta\setminus B)\cap\{i<\delta: h(i)<g_\alpha(i)\}\neq \emptyset\ \mod\,
D.\]
Choosing $i$ in this (non-empty) intersection, one obtains $g_\alpha(i)\leq
\xi_i<\xi_i+1=h(i)<g_\alpha(i)$ (the first inequality -- see above, the third
equality -- see choice of $h$, the last inequality -- see choice of
$i$), a contradiction. So $B\in D$, proving the claim.
\hfill$\square_{\ref{3.1A}}$
\medskip

Remember that $|\{f_\alpha\rest i:\alpha<\lambda\}|<\lambda_i$ for each $i<
\delta$, and $\cf(\prod\limits_{i<\delta}\lambda_i/D)=\lambda$, $D$ contains
all co-bounded subsets of $\delta$. By our hypothesis,
\[A=:\bigcap\limits_{\vare\in u} A_\vare\cap\bigcap\limits_{\vare\in v}(
\delta\setminus A_\vare)\neq\emptyset\ \mod\, D,\]
so $C=:\{i<\delta: i^*<i\}\cap A\cap B\neq\emptyset\ \mod\,D$, and one can
choose $i\in C$. For each $\xi<\lambda_i$ choose $\alpha_\xi<\lambda$ such
that $g_{\alpha_{\xi}}(i)>\xi$. Then easily for some unbounded $S\subseteq
\lambda_i$ we have:
\[\xi_1<\xi_2\in S\ \& \ \vare_1,\vare_2\in u\cup v\quad\Rightarrow\quad
f_{\gamma(\vare_1,\alpha_{\xi_1})}(i)<f_{\gamma(\vare_2,\alpha_{\xi_2})}(i).\]
Without loss of generality the sequence $\big\langle\langle f_{\gamma(\vare,
\alpha_{\xi})}\rest i:\vare\in u\cup v\rangle:\xi\in S\big\rangle$ is constant
(by hypothesis (D) of \ref{3.1}(1)). The conclusion should be clear now (look
at the definition of $<_\zeta^*$ and the choice of $i$ being in
$\bigcap\limits_{\vare\in u}A_\vare\setminus\bigcup\limits_{\vare\in v}
A_\vare$). 
\medskip

\noindent 2)\ \ \ We will state the requirements and prove them one by
one.
\begin{claim}
\label{cl3}
The linear orders constructed in the first part have exact density $\mu$.
\end{claim}

\noindent{\em Proof of the claim:}\hspace{0.2in} Let us consider
$\cI=(\cI,<_A)$. For each $i<\delta$ choose a set $X_i\subseteq\lambda$
such that $|X_i|\leq\lambda_i$ and $\{f_\alpha\rest i: \alpha<\lambda\}=\{
f_\alpha\rest i:\alpha\in X_i\}$. Then $\{f_\alpha:\alpha\in\bigcup\limits_{
i<\delta} X_i\}$ is a dense subset of $(\cI,<_A)$ (and its size is
$\leq\mu$).\\
Suppose now that $\cJ\subseteq\cI$, $|\cJ|=\lambda$ and assume that
$\cJ_0\subseteq\cJ$ is a dense subset of $\cJ$, $|\cJ_0|<\mu$. The set $X=
\{\alpha<\lambda:f_\alpha\in\cJ\}$ has cardinality $\lambda$, so it is
unbounded in $\lambda$. Let $i(*)=\min\{i<\delta:\lambda_i>|\cJ_0|\}$.
Then
\[(\forall i\geq i(*))(\gamma_i=:\sup\{f_\alpha(i)+1:f_\alpha\in\cJ_0\}<
\lambda_i)\]
(as $\gamma_i$ is a supremum of a set of $|\cJ_0|<\lambda_i=\cf(\lambda_i)$
ordinals $<\lambda_i$). Let $\gamma_i=0$ for $i<i(*)$. Then $\langle\gamma_i:
i<\delta\rangle\in\prod\limits_{i<\delta} \lambda_i$ and, as $\langle f_\alpha
:\alpha<\lambda\rangle$ is cofinal in $(\prod\limits_{i<\delta}\lambda_i,
<_D)$, for some $\alpha(*)<\lambda$ we have $\langle\gamma_i:i<\delta\rangle
<_D f_{\alpha(*)}$. Since $(\forall\alpha\in X\setminus\alpha(*))(f_{\alpha
(*)}\leq_D f_\alpha)$, we have
\[(\forall\alpha\in X\setminus\alpha(*))(\{i<\delta:\gamma_i<f_\alpha(i)\}\in
D).\]
Consequently, for each $\alpha\in X\setminus \alpha(*)$ we find $i_\alpha\in
(i(*),\delta)$ such that $f_\alpha(i_\alpha)>\gamma_{i_\alpha}$. As $\lambda=
\cf(\lambda)>|\delta|$, there is $j\in (i(*),\delta)$ such that the set
\[X'=:\{\alpha\in X: \alpha>\alpha(*)\ \&\ i_\alpha=j\}\]
is unbounded in $\lambda$. Since $|\{f_\alpha\rest j:\alpha<\lambda\}|<
\lambda_j<\lambda=\cf(\lambda)$, for some unbounded set $X''\subseteq
X'$ and a sequence $\nu$ we have $(\forall\alpha\in X'')(f_\alpha\rest
j=\nu)$. But now note that the convex hull of $\{f_\alpha:\alpha\in
X''\}$ in $(\cI,<_A)$ is disjoint from $\cJ_0$, a contradiction.
\hfill$\square_{\ref{cl3}}$

\begin{claim}
\label{cl4}
$(\cI,<_A)$ is positively $\sigma$-entangled.
\end{claim}

\noindent{\em Proof of the claim:}\hspace{0.2in}  Like in part (1).

\begin{claim}
\label{cl5}
If the sequence $\langle f_\alpha:\alpha<\lambda\rangle$ is $\mu$--free
and the set $A\subseteq\delta$ is neither bounded nor co-bounded then the
linear order $(\cI,<_A)$ has exact density $(\mu^+,\mu^+,\mu)$.
\end{claim}

\noindent{\em Proof of the claim:}\hspace{0.2in} Suppose that
$\cJ\subseteq\cI$ is of size $\geq\mu^+$. By \ref{cl3} its density is
$\leq\mu$. For the other inequality suppose that $\cJ_0$ is a dense
subset of $\cJ$ of cardinality $<\mu$. Let
\[\begin{array}{ll}
\cJ'=\{f_\alpha\in\cJ: &\mbox{for each }i<\delta\mbox{ there are }\beta_1,
\beta_2\mbox{ such that }f_{\beta_1},f_{\beta_2}\in\cJ\mbox{ and}\\
\ &f_{\beta_1}\rest i=f_\alpha\rest i= f_{\beta_2}\rest i\ \&\ f_{\beta_1}
<_A f_\alpha<_A f_{\beta_2})\}.
\end{array}\]
Plainly $|\cJ\setminus\cJ'|=\mu$, so $|\cJ'|\geq\mu^+$. Since $\theta=:
|\cJ_0|<\mu$ and $\mu$ is a limit cardinal, we have $\sigma=(\theta+|\delta|
)^+<\mu$. Let $X=\{\alpha: f_\alpha\in\cJ'\}$ and choose $X_1\subseteq X$ of
size $\sigma$. Now we may find $\langle B_\alpha:\alpha\in X_1\rangle\subseteq
D$ such that for each $j<\delta$ the sequence $\langle f_\alpha(j):\alpha\in
X\ \&\ j\in B_\alpha\rangle$ is strictly increasing (or just without
repetitions). Then for some $i(*)<\delta$ the set $X_2=\{\alpha\in X_1: i(*)
\in B_\alpha\}$ has cardinality $\sigma$. But then the set
\[X=:\{\alpha\in X_2: \neg(\exists f_\beta\in \cJ_0)(f_\beta\rest(i(*)+1)=
f_\alpha\rest(i(*)+1))\}\]
is of size $\sigma$ (remember $|\cJ_0|<\sigma=\cf(\sigma)$), a
contradiction with the choice of $\cJ'$. \hfill$\square_{\ref{cl5}}$

This finishes the proof of the lemma. \QED$_{\ref{3.1}}$

\begin{lemma}
\label{3.2}
Suppose that $\ga$ is a set of regular cardinals satisfying
\[|\ga|^+<\min(\ga),\quad \lambda=\max\pcf(\ga)\ \mbox{ and }\ [\theta\in
\ga\ \Rightarrow\ \theta>(\max\pcf(\theta\cap\ga))^{<\sigma}].\]
\begin{enumerate}
\item Assume that $\kappa=|\ga|$, $\kappa=\kappa^{<\sigma}$ and for
$\vare<\kappa$, $J$ is a $\sigma$-complete ideal on $\ga$ extending
$J_{<\lambda}[\ga]$ and $\ga_\vare\subseteq\ga$ are pairwise disjoint not in
$J$. If $2^\kappa\ge\lambda$ or just $2^\kappa\geq\sup(\ga)$ then there is
a $\sigma$-entangled linear order of power $\lambda$.
\item We can replace ``$\kappa=|\ga|$'' by ``$\cf(\sup\ga)\leq\kappa$''.

Clearly in parts (1) and (2) we have:\quad $\ga$ has no last element and
$|\ga|\geq\cf(\sup\ga)\geq\sigma$.
\item If in (1) we omit the $\ga_\vare$, still there is a positively
$\sigma$-entangled linear order of power $\lambda$.
\item The linear order above has the exact density $\mu=:\sup\ga$. If there is
a $\mu$-free sequence $\langle f_\alpha:\alpha<\lambda\rangle$ which is
$<_{J_{<\lambda}[\ga]}$--increasing and cofinal (see \cite[\S 1]{Sh:355})
then the linear order has the exact density $(\mu^+,\mu)$ (see
Definition \ref{1.9}).
\end{enumerate}
\end{lemma}

\Proof 1)\ \ \ It follows from part (2) as $\cf(\sup\ga)\leq|\ga|\leq\kappa$.
\medskip

\noindent 2)\ \ \ Let $\langle f_\alpha:\alpha<\lambda\rangle$ be $\leq_{
J_{<\lambda}[\ga]}$-increasing cofinal in $\prod\ga/J_{<\lambda}[\ga]$ with
\[|\{f_\alpha\rest\theta: \alpha<\lambda\}|\leq\max\pcf(\ga\cap\theta)\qquad
\mbox{ for }\ \theta\in\ga\]
(exists by \cite[3.5]{Sh:355}). For each $\theta\in\ga$ we can find sets
$F_{\theta,\zeta}$ (for $\zeta<\kappa$) such that $F_{\theta,\zeta}\subseteq
\{f_\alpha\rest\theta:\alpha<\lambda\}$, and for any disjoint subsets $X,Y$
of $\{f_\alpha\rest\theta:\alpha<\lambda\}$ of cardinality $<\sigma$, for
some $\zeta<\kappa$, $F_{\theta,\zeta}\cap(X\cup Y)=X$ (possible as
$\kappa=\kappa^{<\sigma}$ and $2^\kappa\geq|\{f_\alpha\rest\theta:\alpha<
\lambda\}|$ -- by \cite{EK} or see \cite[AP1.10]{Sh:g}). Clearly $\ga$
has no last element (as $J^{\rm bd}_{\ga}\subseteq J$ and by the
existence of the $\ga_\vare$'s) and $\cf(\sup\ga)\leq\kappa$, so there
is an unbounded $\gb\subseteq\ga$ of cardinality $\leq\kappa$. As $\ga$ can
be partitioned to $\kappa$ pairwise disjoint sets each not in $J$ (and
$J^{\rm bd}_{\ga}\subseteq J$), we can find a sequence
$\langle(\theta_{\Upsilon},\zeta_\Upsilon):\Upsilon\in\ga\rangle$ such that
\begin{itemize}
\item for each $\Upsilon\in\ga$:\qquad $\theta_\Upsilon\in\ga$, $\Upsilon\geq
\theta_\Upsilon$, $\zeta_\Upsilon<\kappa$,\quad and
\item for each $\theta\in\gb$, $\zeta<\kappa$ the set $\{\Upsilon\in\ga:
\theta_\Upsilon=\theta,\zeta_\Upsilon=\zeta\}$ is $\neq\emptyset\ \mod\,J$.
\end{itemize}
Now we define a linear order $<_{et}$ on $\{f_\alpha:\alpha<\lambda\}$ as
follows:

$f_\alpha<_{et} f_\beta$\qquad if and only if \qquad for some
$\Upsilon\in\ga$ we have

\qquad\quad $f_\alpha\rest(\ga\cap\Upsilon)=f_\beta\rest(\ga\cap\Upsilon)$,
$f_\alpha(\Upsilon)\neq f_\beta(\Upsilon)$ and

\qquad\quad $f_\alpha(\Upsilon)<f_\beta(\Upsilon)\quad\Leftrightarrow\quad
f_\alpha\rest\theta_\Upsilon\in F_{\theta_\Upsilon,\zeta_\Upsilon}$.

\noindent Readily $<_{et}$ is a linear order on the set $\cI=\{f_\alpha:
\alpha<\lambda\}$. We are going to show that it is as required (note that in
the definition of $f_\alpha<_{et}f_\beta$ we have $f_\alpha\rest\theta_\Upsilon
=f_\beta\rest\theta_\Upsilon$ as $\theta_\Upsilon\le\Upsilon$). Suppose that
$\vare(*)<\sigma$, $u\subseteq\vare(*)$,$v=\vare(*)\setminus u$,
$t^\vare_\alpha=f_{\gamma(\vare,\alpha)}$ (for $\vare<\vare(*)$, $\alpha<
\lambda$) and $\gamma(\vare,\alpha)$'s are pairwise distinct. For each $\alpha
<\lambda$ take $\theta_\alpha\in\gb$ such that $\{f_{\gamma(\vare,\alpha)}
\rest\theta_\alpha:\vare<\vare(*)\}$ is with no repetitions (possible as
$\vare(*)<\sigma\leq\cf(\sup\ga)\leq\kappa$, $\gb\subseteq\ga$ is
unbounded). Since $\lambda$ is regular, $\lambda>\kappa$, we may assume that
for each $\alpha<\lambda$, $\theta_\alpha=\theta^*\in\gb$. We have that
$|\{f_\alpha\rest\theta^*: \alpha<\lambda\}|\leq\max\pcf(\ga\cap\theta^*)$ and
$(\max\pcf(\ga\cap\theta^*))^{<\sigma}<\theta^*$ and hence we may assume that
for some $\langle g_\vare:\vare<\vare(*)\rangle$:
\[(\forall\alpha<\lambda)(\forall\vare<\vare(*))(f_{\gamma(\vare,\alpha)}
\rest\theta^*=g_\vare).\]
Let $X=\{g_\vare:\vare\in u\}$, $Y=\{g_\vare:\vare\in\vare(*)\setminus u\}$
and let $\zeta<\kappa$ be such that $F_{\theta^*,\zeta}\cap(X\cup Y)=X$. Like
in the proof of \ref{3.1} one can show that
\[\{\mu\in\ga:(\forall\xi<\mu)(|\{\alpha<\lambda:g_\alpha(\mu)>\xi\}|=
\lambda)\}=\ga\ \mod\,J,\]
where $g_\alpha(\mu)=\min\{f_{\gamma(\vare,\alpha)}(\mu):\vare<\vare(*)\}$.
Thus we can find $\Upsilon\in\ga$ such that $\theta^*=\theta_\Upsilon$,
$\zeta=\zeta_\Upsilon$ and $(\forall\xi<\Upsilon)(|\{\alpha<\lambda:
g_\alpha(\Upsilon)>\xi\}|=\lambda)$. Next, as in \ref{3.1}, we can find
$\alpha_\xi<\lambda$ (for $\xi<\Upsilon$) and $S\in [\Upsilon]^\Upsilon$
such that for each $\xi<\Upsilon$ we have $\xi<g_{\alpha_\xi}(\Upsilon)$,
$(\forall\zeta<\xi)(\alpha_\zeta<\alpha_\xi)$ and
\[(\forall\vare<\vare(*))(\forall\zeta<\xi)\left(f_{\gamma(\vare,
\alpha_\zeta)}(\Upsilon)<g_{\alpha_\xi}(\Upsilon)\right)\]
and the sequence $\left\langle\langle f_{\gamma(\vare,\alpha_\xi)}\rest
\Upsilon:\vare<\vare^*\rangle:\xi\in S\right\rangle$ is constant. Choose any
$\xi_1,\xi_2\in S$, $\xi_1<\xi_2$ and note that for every $\vare<\vare(*)$
we have
\[f_{\gamma(\vare,\alpha_{\xi_1})}\rest\Upsilon=f_{\gamma(\vare,\alpha_{
\xi_2})}\rest\Upsilon,\quad f_{\gamma(\vare,\alpha_{\xi_1})}(\Upsilon)<
f_{\gamma(\vare,\alpha_{\xi_2})}(\Upsilon),\]
and $f_{\gamma(\vare,\alpha_{\xi_1})}\rest\theta_\Upsilon=g_\vare=f_{\gamma
(\vare,\alpha_{\xi_2})}\rest\theta_\Upsilon$. Thus $\alpha_{\xi_1}<
\alpha_{\xi_2}$ satisfy the condition given by entangledness for
$t^\vare_\alpha$'s.
\medskip

\noindent 3)\ \ \ Let $\langle f_\alpha:\alpha<\lambda\rangle$ be as in
the proof of part (2). We define a linear order $<_{pet}$ on $\{f_\alpha:
\alpha<\lambda\}$ as follows:

$f_\alpha<_{pet} f_\beta$\qquad if and only if \qquad for some
$\Upsilon\in\ga$ we have

\qquad\quad $f_\alpha\rest(\ga\cap\Upsilon)=f_\beta\rest(\ga\cap\Upsilon)$ and
$f_\alpha(\Upsilon)<f_\beta(\Upsilon)$.

\noindent The rest is even simpler than in the proof of part (2) after
defining $<_{et}$ (remember \ref{1.2}(6)).
\medskip

\noindent 4)\ \ \ It is similar to the proof of \ref{3.1}(2), noting that
\begin{quotation}
\noindent if $\bar{f}^\ell=\langle f^\ell_\alpha:\alpha<\lambda\rangle$ is
$<_J$--increasing cofinal in $\prod\ga$ for $\ell=1,2$,
$\lambda=\cf(\lambda)$ and $\bar{f}^1$ is $\mu$--free

\noindent then for some $X\in [\lambda]^{\textstyle \lambda}$, $\bar{f}^2\rest
X$ is $\mu$--free. \QED$_{\ref{3.2}}$
\end{quotation}

\begin{proposition}
\label{3.3}
\begin{enumerate}
\item Assume
\begin{description}
\item[(a)] $\ens_\sigma(\lambda_i,\kappa_i)$ for $i<\delta$,
\item[(b)] $\lambda_i$ are regular cardinals for $i<\delta$, $\langle
\lambda_i:i<\delta\rangle$ is strictly increasing, $\delta<\lambda_0$,
\item[(c)] $J$ is a $\sigma$--complete ideal on $\delta$ extending
$J^{\rm bd}_{\delta}$,
\item[(d)] $\kappa<T^+_J(\langle\kappa_i:i<\delta\rangle)$ ($=:\sup\{|F|^+:
F\subseteq\prod\limits_{i<\delta}\kappa_i$ and $f\neq g\in F\Rightarrow f
\neq_J g\}$),
\item[(e)] $\prod\limits_{i<\delta}\lambda_i/J$ has the true cofinality
$\lambda$ as exemplified by $\{f_\alpha:\alpha<\lambda\}$ and for each
$i<\delta$, $\lambda_i>|\{f_\alpha\rest i:\alpha<\lambda\}|^{<\sigma}$

(if for each $i$, $\max\pcf(\{\lambda_j: j<i\})<\lambda_i$ then we have such
$f_\alpha$'s).
\end{description}
Then $\ens_\sigma(\lambda,\kappa)$.
\item Assume that in part (1) we omit (d) but in addition we have
\begin{description}
\item[(f)] for each $i<\delta$, $\kappa_i\geq|\{f_\alpha\rest i:\alpha<
\lambda\}|$
\end{description}
[and we have such $f_\alpha$'s e.g.~if $\kappa_i\geq\max\pcf(\{\lambda_j:
j<i\})$], or at least $\lim\inf\limits_J(\kappa_i)=\sup\limits_{i<\delta}
\lambda_i$, or
\begin{description}
\item[(f')] $\delta$ can be partitioned\footnote{If $\kappa_j$ is
non-decreasing then partition to $\cf(\delta)$ sets suffices.} to $|\delta|$
many $J$-positive sets and for each $i<\delta$ for $J$-almost all $j<\delta$
we have $\kappa_j\geq |\{f_\alpha\rest i:\alpha<\lambda\}|$

(if $\kappa_j$ strictly increasing this means ``every large enough $j$'').
\end{description}
Then there is a $\sigma$-entangled linear order $\cI$ of cardinality $\lambda$.
\item Assume (f') or (f) + (g), where
\begin{description}
\item[(g)] there is a decreasing sequence $\langle B_\vare:\vare<\sigma
\rangle$ of elements of $J^+$ with empty intersection.
\end{description}
Then in (2) we can get:\qquad $\cI=\cI_1+\cI_2$ such that for any uniform
ultrafilter $D$ on $\sigma$ the orders $(\cI_1)^\sigma/D$, $(\cI_2)^\sigma/D$
have isomorphic subsets of cardinality $\lambda$ (see \ref{1.8} for a
conclusion).
\item The linear order has exact density $\mu=:\sup\limits_{i<\delta}
\lambda_i$, if $\langle f_\alpha:\alpha<\lambda\rangle$ is $\mu$--free even
exact density $(\mu,\ldots)$.
\item In \ref{3.3}(1) we can weaken clause (d) to:
\begin{description}
\item[(d)$^-$] for some $F\subseteq\prod\limits_{i<\delta}\kappa_i$, $|F|=
\kappa$ and for every $F'\subseteq F$ of cardinality $<\sigma$ we have
\[\{i<\delta: \langle g(i): g\in F'\rangle\mbox{ is with no repetition}\}
\in J^+.\]
\end{description}
\item In \ref{3.3}(2) we can replace (f') by
\begin{description}
\item[($\oplus$)] there is $h:\delta\longrightarrow\delta$ such that $i\geq
h(i)$, $\kappa_i\geq |\{f_\alpha\rest h(i): \alpha<\lambda\}|$ and for every
$i<\delta$, $\{j<\delta: h(j)\geq i\}\in J^+$.
\end{description}
\end{enumerate}
\end{proposition}

\Proof Similar. (About part (3) look at the proof of \ref{2.2}(3)). However
we will give some details.

\noindent 1)\ \ \ As $\ens_\sigma(\lambda_i,\kappa_i)$ (by clause (a)), we
can find linear orders $<^i_\alpha$ of $\lambda_i$ (for $\alpha<\kappa_i$)
such that the sequence $\langle(\lambda_i,<^i_\alpha): \alpha<\kappa_i\rangle$
is $\sigma$-entangled. By clause (d) we can find $g_\zeta\in\prod\limits_{i<
\delta}\kappa_i$ (for $\zeta<\kappa$) such that
\[\vare<\zeta<\kappa\quad\Rightarrow\quad\{i<\delta: g_\vare(i)=g_\zeta(i)\}
\in J.\]
Now for each $\zeta<\kappa$ we define a linear order $\cI_\zeta=(F,<^*_\zeta)$
with the set of elements $F=:\{f_\alpha:\alpha<\lambda\}$ as follows:

$f_\alpha<^*_\zeta f_\beta$\qquad if and only if\qquad for some $i<\delta$ we
have

\qquad\quad $f_\alpha(i)\neq f_\beta(i)$, $f_\alpha\rest i= f_\beta\rest i$
and $f_\alpha(i)<^i_{g_\zeta(i)} f_\beta(i)$.

\noindent It is easy to check that $<^*_\zeta$ is a linear order of $F$. For
the relevant part of (4) note that its density is $\leq|\{f_\alpha\rest i:
\alpha<\lambda$ and $i<\delta\}|=\mu=:\sum\limits_{i<\delta}\lambda_i$. As in
the proof of \ref{3.1},
\[A=:\{i<\delta: (\exists f\in\prod\limits_{j<i}\lambda_j)(\lambda_i=|\{
f_\alpha(i):\alpha<\lambda\ \&\ f=f_\alpha\rest i\}|)\}=\delta\ \mod\;J\]
and for $i\in A$ let $f^*_i$ exemplify $i\in A$. If $G\subseteq F$ is dense,
$|G|$ cannot be $<\mu$ as then it is $<\lambda_{i(*)}$ for some $i(*)\in A$ and
so for some $\gamma<\lambda_{i(*)}$ for no $\alpha<\lambda$ is
\[f_\alpha\rest i(*)=f^*_{i(*)}\ \&\ \gamma\leq f_\alpha(i(*))<\lambda_{i(*)
},\]
thus proving $\mu=\dens(\cI_\zeta)$. The part ``if $\{f_\alpha:\alpha<\lambda
\}$ is $\mu$-free then any $\cJ\subseteq\cI_j$ of cardinality $\geq\mu$ has
density $\mu$ (i.e.~has exact $(\mu^+,\mu)$-density)'' can be proven similarly.
\medskip

\noindent Finally ``$\langle\cI_\zeta:\zeta<\kappa\rangle$ is
$\sigma$-entangled'' is proved as in the proof of \ref{3.1}. Assume $u\cup v
=\vare(*)<\sigma$, $u\cap v=\emptyset$, and $\zeta_\vare<\kappa$ for $\vare<
\vare(*)$ are pairwise distinct. Now
\[A=:\{i<\delta:\langle g_\vare(i): \vare<\vare(*)\rangle\mbox{ are pairwise
distinct}\}\neq\emptyset\ \mod\;J\]
(as $J$ is $\sigma$-complete, (d)\ $\Rightarrow$\ (d)$^-$). We continue as in
the proof of \ref{3.1} (only with $A$ as here) and using $\langle(\lambda_i,
<^i_\alpha):\alpha<\kappa_i\rangle$ is $\sigma$-entangled.
\medskip

\noindent 2)\ \ \ First we assume clause (f). As $\ens_\sigma(\lambda_i,
\kappa_i)$ and $\kappa_i\geq|\Pi_i|$ where $\Pi_i=:\{f_\alpha\rest i:\alpha<
\lambda\}$, we can find linear orders $<^i_\eta$ of $\lambda_i$ (for $\eta\in
\Pi_i$) such that $\langle(\lambda_i,<^i_\eta):\eta\in\Pi_i\rangle$ is
$\sigma$-entangled. We define the linear order $<^*$ of $F=:\{f_\alpha:\alpha
<\lambda\}$ as follows:

$f_\alpha<^* f_\beta$\qquad if and only if\qquad for some $i<\delta$ we have
\[f_\alpha(i)\neq f_\beta(i),\quad f_\alpha\rest i= f_\beta\rest
i,\quad\mbox{and }\ f_\alpha(i)<^i_{f_\alpha\rest i} f_\beta(i).\]
The rest is as in \ref{3.3}(1).
\medskip

\noindent Next we assume clause (f') instead of (f). So let $\langle A_i:
i<\delta\rangle$ be a partition of $\delta$ with every $A_j$ in $J^+$ and so
necessarily 
\[A_i'=:\{j\in A_i: \kappa_j\geq|\Pi_i|\}=A_i\ \mod\; J.\]
Then we can choose\footnote{Actually we can replace the assumption (g)
(in \ref{3.3}(3)) by the existence of such $h$.} a function $h$ such
that
\begin{description}
\item[$(\otimes)$] $h:\delta\longrightarrow\delta$, $h(i)\leq i$, $\kappa_i
\geq |\Pi_{h(i)}|$ and for every $i<\delta$ we have 
\[\{j<\delta: h(j)\geq i\}\in J^+.\]
\end{description}
Let $\langle(\lambda_i,<^i_\eta):\eta\in\Pi_{h(i)}\rangle$ be a
$\sigma$-entangled sequence of linear orders. Now we define a linear order on
$F=:\{f_\alpha: \alpha<\lambda\}$:

$f_\alpha<^* f_\beta$\qquad if and only if\qquad for some $i<\delta$ we have
\[f_\alpha(i)\neq f_\beta(i),\quad f_\alpha\rest i= f_\beta\rest i,\quad
\mbox{and }\ f_\alpha(i)<^i_{f_\alpha\rest h(i)} f_\beta(i).\]
\medskip

\noindent 3)\ \ \ Without loss of generality, for each $\xi<\sigma+\sigma$ the
set $J_\xi=:\{\alpha<\lambda:f_\alpha(0)=\xi\}$ has cardinality $\lambda$. Let
$F=\{f_\alpha: \alpha<\lambda\ \&\ f_\alpha(0)<\sigma+\sigma\}$. We can find
$B_\vare\in J^+$ (for $\vare<\sigma$), decreasing with $\vare$ and such that
$\bigcap\limits_{\vare<\sigma}B_\vare=\emptyset$ and in the proof from (f')
replace $(\otimes)$ by
\begin{description}
\item[$(\otimes')$] $h:\delta\longrightarrow\delta$, $h(i)\leq i$, $\kappa_i
\geq |\Pi_{h(i)}|$ and for every $i<\delta$, $\xi<\sigma+\sigma$ we have
\[\{j\in B_\xi: h(j)\geq i\}\in J^+,\]
\end{description}
and define $<^*$ as

$f_\alpha<^* f_\beta$\qquad if and only if\qquad for some $i<\delta$ we have
\[\begin{array}{l}
f_\alpha(i)\neq f_\beta(i),\quad f_\alpha\rest i=f_\beta\rest i,\quad\mbox{
and}\\
{[i=0\ \Rightarrow\ f_\alpha(0)<f_\beta(0)]},\quad [0<i\in B_{f_\alpha(0)}\
\Rightarrow\ f_\alpha(i)<^i_{f_\alpha\rest h(i)} f_\beta(i)]\quad\mbox{ and}\\
{[0<i\notin B_{f_\alpha(0)}\ \Rightarrow\ f_\alpha(i)< f_\beta(i)]}.
\end{array}\]
The rest is as before (we can replace $\sigma$ by other cardinal $\geq\sigma$
but $\leq\lambda_0$).
\medskip

\noindent 4)\ \ \ For \ref{3.3}(1) see in its proof, other cases similar.
\medskip

\noindent 5)\ \ \ Really included in the proof of \ref{3.3}(1).
\medskip

\noindent 6)\ \ \ Really included in the proof of \ref{3.3}(2).
\QED$_{\ref{3.3}}$

\begin{proposition}
\label{3.4}
\begin{enumerate}
\item Assume that:
\begin{description}
\item[(a)] $\ens_\sigma(\lambda_i,\mu_i,\kappa_i)$ for $i<\delta$,
\item[(b)] $\langle\lambda_i:i<\delta\rangle$ is a strictly increasing
sequence of regular cardinals, $2^{|\delta|}<\lambda_0$,
\item[(c)] $J$ is a $\sigma$-complete ideal on $\delta$ extending
$J^{\rm bd}_\delta$,
\item[(d)] $\kappa<T^+_J(\langle\kappa_i:i<\delta\rangle)$,
$\sum\limits_{i<\delta}\lambda_i\leq\mu<\lambda$, $\mu=\cf(\mu)$ and
\[(\forall\alpha<\mu)(|\alpha|^{<\delta}<\mu),\]
\item[(e)] $F=\{f_\alpha:\alpha<\lambda\}\subseteq\prod\limits_{i<\delta}
\lambda_i$, $f_\alpha\neq_J f_\beta$ for $\alpha\neq\beta$, and 
\[|\{f_\alpha\rest i:\alpha<\lambda\}|^{<\sigma}<\mu_i,\]
\item[(f)] if $\mu'_i=\cf(\mu'_i)<\lambda_i$, $A\in J^+$ then $\tcf(
\prod\limits_{i\in A}\mu'_i/J)=\mu$ is impossible,
\item[(g)] $(\forall\alpha<\mu)(|\alpha|^{<\sigma}<\mu)$.
\end{description}
Then $\ens_\sigma(\lambda,\mu,\kappa)$.
\item Assume in addition
\begin{description}
\item[(h)] $\kappa_i\geq |\{f_\alpha\rest i:\alpha<\lambda\}|$,
\end{description}
or at least
\begin{description}
\item[(h')] $\cf(\delta)=\omega$ and $\lim\inf_J\kappa_i=\sup\limits_{i<
\delta}\lambda_i$,
\end{description}
or at least
\begin{description}
\item[(h'')] there is $h:\delta\longrightarrow\delta$ such that $i\geq h(i)$,
$\kappa_i\geq|\{f_\alpha\rest h(i): \alpha<\lambda\}$ and $\delta=
\lim\sup\limits_{i<\delta} h(i)$.
\end{description}
Then there is a $(\mu,\sigma)$-entangled linear order $\cI$ of cardinality
$\lambda$.
\item Suppose also that
\begin{description}
\item[(i)] we can partition $\delta$ into $\sigma$ sets from $J^+$ (or
clause (g) from \ref{3.3}(3) holds).
\end{description}
Then we can get:\quad for any uniform ultrafilter $D$ on $\sigma$,
$\cI^\sigma/D$ has two isomorphic subsets with disjoint convex hulls of
cardinality $\lambda$.
\end{enumerate}
\end{proposition}

\Proof Similar proof but for reader's convenience we will present some
details.

\noindent 1)\ \ \ We repeat the proof of \ref{3.3}(1) up to proving
entanglness. To show ``$\langle\cI_\zeta:\zeta<\kappa\rangle$ is
$(\mu,\sigma)$--entangled'' suppose that $u\cup v=\vare(*)<\sigma$,
$u\cap v=\emptyset$ and $\langle\zeta_\vare:\vare<\vare(*)\rangle$ is a
sequence of pairwise distinct ordinals $<\kappa$ and $\langle\gamma(
\beta,\vare):\beta<\mu,\vare<\vare(*)\rangle\subseteq\lambda$ is such
that 
\[(\forall\vare<\vare(*))(\forall\beta_1<\beta_2<\lambda)(\gamma(\beta_1,
\vare)\neq\gamma(\beta_2,\vare)).\]
We want to find $\beta_1<\beta_2<\mu$ such that 
\[(\forall\vare<\vare(*))(\gamma(\beta_1,\vare)<_{\zeta_\vare}\gamma(\beta_2,
\vare)\ \Leftrightarrow\ \vare\in u).\]
\begin{claim}
\label{3.4A}
Assume that
\begin{description}
\item[($\alpha$)] $\langle\lambda_i:i<\delta\rangle$ is a strictly
increasing sequence of regular cardinals,
\item[($\beta$)]  $J$ is a $\sigma$-complete ideal on $\delta$ extending
$J^{\rm bd}_\delta$,
\item[($\gamma$)] a sequence $\bar{f}=\langle f_\alpha:\alpha<\lambda\rangle
\subseteq\prod\limits_{i<\delta}\lambda_i$ is $<_J$--increasing,
\item[($\delta$)] $\sum\limits_{i<\delta}\mu_i<\mu=\cf(\mu)\leq\lambda$
and $(\forall\alpha<\mu)(|\alpha|^{<\sigma}<\mu)$,
\item[($\vare$)]  one of the following occurs:
\begin{description}
\item[(i)]  $2^{|\delta|}<\lambda_0$,\quad or
\item[(ii)] $\bar{f}$ is $\mu$--free,
\end{description}
\item[($\zeta$)]  if $\mu'_i<\mu_i$ for $i<\delta$ and $A\in J^+$ then
$\tcf(\prod\limits_{i\in A}\mu_i/J)=\mu$ is impossible,
\item[($\theta$)] a sequence $\langle\gamma(\beta,\vare):\beta<\mu,\ \vare<
\vare(*)<\sigma\rangle$ of ordinals $<\lambda$ satisfies
\[(\beta_1,\vare_1)\neq (\beta_2,\vare_2)\quad\Rightarrow\quad\gamma(\beta_1,
\vare_1)\neq\gamma(\beta_2,\vare_2).\]
\end{description}
Then there are a set $X\in [\mu]^{\textstyle \mu}$ and a sequence $\langle
h_\vare:\vare<\vare(*)\rangle\subseteq\prod\limits_{i<\delta}\lambda_i$
such that
\begin{description}
\item[(a)] for each $\vare<\vare(*)$ the sequence $\langle\gamma(\beta,\vare):
\beta\in X\rangle$ is strictly increasing,
\item[(b)] $(\forall\beta\in X)(\forall\vare<\vare(*))(f_{\gamma(\beta,\vare)}
<_J h_\vare)$,
\item[(c)] $\langle h_\vare(i):i<\delta\rangle$ is the $<_J$-eub of
$\langle f_{\gamma(\beta,\vare)}:\beta\in X\rangle$,
\item[(d)] $B^*=\{i<\delta: (\forall\vare<\vare(*))(\cf(h_\vare(i))\geq\mu_i)
\}=\delta\ \mod\;J$.
\end{description}
\end{claim}

\noindent{\em Proof of the claim:}\hspace{0.2in} Since $(\forall\alpha<\mu)(
|\alpha|^{<\sigma}<\mu)$ and $\cf(\mu)=\mu$ we know that for some $X\in
[\mu]^{\textstyle\mu}$ we have
\[(\forall\vare<\vare(*))(\mbox{the sequence }\langle\gamma(\beta,\vare):\beta
\in X\rangle\mbox{ is strictly increasing}).\]
[Why? For $\beta<\mu$, $\vare<\vare(*)$ define $f(\beta,\vare)$ as follows:
\begin{quotation}
\noindent if there exists $\delta$ such that $\gamma(\delta,\vare)>\gamma(
\beta,\vare)$ and $[\gamma(\beta,\vare),\gamma(\delta,\vare)]\cap\{\gamma
(\alpha,\vare):\alpha<\beta \}=\emptyset$ then $f(\beta,\vare)$ is this unique
$\delta$,

\noindent otherwise $f(\beta,\vare)=-1$.
\end{quotation}
By Fodor Lemma, there is a stationary set $S\subseteq\mu$ such that $\sup f[S
\times\vare(*)]<\mu$. Since $|\alpha|^{\vare(*)}<\mu$, on a stationary set
$X\subseteq S$ the sequence $\langle f(\alpha,\vare):\vare<\vare(*)\rangle$
does not depend on $\alpha$. This $X$ is as required.]\\
By renaming we may assume that $X=\mu$. Consequently, for each $\vare<
\vare(*)$ the sequence $\langle f_{\gamma(\beta,\vare)}:\beta<\mu\rangle$
is $<_J$--increasing. Since $\mu=\cf(\mu)>2^{|\delta|}$\ \ or\ \ the
sequence is $\mu$-free, we may use \cite[\S 1]{Sh:355} to conclude that
it has a $<_J$--eub, call it $h_\vare$. We may assume that, for each
$i<\delta$, $h_\vare(i)$ is a limit ordinal. Since $h_\vare<_J\langle
\lambda_i:i<\delta\rangle$, wlog $(\forall i<\delta)(h_\vare(i)<\lambda_i)$.
Also $\mu=\tcf(\prod\limits_{i<\delta}\cf(h_\vare(i))/J)$ and
$\cf(h_\vare(i))\leq h_\vare(i)<\lambda_i$, so by the assumption $(\vare)$ we
have
\[\{i<\delta:\cf(h_\vare(i))<\mu_i\ \mbox{ or }\ \max\pcf(\{\lambda_j:
j<i\}<\cf(h_\vare(i))\})=\emptyset\ \mod\;J.\]
Since $J$ is $\sigma$-complete we conclude
\[B^*=\{i<\delta:(\forall\vare<\vare(*))(\cf(h_\vare(i))\geq\mu_i)\}=\delta\
\mod\; J,\]
finishing the proof of the claim. \hfill$\square_{\ref{3.4A}}$
\medskip

Note that we may assume that for some $i^\otimes<\delta$, for every
$\beta<\mu$ the sequence $\langle f_{\gamma(\beta,\vare)}\rest i^\otimes:
\vare<\vare(*)\rangle$ is with no repetition and does not depend on
$\beta$. Now we may apply \ref{3.4A} to find $X\in [\mu]^{\textstyle\mu}$
and $\langle h_\vare:\vare<\vare(*)\rangle$ as there. We shall continue
like in the proof of \ref{3.1} with some changes, however. We let
\[G_i=\{\langle f_{\gamma(\beta,\vare)}(i):\vare<\vare(*)\rangle:\beta<\mu
\}\subseteq\prod_{\vare<\vare(*)}h_\vare(i),\qquad\mbox{ and}\]
\[\begin{array}{ll}
B=\{i\in B^*:&\mbox{ for each }\langle\xi_\vare:\vare<\vare(*)\rangle\in
\prod_{\vare<\vare(*)}h_\vare(i)\\
\ &\mbox{ there are $\mu$ many ordinals }\beta<\mu\mbox{ such that}\\
\ &\quad (\forall\vare<\vare(*))(\xi_\vare<f_{\gamma(\beta,\vare)}(i))\}.
\end{array}\]
We have to show the following.

\begin{claim}
\label{3.4B}
$B=\delta\ \mod\; J$.
\end{claim}

\noindent{\em Proof of the claim:}\hspace{0.2in} Assume not. Then, as
$B^*=\delta\ \mod\; J$, necessarily $B^*\setminus B\neq\emptyset\ \mod\;
J$. For each $i\in B^*\setminus B$ choose a sequence $\langle\xi^i_\vare:
\vare<\vare(*)\rangle\in \prod\limits_{\vare<\vare(*)} h_\vare(i)$ and
an ordinal $\beta_i<\mu$ exemplifying $i\notin B$. Thus

if $i\in B^*\setminus B$ and $\beta\in [\beta_i,\mu)$ then $(\exists
\vare<\vare(*))(\xi^i_\vare\geq f_{\gamma(\beta,\vare)}(i))$.

\noindent For $\vare<\vare(*)$ define a function $h^\vare\in\prod\limits_{
i<\delta}h_\vare(i)$ by
\[h^\vare(i)=\left\{\begin{array}{ll}
\xi^i_\vare+1 &\mbox{if }\ i\in B^*\setminus B,\\
0             &\mbox{if }\ i\in \delta\setminus B^* \mbox{ or } i\in B.
\end{array}\right.\]
Now, for each $\vare$, for every sufficiently large $\beta<\mu$ we have
$h^\vare<_J f_{\gamma(\beta,\vare)}$. Consequently, we find $\beta^*<\mu$
such that
\[\vare<\vare(*)\ \ \&\ \ \beta\in [\beta^*,\mu)\quad \Rightarrow\quad
h^\vare<_J f_{\gamma(\beta,\vare)}.\]
But the ideal $J$ is $\sigma$-complete, so for each $\beta\in [\beta^*,
\mu)$
\[B_\beta=:\{i<\delta: (\forall\vare<\vare(*))(h^\vare(i)<f_{\gamma(
\beta,\vare)}(i))\}=\delta\ \mod\; J.\]
Now we may take $\beta\in [\beta^*,\mu)$ and then choose $i\in B_\beta\cap
(B^*\setminus B)$ and get a contradiction as in the proof of \ref{3.1A}.
\hfill$\square_{\ref{3.4B}}$
\medskip

Remember that
\[|\{\langle f_{\gamma(\beta,\vare)}\rest i: \vare<\vare(*)\rangle:
\beta<\mu\}|\leq |\{f_\alpha\rest i:\alpha<\lambda\}|^{<\sigma}<\mu_i\]
(not just $<\lambda_i$), see clause (e) of the assumptions of
\ref{3.4}(1). Hence there is $\beta^\otimes<\mu$ such that
\[(\forall j<\delta)(|\{\langle f_{\gamma(\beta,\vare)}\rest j:\vare<\vare(*)
\rangle:\beta\in [\beta^\otimes,\mu)\}|=\mu).\]
For each $i\in B$ (defined above) we know that $(\forall\vare<\vare(*))(
\cf(h_\vare(i))\geq\mu_i)$, and hence $\prod\limits_{\vare<\vare(*)}
h_\vare(i)$ is $\mu_i$-directed and the set
\[\{\langle f_{\gamma(\beta,\vare)}(i):\vare<\vare(*)\rangle:\beta\in
[\beta^\otimes,\mu)\}\] 
is cofinal in $\prod\limits_{\vare<\vare(*)}h_\vare(i)$. Putting these
together, there are $\nu^i_\vare\in\prod\limits_{j<i}\lambda_j$ (for
$\vare<\vare(*)$) such that 
\begin{quotation}
\noindent for every
$\langle\xi_\vare:\vare<\vare(*)\rangle\in\prod\limits_{\vare<\vare(*)}
h_\vare(i)$ for some $\beta\in [\beta^\otimes,\mu)$ we have
$(\forall\vare<\vare(*))(f_{\gamma(\beta,\vare)}\rest i=\nu^i_\vare\
\&\ \xi_\vare<f_{\gamma(\beta,\vare)}(i))$.
\end{quotation}
Now take any $i\in B$, $i>i^\otimes$ such that the sequence $\langle g_{
\zeta_\vare}(i):\vare<\vare(*)\rangle$ is with no repetition. Again, as
$\prod\limits_{\vare<\vare(*)}h_\vare(i)$ is $\mu_i$--directed we can choose
by induction on $\alpha<\mu_i$ a sequence $\langle\beta_\alpha:\alpha<\mu_i
\rangle\subseteq\mu_i$ such that for each $\alpha<\mu_i$
\[\sup\{\beta_{\alpha'}:\alpha'<\alpha\}<\beta_\alpha\quad\mbox{and}\quad
(\forall\alpha'\!<\!\alpha)(\forall\vare\!<\!\vare(*))(f_{\gamma(\beta_{
\alpha'},\vare)}(i)<f_{\gamma(\beta_\alpha,\vare)}(i)).\]
Now remember that the sequence $\langle<^i_\zeta:\zeta<\kappa_i\rangle$
exemplifies $\ens_\sigma(\lambda_i,\mu_i,\kappa_i)$. So we apply this to 
$\langle g_{\zeta_\vare}(i):\vare<\vare(*)\rangle$ and $\big\langle
\langle f_{\gamma(\beta_\alpha,\vare)}:\vare<\vare(*)\rangle:\alpha<\mu\big
\rangle$, and we find $\alpha_1<\alpha_2<\mu$ such that
\begin{quotation}
\noindent $\vare\in u\quad\Rightarrow\quad f_{\gamma(\alpha_1,\vare)}(i)
<^i_{g_{\zeta_\vare}(i)}f_{\gamma(\alpha_2,\vare)}(i)$,

\noindent $\vare\in v\quad\Rightarrow\quad f_{\gamma(\alpha_2,\vare)}(i)
<^i_{g_{\zeta_\vare}(i)}f_{\gamma(\alpha_1,\vare)}(i)$,
\end{quotation}
so we are done.
\medskip

\noindent 2)\ \ \ The proof is exactly like that of \ref{3.4}(1) except
of two points. First we have to define a linear order $\cI$ (rather than
$\cI_\zeta$ for $\zeta<\kappa$). Assuming that the clause (h) of the
assumptions holds, for each $i$ we can find linear orders $<^i_\eta$ on
$\lambda_i$ (for $\eta\in T_i=\{f_\alpha\rest i: \alpha<\lambda\}$)
such that $\langle (\lambda_i,<^i_\eta): \eta\in T_i\rangle$ is a
$(\mu_i,\sigma)$-entangled sequence of linear orders. This does no
affect the proof except in the very end when we use the entangledness.
\medskip

\noindent 3)\ \ \ Combine the proofs above. \QED$_{\ref{3.4}}$

\begin{remark}
{\em
\begin{enumerate}
\item We can also vary $\sigma$.
\item The ``$2^{|\delta|}<\lambda_0$'' rather than just ``$|\delta|<
\lambda_0$'' is needed just to have $<_J$-eub (to use that ``if $\mu_i=
\cf(\mu_i)<\mu_i\ldots$'') so if $\{f_\alpha:\alpha<\lambda\}$ is
$(\sum\limits_{i<\delta}\lambda_i)$-free, we can weaken ``$2^{|\delta|}<
\lambda_0$'' to ``$|\delta|<\lambda_0$''.
\item Instead of $T^+_J(\langle\kappa_i:i<\delta\rangle)$ we may use any
$\chi=|G|$, $G\subseteq\prod\limits_{i<\delta}\kappa_i$ such that for every
sequence $\langle g_\vare:\vare<\vare(*)\rangle$ of distinct elements of
$G$ the set
\[\{i<\delta:\langle g_\vare(i):\vare<\vare(*)\rangle\mbox{ is with no
repetition}\}\]
belongs to $J^+$. But then in \ref{3.4}(3) we have to change (i).
\end{enumerate}
}
\end{remark}

\begin{proposition}
\label{3.5}
\begin{enumerate}
\item Suppose $\mu=\mu^{<\sigma}$. Then the set
\[\{\delta<\mu^{+4}:\mbox{ if }\cf(\delta)\geq\sigma\mbox{ then }
\ens_\sigma(\mu^{+\delta+1},2^{\cf(\delta)})\mbox{ or }\ens_\sigma(\mu^{
+\delta+1},2^{\mu^+})\}\]
contains a club of $\mu^{+4}$.
\item If in addition $2^\mu\geq\aleph_{\mu^{+4}}$ (or $2^{\mu^+}\geq
\aleph_{\mu^{+4}}$) then the set
\[\begin{array}{ll}
\{\delta\leq\mu^{+4}:&\mbox{ if }\cf(\delta)\geq\mu\mbox{ (or }\cf(\delta)
\geq\mu^+\mbox{ )\ \ then}\\
\ &\mbox{ there is a $\sigma$-entangled linear order in }\mu^{+\delta+1}\}
\end{array}\]
contains a club of $\mu^{+4}$ and $\mu^{+4}$ itself. (We can weaken the
assumptions.)
\item  We can add in part (2) the conditions needed for \ref{1.7}. Also in
parts (1), (2) the exact density of the linear orders is $\mu^{+\delta}$
provided $\cf(\delta)\le\mu^+$.
\end{enumerate}
\end{proposition}

\Proof 1)\ \ \ By \cite[\S 4]{Sh:400}, for some club $C$ of $\mu^{+4}$,
\[(*)\qquad\qquad\alpha<\delta\in C\quad \Rightarrow\quad \mu^{+\delta}>
\cov(\mu^{+\alpha},\mu^+,\mu^+,2),\]
and hence, if $\cf(\delta)\ge\sigma$ and $\delta\in C$ then $(\mu^{+\delta}
)^{<\sigma}=\mu^{+\delta}$.

\noindent Let $\delta$ be an accumulation point of $C$ of cofinality
$\geq\aleph_1$ and let $A\subseteq\delta$ be a closed unbounded set such that
$\prod\limits_{\alpha\in A}\mu^{+\alpha+1}/J^{\rm bd}_A$ has the true
cofinality $\mu^{+\delta+1}$ and $\otp(A)=\cf(\delta)$ (exists by
\cite[2.1]{Sh:355}). Now for $\beta\in A$ we have $\max\pcf(\{\mu^{+\alpha+1}:
\alpha\in\beta\cap A\})$ is $<\mu^{+\delta}<\mu^{+\delta+1}$ if $\cf(\delta)
\leq\mu$ by $(*)$ (and \cite[5.4]{Sh:355}). Hence for some closed unbounded
set $B\subseteq A$ we have
\[\alpha\in\nacc(B)\quad\Rightarrow\quad\cov(\mu^{+\sup(A\cap\alpha)+1},
\mu^+,\mu^+,2)<\mu^{+\alpha}.\]
Hence, if $\alpha\in\nacc(B)$ then  $\max\pcf(\{\mu^{+(\beta+1)}:\beta\in
B\cap\alpha\})<\mu^{+\alpha+1}$. Wlog $\otp(B)=\cf(\delta)$, $B=A$. Now if
$\sigma\leq\cf(\delta)\leq\mu^+$ then we may apply \ref{3.1} to
$\{\mu^{+\alpha+1}:\alpha\in\nacc(A)\}$ (for $\bar{\lambda}$) and get
$\ens_\sigma(\mu^{+\delta+1}, 2^{\cf(\delta)})$. We still have to deal
with the case $\cf(\delta)>\mu^+$. We try to choose by induction on $i$
ordinals $\alpha_i\in A\setminus\bigcup\limits_{j<i}(\alpha_j+1)$ such that
\[\mu^{+\alpha_i}>\max\pcf(\{\mu^{+\alpha_j}:j<i\})<\sigma,\quad
\cf(\alpha_i)=\mu^+.\]
For some $\gamma$ the ordinal $\alpha_i$ is defined if and only if $i<\gamma$.
Necessarily $\gamma$ is limit and by $(*)$ + \cite[1.1]{Sh:371} we have
$\cf(\gamma)\geq\mu^+$. Now for each $i<\gamma$, $\ens_\sigma(\mu^{+\alpha_i
+1},2^{\mu^+})$ holds by the part we have already proved (as $\cf(\alpha_i)
=\mu^+$). So we may apply \ref{3.3}(1).
\medskip

\noindent 2)\ \ \ Let $C$ be as in the proof of part (1) and exemplifying
its conclusion. If $\delta=\sup(C\cap\delta)<\mu^{+4}$, $\cf(\delta)\leq
\mu^+$, $2^{\cf(\delta)}\geq\mu^{+\delta}$ then we can apply \ref{3.2}(1)
(and the proof of part (1)). So we are left with the case $\cf(\delta)>\mu^+$.
Now we can repeat the proof of part (1) the case $\cf(\delta)>\aleph_0$.
Choose $A$ as there and also $\alpha_i\in A$, but demand in addition
$\ens_\sigma(\mu^{+(\delta+1)},2^{\mu^+})$, $\alpha_i\in\mbox{acc}(C)$ and
$\cf(\alpha_i)=\mu^+$, hence $\ens_\sigma(\mu^{+(\delta+1)},2^{\mu^+})$. In
the end apply \ref{3.3}(2) to $\langle\mu^{+(\alpha_i+1)}:i<\gamma\rangle$.
\medskip

\noindent 3)\ \ \ Similar to the proof of \ref{3.3}(4), \ref{3.3}(5).
\QED$_{\ref{3.5}}$

\begin{conclusion}
\label{3.6}
Assume $\sigma\geq\aleph_0$. Then for arbitrarily large cardinals $\lambda$
there is are $\sigma$-entangled linear orders of cardinality $\lambda^+$.
\end{conclusion}

\Proof Let $\chi>\sigma$ be given. We choose by induction on $i<\sigma$
regular cardinals $\lambda_i>\chi$ such that $\ens_\sigma(\lambda_i,\aleph_0
+\prod\limits_{j<i}\lambda_j)$ holds and $\lambda_i>\prod\limits_{j<i}
\lambda_j$. The inductive step is done by \ref{3.5}. Now for some
$\sigma$-complete ideal $I$ on $\sigma$ extending $J^{\rm bd}_\sigma$,
$\prod\limits_{j<i}\lambda_j/I$ has a true cofinality, say $\lambda$. By
\ref{3.3} there is a $\sigma$-entangled linear order of cardinality $\lambda$,
so if $\lambda$ is a successor cardinal we are done (as $\lambda>\chi$). If
not, necessarily $\lambda$ is inaccessible and letting $\mu=\sum\limits_{j<i}
\lambda_j$ clearly $\mu=\mu^{<\sigma}<\lambda\leq\mu^\sigma$. Now we use
\ref{3.5}(2) to find $\lambda_i\in(\mu,\aleph_{\mu^{+4}})$ such that there is
an entangled linear order in $\lambda^+_i$, so in any case we are done.
\QED$_{\ref{3.6}}$

\section{Boolean Algebras with neither pies nor chains}
Let us recall the following definition.
\begin{definition}
\label{4.1}
Let $B$ be a Boolean algebra.
\begin{description}
\item[(a)] We say that a set $Y\subseteq B$ is a {\em chain} of $B$ if
\[(\forall x,y\in Y)(x\neq t\ \Rightarrow\ x<_B y\mbox{ or } y<_B x).\]
\item[(b)] We say that a set $Y\subseteq B$ is a {\em pie} of $B$ if
\[(\forall x,y\in Y)(x\neq y\ \Rightarrow x\not\leq y\mbox{ and }
y\not\leq x).\]
\item[(c)] $\pi(B)$, {\em the (algebraic) density of $B$}, is
\[\min\{|X|: X\subseteq B\setminus\{0\}\mbox{ and }(\forall y)(\exists x
\in X)(0<_B y\in B\ \Rightarrow\ x\leq_B y)\}.\]
\end{description}
\end{definition}

\begin{lemma}
\label{4.2}
\begin{enumerate}
\item Suppose that
\begin{description}
\item[(a)] $\langle\lambda_i:i<\delta\rangle$ is a strictly increasing
sequence of regular cardinals, $\lambda$ is a regular cardinal,
\item[(b)] $J$ is a $\sigma$-complete ideal on $\delta$ extending
$J^{\rm bd}_\delta$,
\item[(c)] $\langle f_\alpha:\alpha<\lambda\rangle$ is a $<_J$-increasing
sequence of functions from $\prod\limits_{i<\delta}\lambda_i$, cofinal in
$(\prod\limits_{i<\delta}\lambda_i,<_J)$,
\item[(d)] for every $i<\delta$, $|\{f_\alpha\rest i:\alpha<\lambda\}|^{
<\sigma}<\lambda_i$,
\item[(e)] $\langle A_\zeta:\zeta<\kappa\rangle\subseteq\delta$ is a
sequence of pairwise disjoint sets such that for every $B\in J$ and
$\zeta<\kappa$ there is $i\in\delta$ such that $\{2i,2i+1\}\subseteq
A_\zeta\setminus B$,
\item[(f)] $2^\kappa\geq\mu=:\sup\limits_{i<\delta}\lambda_i$ and
$\kappa=\kappa^{<\sigma}$ (so $\kappa\leq |\delta|\leq\mu$) and $\cf(\delta)
\leq\kappa$.
\end{description}
Then there is a Boolean algebra $B$ of cardinality $\lambda$ such that:
\begin{description}
\item[$(\oplus)^B_\lambda$] $B$ has neither a chain of cardinality $\lambda$
nor a pie of cardinality $\lambda$ (i.e. $\inc^+(B)\leq\lambda$, Length$^+(B)
\leq\lambda$).
\item[$(\otimes)^B_\mu$] $B$ has the algebraic density $\pi(B)=\mu$ (in fact,
for $a\in B\setminus\{0\}$, $\pi(B\rest a)=\mu$)
\end{description}
This applies also to the $\sigma$-complete algebra which $B$ generates,
provided $(\forall\alpha<\lambda)(|\alpha|^{<\sigma}<\lambda)$.
\item Suppose $2\leq n^*\leq\omega$ and that in part (1) we replace (e) by
\begin{description}
\item[(e)$^{n^*}_\kappa$]
\begin{description}
\item[($\alpha$)] $A_\zeta\subseteq\delta$ for $\zeta<\kappa$, are pairwise
disjoint,
\item[($\beta$)]  $e$ is an equivalence relation on $\delta$ such that each
equivalence class is a finite interval,
\item[($\gamma$)] for every $n<n^*$, $B\in J$ and $\zeta<\kappa$ for some
$\alpha<\delta$ we have:\quad $\alpha/e\subseteq A_\zeta\setminus B$ and
$|\alpha/e|\geq n$.
\end{description}
\end{description}
Then there is a Boolean algebra $B$ of cardinality $\lambda$ as in part (1)
but $(\oplus)^B_\lambda$ is strengthened to
\begin{description}
\item[$(\circledast)^B_{\lambda,n}$]  if:
\begin{description}
\item[($\alpha$)] $a_\alpha\in B$ for $\alpha<\lambda$ are distinct,
\item[($\beta$)]  $n<n^*$,
\item[($\gamma$)] $B^*_n$ is the finite Boolean algebra of subsets of $n
\times (n+1)$, and for $\ell\leq n$, $f_\ell:n\longrightarrow n$ is a function
such that
\[\ell<m<n\ \&\ i<n\quad\Rightarrow\quad f_\ell(i)\neq f_m(i)\qquad\mbox{
and}\]
\[x_\ell=\{(i,j):i<n,\ j<n+1,\ j\leq f_\ell(i)\}$$
\end{description}
then for some $\alpha_0<\cdots<\alpha_{n-1}<\lambda$, the quantifier free type
which $\langle a_{\alpha_0},\ldots a_{\alpha_{n-1}}\rangle$ realizes in $B$
is equal to the quantifier free type which $\langle x_0,\ldots,x_{n-1}\rangle$
realizes in $B^*$.
\end{description}
\end{enumerate}
\end{lemma}

\noindent{\bf Remark:}\qquad 1)\ \ The case $\delta=\sup\limits_{i<\delta}
\lambda_i$ is included.

\noindent 2)\ \ Of course, no order Boolean Algebra of cardinality $\lambda$
can satisfy $(\oplus)^B_\lambda$.

\noindent 3)\ \ Again, $B$ has density $\mu$ and if $\langle f_\alpha:\alpha<
\lambda\rangle$ is $\mu$-free then $B$ has the exact density $\mu$. 
\medskip

\Proof We shall prove only  part (1) as the proof of part (2) is similar.
Without loss of generality $\delta$ is additively indecomposable.\\
We define $B$ as an algebra of subsets of $Y=\bigcup\limits_{i<\delta}Y_{2i}$
where $Y_i=\{f_\alpha\rest i:\alpha<\lambda\}$ for $i<\delta$. For each $i
<\delta$ we can find subsets $F_{2i,\zeta}$ (for $\zeta<\kappa$) of $Y_{2i}$
such that for any disjoint subsets $X_1, X_2$ of $Y_{2i}$, each of cardinality
$<\sigma$, for some $\zeta<\kappa$ we have $F_{2i,\zeta}\cap(X_1\cup X_2)=X_1$
(possible as $\kappa=\kappa^{<\sigma}$ and $2^\kappa\geq\lambda_{2i+1}>
|Y_{2i}|$, by \cite{EK} or see \cite[Appendix 1.10]{Sh:g}). We can find a
sequence $\langle (j_i,\zeta_i):i<\delta\rangle$ such that
\begin{quotation}
\noindent $j_i\leq i<\delta$, $\zeta_i<\kappa$ and for an unbounded set of
$j<\delta$ for every $\zeta<\kappa$ for some $\vare<\kappa$ we have
\[\{2i,2i+1:\ 2i<\delta\ \&\ 2i,2i+1\in A_\vare\}\subseteq\{2i,2i+1<\delta:\
(j_i,\zeta_i)=(j,\zeta)\}\]
\end{quotation}
(we use: $\cf|\delta|=\cf(\mu)\leq\kappa$). Now for each $\alpha<\lambda$ we
define a set $Z_{\alpha}\subseteq Y$:

$f\in Z_{\alpha}$\qquad if and only if\qquad for some $i<\delta$
\begin{description}
\item[$(\alpha)$] \ \ \ $f\rest(2i)=f_\alpha\rest(2i)$,
\item[$(\beta)$]  \ \ \ $f\rest(2i+2)\neq f_\alpha\rest(2i+2)$,
\item[$(\gamma)$] \ \ \ $f\rest(2j_i)\in F_{2j_i,\zeta_i}\quad\Rightarrow\quad
f(2i)\leq f_\alpha(2i)\ \&\ f(2i+1)\leq f_\alpha(2i+1)$,
\item[$(\delta)$] \ \ \ $f\rest(2j_i)\notin F_{2j_i,\zeta_i}\quad\Rightarrow
\quad f(2i)\geq f_\alpha(2i)\ \&\ f(2i+1)\geq f_\alpha(2i+1)$.
\end{description}
Let $B$ be the $\sigma$-complete Boolean algebra of subsets of $Y$
generated by the family $\{Z_\alpha:\alpha<\lambda\}$. For $f\in Y$ let
$[f]=\{g\in Y:g$ extends $f\}$. For notational simplicity let $\sigma=
\aleph_0$.

\begin{claim}
\label{4.2B}
If for $\ell=1,2$, $Z^\ell\in B$ are the same Boolean combinations of
$Z_{\alpha^\ell_0},\ldots,Z_{\alpha^\ell_{n-1}}$, say $Z^\ell=\tau
[Z_{\alpha^\ell_0},\ldots,Z_{\alpha^\ell_{n-1}}]$ (where $\tau$ is a Boolean
term) and $i<\delta$ is such that $\langle f_{\alpha^\ell_m}\rest (2i): m<n
\rangle$ is with no repetition and $(\forall m<n)(f_{\alpha^1_m}\rest(2i)=
f_{\alpha^2_m}\rest (2i))$ then:
\begin{description}
\item[(a)] $Z^1\setminus\bigcup\limits_{m<n}[f_{\alpha^1_m}\rest(2i)]=
Z^2\setminus\bigcup\limits_{m<n} [f_{\alpha^2_m}\rest (2i)]$,
\item[(b)] for each $m<n$:
\begin{description}
\item[($\alpha$)] either $Z^1\cap [f_{\alpha^1_m}\rest (2i)]=[f_{\alpha^1_m}
\rest (2i)]\cap Z_{\alpha^1_m}$ and

$Z^2\cap [f_{\alpha^2_m}\rest (2i)]=[f_{\alpha^2_m}\rest (2i)]\cap Z_{
\alpha^2_m}$,
\item[($\beta)$] or $Z^1\cap[f_{\alpha^1_m}\rest (2i)]=[f_{\alpha^1_m}
\rest(2i)]\setminus Z_{\alpha^1_m}$ and

$Z^2\cap[f_{\alpha^2_m}\rest(2i)]=[f_{\alpha^2_m}\rest(2i)]\setminus
Z_{\alpha^2_m}$.
\end{description}
\end{description}
\end{claim}

\noindent{\em Proof of the claim:}\hspace{0.2in} Check the definition of
$Z_\alpha$. \hfill$\square_{\ref{4.2B}}$
\medskip

Clearly $B$ is a Boolean algebra of cardinality $\lambda$. Now the proof of
``$B$ has no chain of cardinality $\lambda$'' is similar to the proof of
\ref{3.1}, \ref{3.2} noting that for each $i$:
\begin{description}
\item[$(\ast)$] if\qquad ($\circledcirc$) $\Gamma\subseteq\lambda_{2i}\times\lambda_{2i+1}$
and for arbitrarily large $\alpha<\lambda_{2i}$, for arbitrarily large $\beta
<\lambda_{2i+1}$ we have $(\alpha,\beta)\in\Gamma$

then\qquad we can find $(\alpha_1,\beta_1)\in\Gamma$, $(\alpha_2,\beta_2)
\in\Gamma$ such that $\alpha_1<\alpha_2\ \&\ \beta_1<\beta_2$.
\end{description}
Up to now, the use of the pairs $2i, 2i+1$ was not necessary. But in the proof
of ``$B$ has no pie of cardinality $\lambda$'', instead of $(*)$ we use:
\begin{description}
\item[$(**)$] if\qquad ($\circledcirc$) of $(*)$ holds

then\qquad we can find $(\alpha_1,\beta_1)\in\Gamma$ and $(\alpha_2,
\beta_2)\in \Gamma$ such that $\alpha_1<\alpha_2\ \&\ \beta_1>\beta_2$,
\end{description}
easily finishing the proof \QED$_{\ref{4.2}}$

\begin{conclusion}
\label{4.3}
For a class of cardinals $\lambda$, there is a Boolean algebra $B$ of
cardinality $\lambda^+$, with no chain and no pie of cardinality $\lambda^+$.
[We can say, in fact, that this holds for many $\lambda$.]
\end{conclusion}

\Proof For any regular $\kappa$, if $2^\kappa>\aleph_{\kappa^{+4}}$ then (by
\cite[\S 4]{Sh:400}) for some club $E$ of $\kappa^{+4}$,
\[\ga\subseteq (\kappa^{+5},\aleph_{\kappa^{+4}})\ \&\ |\ga|\leq\kappa\ \&\
\sup(\ga)<\aleph_{\delta}\ \&\ \delta\in E\quad \Rightarrow\quad\max\pcf(
\ga)<\aleph_\delta.\]
Next we choose by induction on $i<\kappa$ cardinals $\lambda_i\in\Reg
\cap [\kappa^{+4},\aleph_{\kappa^{+4}}]$, $\lambda_i>\max\pcf(\{\lambda_j:
j<i\})$. Let $\mu$ be minimal such that $\mu\geq\sup\limits_{i<\kappa}
\lambda_i$ and $\mu\in\pcf(\{\lambda_i: i\in\kappa\})$. Then (by
\cite[1.8]{Sh:345a}, replacing $\{\lambda_i: i<\kappa\}$ by a subset of the
same cardinality, noting $\{\lambda_i: i\in A\}\in J_{<{\mu}}[\{\lambda_i:
i<\kappa\}]$ when $A\in [\kappa]^{<\kappa}$) we have $\mu=\tcf(\prod\limits_{
i<\kappa}\lambda_i/J^{\rm bd}_{\kappa})$. Also, as $\mu\in [\kappa,\aleph_{
\kappa^{+4}}]$ is regular, it is a successor cardinal; now the conclusion
follows by \ref{4.2}. Generally, for any $\kappa$ let $\alpha_0=0$, $\lambda^0
=\kappa^+$. Choose by the induction on $n$, $\alpha_{n+1}$ and $\langle
\lambda_i: \alpha_n\leq i<\alpha_{n+1}\rangle$, $u_n$ and regular
$\lambda^{n+1}$ such that $\alpha_{n+1}=\alpha_n+\lambda^n$, $\langle
\lambda_i: \alpha_n\leq i<\alpha_{n+1}\rangle$ is a strictly increasing
sequence of regular cardinals in $[\lambda^n,\aleph_{(\lambda^n)^{+4}}]$ such
that $\lambda_i>\max\pcf(\{\lambda_j: \alpha_n\leq j<i\})$ (possible by
\cite[\S 2]{Sh:400} as above). Let $u_n\subseteq [\alpha_n,\alpha_{n+1})$,
$|u_n|=\lambda^n$ be such that $\prod\limits_{i\in u_n}\lambda_i/J^{\rm
bd}_{u_n}$ has a true cofinality which we call $\lambda^{n+1}$. Lastly, for
some infinite $v\subseteq\omega$, $\prod\limits_{n\in v}\lambda_n/J^{\rm
bd}_v$ has a true cofinality, which we call $\lambda$. By renaming
$v=\omega$, $u_n=[\alpha_n,\alpha_{n+1})$. Then $\delta=:\sup_n \alpha_n$,
$\lambda$, $\langle\lambda_i:i<\delta\rangle$ are as required in \ref{4.2},
if we let:
\[J=\{u\subseteq\delta:\mbox{ for every large enough }n,\ \ \sup(u\cap
\alpha_n)<\alpha_n\}.\]
One point is left: why is $\lambda$ a successor cardinal? Because it is in
$[\sup\limits_{n<\omega}\lambda_n,\prod\limits_{n<\omega}\lambda_n]$ and
either
\[\prod\limits_{n<\omega}\lambda_n\leq [\sup\limits_{n<\omega}\lambda_n]^{
\aleph_0}\leq 2^{\sum\limits_n\lambda_n}<\aleph_{(\sum\limits_{n<\omega}
\lambda_n)^{+4}},\]
or the first attempt succeeds for $\kappa=\sum\limits_{n<\omega}\lambda^n$.
\QED$_{\ref{4.3}}$
\medskip

We have actually proved the existence of many such objects. If we waive some
requirements, even more.

\begin{proposition}
\label{4.3A} For any regular cardinal $\theta$ we can find $\delta$, $J$,
$\lambda$, $\lambda_i$ (for $i<\delta$) as in \ref{4.2} and such that:
\begin{description}
\item[$(x)$] $\lambda$ is a successor cardinal,
\item[$(y)$] for each $i$ for some regular cardinal $\mu_i$ we have $\lambda_i
=\mu^+_i$ and $(\mu_i)^{\theta}=\mu_i$,
\item[$(z)$] one of the following occurs:
\begin{description}
\item[(i)] $\delta$ is a regular cardinal $<\lambda_0$, $\delta>0$ and
$J=J^{\rm bd}_\delta$,
\item[(ii)] $\delta=\beth_\delta$ has cofinality $\theta$, and for some
$\lambda^j$ ($j<\cf(\delta)$) we have:
\[J=\{a: a\subseteq\delta, \{j<\cf(\delta): a\cap\lambda^j\notin J^{\rm
bd}_{\lambda^j}\}\in J^{\rm bd}_{\cf(\delta)}\}\]
and $\mu_\alpha=(\mu_\alpha)^{\sup\{\lambda_i :i<\alpha\}}$. \QED
\end{description}
\end{description}
\end{proposition}

\begin{lemma}
\label{4.4}
Assume $\langle\lambda_i:i<\delta\rangle$, $\lambda$, $J$, $\sigma$, $\langle
f_\alpha:\alpha<\lambda\rangle$ are as in (a)--(d) of \ref{4.2}(1) and
\begin{description}
\item[(e')] for every $B\in J$ for some $i<\delta$ we have:\quad $2i\notin B
\ \&\ 2i+1\notin B$,
\item[(f')] for every $i<\delta$ we have $\ens_\sigma(\lambda_i,\lambda_i)$
or at least for some club $C$ of $\delta$, if $i<\delta$ and $i=\sup\{j\in C:
|i\cap C\setminus j|\geq 1\}$ then $\ens_\sigma(\lambda_{2i},|\{f_\alpha\rest
(2i): \alpha<\lambda\}|)$. 
\end{description}
Then the conclusion of \ref{4.2}(1) holds.\\
\relax [We can weaken (f') as in \ref{3.3}(6).]
\end{lemma}

\Proof For each $i<\delta$ let $\langle(\lambda_{2i},<^{2i}_\eta):\eta=
f_\alpha\rest 2i$ for some $\alpha<\lambda\rangle$ be a $\sigma$-entangled
sequence of linear orders (each of cardinality $\lambda_{2i})$.

Now repeat the proof of \ref{4.2} with no $F_{2i}$'s, but defining $Z_\alpha$
we let:

$f\in Z_\alpha$\qquad if and only if\qquad for some $i<\delta$, letting $j=2i$
or be as in clause (f') for $2i$ we have:
\[\begin{array}{l}
f\rest (2i)=f_\alpha\rest (2i),\quad f\rest(2i+2)\neq f_\alpha\rest
(2i+2),\quad\mbox{ and}\\
f(2i)\leq^{2i}_{f\rest(2i)} f_\alpha(2i)\quad\mbox{ and }f(2i+1)\leq^{2i}_{
f\rest(2i)}f_\alpha(2i+1).
\end{array}\]
\QED$_{\ref{4.4}}$
\medskip
    
\noindent{\sc Discussion:}\quad Now instead of using on each set
$\{\eta^\frown\!\langle\alpha\rangle:\ \alpha<\lambda_{\lh(\eta)}\}$ a linear
order we can use a partial order; we can combine \ref{4.5} below with
\ref{4.2}(2) or with any of our proofs involving pcf for the existence of
entangled linear order.

\begin{lemma}
\label{4.5}
\begin{enumerate}
\item Assume $\langle\lambda_i:i<\delta\rangle$, $\lambda$, $J$, $\sigma$,
$\langle f_\alpha:\alpha<\lambda\rangle$ are as in (a)--(d) of \ref{4.2}(1)
and
\begin{description}
\item[(e)] for each $i<\delta$ there is a sequence $\bar{P}=\langle P_\vare:
\vare<\kappa_i\rangle$ where $\kappa_i=|\{f_\alpha\rest i:\alpha<\lambda\}|$,
each $P_\vare$ is a partial order, $\bar{P}$ is $(\lambda,\sigma)$-entangled
which means:
\begin{description}
\item[\ ] if $u_0, u_1, u_2$ are disjoint subsets of $\kappa_i$ of cardinality
$<\sigma$ and for $\vare\in u_0\cup u_1\cup u_2$, $t^\vare_\alpha\in P_\vare$
(for $\alpha<\lambda$) are pairwise distinct

then for some $\alpha<\beta$:
\[\begin{array}{lcl}
\vare\in u_0 & \Rightarrow & P_\vare\models t^\vare_\alpha<t^\vare_\beta,\\
\vare\in u_1 & \Rightarrow & P_\vare\models t^\vare_\alpha>t^\vare_\beta,\\
\vare\in u_2 & \Rightarrow & P_\vare\models \mbox{``}t^\vare_\alpha,
t^\vare_\beta\mbox{ are incomparable''.}\\
\end{array}\]
\end{description}
\end{description}
Then the conclusion of \ref{4.2}(1) holds.
\item Assume as in (1) but
\begin{description}
\item[(e')] for some $A\subseteq\delta$, $A\notin J$, $\delta\setminus A
\notin J$ and 
\item[(e')$_1$] like (e) of \ref{4.5}(1) for $i\in A$ with
$u_2=\emptyset$,
\item[(e')$_2$] like (e) of \ref{4.5}(1) for $i\in\delta\setminus A$, with
$u_0=u_1=\emptyset$, $u_2=1$ (so we can use $P_\alpha=(\lambda_i,=))$.
\end{description}
Then the conclusion of \ref{4.2}(1) holds.
\item We can weaken ``$\kappa_i=|\{f_\alpha\rest i:\alpha<\lambda\}|$'' as in
\ref{4.4}(f').
\end{enumerate}
\end{lemma}

\Proof Similar to earlier ones.  \QED$_{\ref{4.5}}$

\section{More on Entangledness}

\begin{proposition}
\label{5.1}
Suppose that $\langle\lambda_i:i<i(*)\rangle$ is a strictly increasing
sequence of regular cardinals, $T_i\subseteq {}^{\lambda_i>} 2$ is closed
under initial segments, $i+1<i(*)\ \Rightarrow\ |T_i|<\lambda_{i+1}$ and
the set
\[B_i=\{\eta\in {}^{\lambda_i} 2:\mbox{ for every }\alpha<\lambda_i,
\eta\rest\alpha\in T_i\}\]
has cardinality $\geq\mu=\cf(\mu)>\lambda_i+|T_i|$ (for each $i<i(*))$.\\
Then $\langle(B_i,<_{\ell x}):i<i(*)\rangle$ is a $(\mu,\aleph_0)$-entangled
sequence of linear orders $(<_{\ell x}$ is the lexicographic order).
\end{proposition}

\begin{remark}
\label{5.1A}
{\em
So if $\mu=\cf(\mu)$, $\theta=|\{\lambda:\lambda<\mu\leq 2^\lambda\mbox{ and
}2^{<\lambda}<2^\lambda\}|$ then $\ens(\mu,\theta)$, see \cite[3.4]{Sh:430}.
}
\end{remark}

\Proof Clearly $|T_i|\geq\lambda_i$ (as $B_i\neq\emptyset$). So let $n<
\omega$, $i_0<i_1<\ldots<i_{n-1}<i(*)$, and $\eta^\ell_\zeta\in B_{i_\ell}$
for $\ell<n$, $\zeta<\mu$ be such that:
\[\zeta<\xi<\mu\ \&\ \ell<n\quad\Rightarrow\quad\eta^\ell_\zeta\neq
\eta^\ell_\xi,\]
and let $u\subseteq n$. We should find $\zeta<\xi<\mu$ such that
$(\forall\ell<n)(\eta^\ell_\zeta<_{\ell x}\eta^\ell_\xi\ \Leftrightarrow\
i\in u)$. To this end we prove by downward induction on $m\leq n$ that
(stipulating $\lambda_n=\mu)$:
\begin{description}
\item[$(*)_m$] there is a set $w\subseteq\mu$ of cardinality $\geq\lambda_{
i_m}$ such that:

if $m\leq\ell<n$ and $\zeta<\xi$ are from $w$ then
$\left[\eta^\ell_\zeta <_{\ell x}\eta^\ell_\xi\right]^{\mbox{if }(\ell\in u)}$.
\end{description}
Note that $(*)_n$ is exemplified by $w=:\mu$ and $(*)_0$ says (more than) that
the conclusion holds, so this suffices. Hence assume $(*)_{m+1}$ is
exemplified by $w^*$ and we shall find $w\subseteq w^*$ exemplifying $(*)_m$,
$|w|=\lambda_{i_m}$. Without loss of generality $m\in u$ (otherwise replace
each $\eta\in T_{i_m}\cup\{\eta^m_\zeta: \zeta\in w^*\}$ by $\langle 1-\eta(
\alpha): \alpha<\lh(\eta)\rangle)$. Let for $\alpha<\lambda_{i_m}$ and $\nu
\in T_i\cap{}^\alpha 2$:
\[w^*_\nu=:\{\zeta\in w^*:\nu=\eta^m_\zeta\rest\alpha,\ \neg(\exists\xi\in
w^*)(\xi<\zeta\ \&\ \eta^m_\xi\rest\alpha=\nu\ \&\ \eta^m_\xi<_{\ell x}
\eta^m_\zeta)\},\]
\[w^*_\alpha=\bigcup\{w^*_\nu:\nu\in T_{i_m}\cap{}^\alpha 2\}.\]
As in $(B_{i_m},<_{\ell x})$ there is no monotonic sequence of length
$\lambda^+_{i_m}$, clearly $|w^*_\nu|\leq\lambda_{i_m}$. Moreover,
\[|w^*_\alpha|\leq|\{\nu\in T_{i_m}:\lh(\nu)=\alpha\}|\times\sup\{|w^*_\nu|:
\nu\in T_{i_m}\cap{}^\alpha 2\}\leq |T_{i_m}|\times\lambda_{i_m},\]
and hence $|\bigcup\limits_{\alpha<\lambda_{i_m}} w^*_\alpha |\leq\lambda_{
i_m}+|T_{i_m}|$. But $|T_{i_m}|<\lambda_{i_{m +1}}$ and $\lambda_{i_m}<
\lambda_{i_{m+1}}$. Hence we find $\zeta(*)\in w^*\setminus\bigcup\limits_{
\alpha<\lambda_{i_m}} w^*_\alpha$. Now, for every $\alpha<\lambda_{i_m}$ let
$\xi_\alpha\in w^*$ exemplify $\zeta(*)\notin w^*_{\eta^m_\zeta(*)\rest
\alpha}\subseteq w^*_\alpha$, so
\[\xi_\alpha<\zeta(*),\quad\eta^m_{\xi_{\alpha}}\rest\alpha=\eta^m_{\zeta(*)}
\rest\alpha\ \mbox{ and }\ \eta^m_{\xi_\alpha}<_{\ell x}\eta^m_{\zeta(*)}.\]
Then some $\gamma_\alpha$, $\alpha\leq\gamma_\alpha<\lambda_{i_m}$, we have
\[\eta^m_{\zeta(*)}\rest\gamma_\alpha=\eta^m_{\xi_\alpha}\rest\gamma_\alpha,
\quad \eta^m_{\zeta(*)}(\gamma_\alpha)=1,\quad \eta^m_{\xi_\alpha}(
\gamma_\alpha)=0.\]
So for some unbounded set $A\subseteq\lambda_i$ the sequence $\langle
\gamma_\alpha:\alpha\in A\rangle$ is strictly increasing in $\alpha$ and also
$\langle \xi_\alpha: \alpha\in A\rangle$ is increasing. Let $w=:\{\xi_\alpha:
\alpha\in A\}\subseteq w^*$. It exemplifies $(*)_m$, hence we finish.
\QED$_{\ref{5.1}}$

\begin{proposition}
\label{5.1B}
\begin{enumerate}
\item Assume that
\begin{description}
\item[(a)] $\lambda=\max\pcf(\ga_\vare)$ for $\vare<\vare(*)$,
\item[(b)] $|\ga_\vare|\leq\kappa<\kappa^*\leq\min(\ga_\vare)$,
\item[(c)] $\theta\in\ga_\vare\ \ \Rightarrow\ \ \theta$ is $(\kappa^*,
\kappa^+,2)$-inaccessible,
\item[(d)] for $n<\omega$ and distinct $\vare_0,\vare_1,\ldots,\vare_n<
\vare(*)$ we have
\[\ga_{\vare_0}\setminus\bigcup_{\ell=1}^n \ga_{\vare_\ell}\notin J_{<
\lambda}[\ga_{\vare_0}].\]
\end{description}
Then $\ens(\lambda,\lambda,\vare(*))$.
\item Assume in addition that
\begin{description}
\item[(d)] if $u\in [\vare(*)]^{<\sigma}$, $\vare\in\vare(*)\setminus u$
then $\ga_\vare\setminus\bigcup\limits_{\zeta\in u}\ga_\zeta\notin J_{<
\lambda}[\ga_\vare]$,
\item[(e)] if $\theta\in\ga_\vare$ then $\max\pcf(\bigcup\limits_{\zeta}
\ga_\zeta\cap\theta)<\theta$,
\item[(f)] $(\forall\alpha<\lambda)(|\alpha|^{<\sigma}<\lambda)$.
\end{description}
Then $\ens_\sigma(\lambda,\lambda,\vare(*))$.
\end{enumerate}
\end{proposition}

\Proof As in \S 3. \QED$_{\ref{5.1B}}$

\begin{proposition}
\label{5.2}
For any cardinal $\lambda$ satisfying $(\forall\kappa<\lambda)(2^\kappa<
2^\lambda)$ there is a successor cardinal $\theta\in [\lambda,2^\lambda]$
such that there is an entangled linear order of cardinality $\theta$.
\end{proposition}

\Proof We prove slightly more, so let $\lambda$ and $\chi\in [\lambda,
2^\lambda]$ be any cardinals (we shall try to find $\theta^+\in [\chi,
2^\lambda]$ for $\chi$ as below; for the proposition $\mu=\lambda$ below).\\
Let $\mu=:\min\{\mu:2^\mu=2^\lambda\}$, so $\mu\leq\lambda$ and $\mu<\cf(
2^\mu)$ and $\kappa<\mu\ \Rightarrow\ 2^\kappa<2^\mu$. First assume
$2^{<\mu}= 2^\mu$. Then necessarily $\mu$ is a limit cardinal. If $\cf(2^{
<\mu})=\cf(\mu)$ we get a contradiction to the previous sentence. Hence
$\langle 2^\theta: \theta<\mu\rangle$ is eventually constant so for some
$\theta<\mu$ we have $2^\theta=2^{<\mu}$ but $2^{<\mu}=2^\mu$, a contradiction
to the choice of $\mu$. Thus we have $2^{<\mu}<2^\mu= 2^\lambda$. Assume
$\chi=\lambda+2^{<\mu}$ or just $2^\lambda>\chi\geq 2^{<\mu}$. The proof
splits to cases:\quad if $\cf(2^\lambda)$ is a successor, use cases B or C or
D, if $\cf(2^\lambda)$ is a limit cardinal (necessarily $>\lambda)$ use case
A.
\medskip

\noindent{\sc Case A:}\qquad $\chi^{+(\mu^{+4})}\leq 2^\lambda$.

\noindent By \ref{3.5}(2).
\medskip

\noindent{\sc Case B:}\qquad $\cf(2^\lambda)$ is a successor, $\mu$ is strong
limit (e.g.~$\aleph_0$).

\noindent Clearly there is a dense linear order of cardinality
$\cf(2^\lambda)$ and density $\mu$, hence there is an entangled linear order
in $\cf(2^\lambda)$, which is as required (\cite{BoSh:210}).
\medskip

\noindent{\sc Case C:}\qquad $\cf(2^\lambda)$ is a successor cardinal,
$\mu$ is regular uncountable.

\noindent Look at \cite[4.3]{Sh:410} (with $\mu$ here standing for $\lambda$
there); conditions (i) + (ii) hold. Now on the one hand, if the assumption
(iii) of \cite[4.3]{Sh:410} fails, we know that there is an entangled linear
order of cardinality $\cf(2^\lambda)$ (as in Case B). But on the other hand,
if (iii) holds, the conclusion of \cite[4.3]{Sh:410} gives more than we asked
for. In both cases there is an entangled linear order in $\cf(2^\lambda)$
(which is a successor).
\medskip

\noindent{\sc Case D:}\qquad $\mu$ is singular, not strong limit.

\noindent By \cite[3.4]{Sh:430} there are regular cardinals $\theta_i$ (for
$i<\cf(\mu)$) such that the sequence $\langle\theta_i:i<\cf(\mu)\rangle$ is
increasing with limit $\mu$, $\mu<2^{\theta_i}$, $\langle 2^{\theta_i}:
i<\cf(\mu)\rangle$ is strictly increasing and
\begin{quotation}
\noindent for all $\sigma$ such that $\sigma=\cf(\sigma)\leq 2^{\theta_i}$
there is a tree $T^i_\sigma$, $|T^i_\sigma|=\theta_i$ or at least $2^{|
T^i_\sigma|}=2^{\theta_i}$ and $T^i_\sigma$ has $\geq\sigma$\ \
$\theta_i$-branches.
\end{quotation}
If for some $i$ either $2^{\theta_i}$ is a successor $\geq\chi$ or $\cf(2^{
\theta_i})$ is a successor $\geq\chi$ we finish as in Case B, if $\cf(2^{
\theta_i})$ is a limit cardinal, we finish as in Case A (or use
\ref{5.2A}). \QED$_{\ref{5.2}}$

\begin{proposition}
\label{5.2A}
Assume $2^\mu$ is singular. Then there is an entangled linear order of
cardinality $(2^\mu)^+$.
\end{proposition}

\Proof Let $\lambda$ be the first singular cardinal $>\mu$ such that
$(\exists\kappa<\lambda)(\pp_\kappa(\lambda)>2^\mu)$. Now,
$\lambda$ is well defined, and moreover $\lambda\leq 2^\mu$ (as $2^\mu$ is
singular so $\pp(2^\mu)\geq (2^\mu)^+>2^\mu$) and (by
\cite[2.3]{Sh:355}, \cite[1.9]{Sh:371})
\[\cf(\lambda)<\chi\in (\mu,\lambda)\setminus\Reg\quad \Rightarrow
\quad \pp_{\cf(\lambda)+\cf(\chi)}(\chi)<\lambda,\]
so $\pp(\lambda)>2^\mu$ and $\cf(\lambda)>\mu$ (otherwise $\pp(
\lambda)\leq\lambda^{\cf(\lambda)}\leq\lambda^\mu\leq (2^\mu)^\mu=2^\mu$).
Lastly apply \ref{3.2}(1). \QED$_{\ref{5.2A}}$

\begin{proposition}
\label{5.3}
If $\kappa^+<\chi_0<\lambda$ and $\kappa^{+4}<\cf(\lambda)<\lambda\leq
2^\kappa$, then (a) or (b) holds:
\begin{description}
\item[(a)] there is a strictly increasing sequence $\langle\lambda^*_i:
i<\delta\rangle$ of regular cardinals from $(\chi,\lambda)$, $\delta=\cf(
\delta)\in [\kappa,\cf(\lambda)]\cap\Reg$ and $\lambda^*_i>\max\pcf(
\{\lambda^*_j: j<i\})$, such that $\lambda^+=\tcf\prod\limits_{i<\delta}
\lambda^*_i/J^{\rm bd}_\delta$,
\item[(b)] there is a strictly increasing sequence $\langle\lambda^*_i:
i<\delta\rangle$ of regular cardinals from $(\chi,\lambda)$ such that
$\lambda^*_i>\max\pcf(\{\lambda^*_j: j<i\})$ and $\ens(\lambda^*_i,
\lambda^*_i)$ and $\lambda^+=\tcf\prod\limits_{i<\delta}\lambda^*_i/I$, $I$
a proper ideal on $\delta$ extending $J^{\rm bd}_\delta$.

Moreover for some $\mu\in(\chi_0,\lambda)$, $\mu<\lambda_0^*$ and there
is a sequence $\langle\gb_{i,j}: i<\delta, j<\kappa^+\rangle$ such that
$\gb_{i,j}\subseteq\Reg\cap\mu\setminus\chi_0$, $|\gb_{i,j}|\leq\kappa$,
each $\theta\in\bigcup\limits_{i,j}\gb_{i,j}$ is $(\chi_0,\kappa^+,
\aleph_0)$-inaccessible (i.e.~$\ga\subseteq\Reg\cap\theta\setminus
\chi_0,\ |\ga|\leq\kappa\ \ \Rightarrow\ \ \max\pcf(\ga)<\theta$) and
\[j_1<j_2<\kappa^+\ \ \Rightarrow\ \ \gb_{i,j_1}\cap\gb_{i,j_2}=\emptyset,
\quad \lambda^*_i=\max\pcf(\gb_{i,j}),\]
\[\mu=\sup(\gb_{i,j}),\quad J^{\rm bd}_{\gb_{i,j}}\subseteq J_{<\lambda^*_i}
[\gb_{i,j}].\]
(This implies $\ens(\lambda^*_i,\lambda^*_i)$.)
\end{description}
\end{proposition}

\begin{remark}
\label{5.3A}
{\em
\begin{enumerate}
\item Why $\lambda^+$ instead of $\lambda^*=\cf(\lambda^*)\in (\lambda,
\pp^+_{J^{\rm bd}_{\cf(\lambda)}}(\lambda))$? To be able to apply
\cite[3.3]{Sh:410} in case III of the proof of \cite[4.1]{Sh:410}. So if
$\langle\lambda_i:i<\cf(\lambda)\rangle$ fits in such a theorem we can get
$\lambda^+$.
\item We could have improved the theorem if we knew that always
\[\cf([\lambda]^{\textstyle\leq\kappa},\subseteq)=\lambda+\sup\{\pp_\kappa
(\mu):\cf(\mu)\leq\kappa<\mu<\lambda\},\]
particularly  $\sigma$-entangledness.
\end{enumerate}
}
\end{remark}

\Proof This is like the proof of \cite[4.1]{Sh:410}. (In case II when
$\sigma=\aleph_0$ imitate \cite[4.1]{Sh:410}.) However, after many doubts,
for reader's convenience we present the proof fully, adopting for our
purposes the proof of \cite[4.1]{Sh:410}.

By \cite[2.1]{Sh:355} there is an increasing continuous sequence $\langle
\lambda_i:i<\cf(\lambda)\rangle$ of singular cardinals with limit $\lambda$
such that $\tcf(\prod\limits_{i<\cf(\lambda)}\lambda_i^+,<_{J^{\rm bd}_{\cf(
\lambda)}}) =\lambda^+$ and $\lambda_0>\chi_0$. The proof will split to cases.
Wlog $\chi_0>\cf(\lambda)$.
\medskip

\noindent{\sc Case I:}\qquad $\max\pcf(\{\lambda_j^+:j<i\})<\lambda$ for $i
<\cf(\lambda)$.

\noindent So for some unbounded $A\subseteq\cf(\lambda)$ we have
\[(\forall i\in A)(\max\pcf(\{\lambda_j^+:j\in A\cap i\})<\lambda_i^+).\]
Consequently $\ga=\{\lambda_i^+: i\in A\}$ satisfies the demands of \ref{3.5}
and hence (a) holds true with $\delta=\cf(\lambda)$, $\lambda^*_i=\lambda_i$.
\medskip

\noindent Thus assume that Case I fails. So there is $\mu$ such that
$\chi_0<\mu<\lambda$, $\cf(\mu)<\cf(\lambda)$ and $\pp_{<\cf(\lambda)}(\mu)
>\lambda$. Choose a minimal such $\mu$. Then, by \cite[3.2]{Sh:410}, we have:
\begin{description}
\item[$(*)$] $[\ga\subseteq\Reg\setminus\chi_0\ \&\ \sup\ga<\mu\ \&\
|\ga|<\cf(\lambda)]\ \ \Rightarrow\ \ \max\pcf(\ga)<\lambda$.
\end{description}
By \cite[2.3]{Sh:355} in the conclusion of $(*)$ we may replace ``$<\lambda$''
by ``$<\mu$'' and we get
\begin{description}
\item[$(*)'$] $[\ga\subseteq\Reg\setminus\chi_0\ \&\ \sup\ga<\mu\ \&\
|\ga|<\cf(\lambda)]\ \ \Rightarrow\ \ \max\pcf(\ga)<\mu$.
\end{description}
Let $\sigma=\cf(\mu)$. Then $\pp(\mu)=\pp_{<\cf(\lambda)}(\mu)$ (and $\pp_{
<\cf(\lambda)}(\mu)>\lambda$). Wlog $\mu<\lambda_0$.
\medskip

\noindent{\sc Case II:}\qquad $\sigma>\kappa$ (and not Case I).

\noindent By \cite[1.7]{Sh:371}, if $\sigma>\aleph_0$ and by
\cite[6.x]{Sh:430} if $\sigma=\aleph_0$ we find a strictly increasing sequence
$\langle\mu_i^*:i<\sigma\rangle$ of regular cardinals, $\mu=\bigcup\limits_{i
<\sigma}\mu^*_i$ and an ideal $J$ on $\sigma$ extending $J^{\rm bd}_\sigma$
(if $\sigma>\aleph_0$ then $J=J^{\rm bd}_\sigma$) such that
\[\lambda^+=\max\pcf(\{\mu^*_i:i<\sigma\})=\tcf\prod\limits_{i<\sigma}
\mu^*_i/J.\]
If $\sigma>\aleph_0$, since we may replace $\langle\mu^*_i: i<\sigma\rangle$
by $\langle\mu^*_i:i\in A\rangle$ for any unbounded $A\subseteq\sigma$, we may
assume that $\mu^*_i>\max\pcf(\{\mu^*_j:j<i\})$. If $\sigma=\aleph_0$ this
holds automatically, so in both cases we can apply \ref{3.5}. So we get
(a), as $\sigma>\kappa$.
\medskip

\noindent{\sc Case III:}\qquad $\sigma\leq\kappa$.

\noindent So $\sigma^{+4}\leq\cf(\lambda)$. Let
\[\begin{array}{ll}
\cP=:\{C\subseteq\cf(\lambda):&\otp(C)=\kappa^{+3},\ \ C\mbox{ is closed in }
\sup(C)\ \mbox{ and}\\
\ &\max\pcf(\{\lambda^+_i: i\in C\})<\lambda\}.
\end{array}\]
[Why there are such $C$'s? For any $\delta<\cf(\lambda)$, $\cf(\delta)=
\kappa^{+3}$ we have a club $C'$ of $\lambda_\delta$ such that $\tcf(\prod
\limits_{\kappa\in C'}\kappa^+/J^{\bd}_{\kappa^{+3}})=\lambda_\delta^+$. Now
$C'\cap \langle\lambda_i:i<\delta\rangle$ will do.]\\
For each $C\in \cP$ try to choose by induction on $i<\kappa^+$,
$\gb_i=\gb_{i,C}$ and $\gamma_i=\gamma_{i,C}$ such that:
\begin{description}
\item[(i)]   $\gb_i\subseteq\Reg\cap\mu\setminus\bigcup\limits_{j<i}
\gb_j\setminus\chi_0$,
\item[(ii)]  $\gamma_i\in C\setminus\bigcup\limits_{j<i}(\gamma_j+1)$,
\item[(iii)] $\lambda^+_{\gamma_i}\in\pcf(\gb_i)$,
\item[(iv)]  $|\gb_i|\leq\sigma$,
\item[(v)]   all members of $\gb_i$ are $(\chi_0,\kappa^+,
\aleph_0)$-inaccessible,
\item[(vi)]  $\gamma_i$ is minimal under the other requirements.
\end{description}
Let $(\gb_{i,C},\gamma_{i,C})$ be defined if and only if $i<i_C(*)$. So
{\em success} in defining means $i_C(*)=\kappa^+$, {\em failure} means
$i_C(*)<\kappa^+$.
\smallskip

\noindent{\bf Subcase IIIA}:\qquad For some $j<\kappa^+$, for every $C\in
\cP$ with $\min(C)\geq j$ we have $i_C(*)<\kappa^+$, so we cannot define
$\gb_{i_C(*),C}$, $\gamma_{i_C(*),C}$.

\noindent Let $C, i_C(*)$ be as above. Let $\gamma^*_C=\bigcup\limits_{i<
i_C(*)}\gamma_{i,C}$, so $\gamma^*_C\in C$. Now, if $\gamma\in C\setminus
\gamma^*_C$ then (by \cite[1.5B]{Sh:355}) as $\pp_\sigma(\mu)\geq\lambda^+>
\lambda^+_\gamma$, there is $\ga_\gamma\subseteq\Reg\cap(\chi,\mu)$,
$|\ga_\gamma|\leq\sigma$ such that $\lambda^+_\gamma\in\pcf(\ga_\gamma)$. By
\cite[3.2]{Sh:410} there is $\gc_\gamma\subseteq\Reg\cap(\chi,\mu)$ of
cardinality $\leq\kappa$ consisting of $(\chi,\kappa^+,\aleph_0)$-inaccessible
cardinals such that $\lambda^+_\gamma\in\pcf(\gc_\gamma)$. Now $\gamma$,
$\gc_\gamma\setminus\bigcup\limits_{i<i_C(*)}\gb_{i,C}$ cannot serve as
$\gamma_{i_C(*),C}$, $\gb_{i_C(*),C}$, so necessarily $\lambda^+_\gamma\notin
\pcf(\gc_\gamma\setminus\bigcup\limits_{i<i_C(*)}\gb_{i,C})$. Hence wlog
$\gc_\gamma\subseteq\bigcup\limits_{i<i_C(*)}\gb_{i,C}$. So
\[\{\lambda^+_i: i\in C\setminus \gamma^*_C\}\subseteq\pcf(\bigcup\limits_{i
<i_C(*)}\gb_{i,C})\ \ \mbox{ and }\ \ |\bigcup\limits_{i<i_C(*)}\gb_{i,C}|\leq
\kappa.\]
By the proof of \cite[4.2]{Sh:400} (or see \cite[\S 3]{Sh:410}) we get a
contradiction (note that $\cf(\lambda)>\kappa^{+4}$ does not disturb).
\smallskip

\noindent{\bf Subcase IIIB}:\qquad For every $j<\cf(\lambda)$ there is $C\in
\cP$ such that $\min(C)>j$ and $i_C(*)=\kappa^+$, i.e.~$\gb_{i,C}$,
$\gamma_{i,C}$ are defined for every $i<\kappa^+$.

\noindent We will show
\begin{description}
\item[($\otimes$)] for every $j(*)<\cf(\lambda)$ there is $\lambda'\in\lambda\cap
\pcf(\{\lambda^+_j:j<\cf(\lambda)\})\setminus\lambda_{j(*)}$ such that
\begin{description}
\item[$(\alpha)$] $\ens(\lambda',\lambda')$ (exemplified by linear order which
has density character $>\chi_0$ in every interval),
\item[$(\beta)$]  for some $\gb\subseteq\mu\cap\Reg\setminus\chi_0$ we
have:\quad $|\gb|\leq\kappa^+$, $\lambda'=\max\pcf(\gb)$. Moreover, $\gb$ can
be divided to $\kappa^+$ subsets of cardinality $\leq\kappa$, no one in
$J_{<\lambda'}[\gb]$ and
\[(\forall\theta\in\gb)(\theta>\max\pcf(\gb\cap\theta))\]
(even $\theta$ is $(\chi_0,\kappa^+,\aleph_0)$-inaccessible).
\end{description}
\end{description}
\smallskip

{\em Why does ($\otimes$) suffice?}

\noindent Suppose that we have proved ($\otimes$) already. So for $i<\cf(
\lambda)$ we can choose $\mu^*_i$, $\lambda_i<\mu^*_i=\cf(\mu^*_i)\in\lambda
\cap\pcf(\{\lambda^+_j: j<\cf(\lambda)\})$ as required in ($\otimes$). Since
$(\forall i)(\mu^*_i<\lambda)$, wlog the sequence $\langle\mu^*_i: i<\cf(
\lambda)\rangle$ is strictly increasing. By induction on $\vare<\cf(\lambda)$
choose strictly increasing $i(\vare)<\cf(\lambda)$ such that $\mu^*_{i(
\vare)}>\max\pcf(\{\mu^*_{i(\zeta)}:\zeta<\vare\})$.

\noindent Let $i(\vare)$ be defined if and only if $\vare<\vare(*)$. So
$\vare(*)$ is limit,
\[\lambda^+=\max\pcf(\{\mu^*_{i(\vare)}:\vare<\vare(*)\}),\quad\mbox{
and }\quad\mu^*_{i(\vare)}>\max\pcf(\{\mu^*_{i(\zeta)}:\zeta<\vare\}),\]
$\mu^*_{i(\vare)}$ is strictly increasing and $\ens(\mu^*_{i(\vare)},
\mu^*_{i(\vare)})$. Thus applying \cite[4.12]{Sh:355} we finish, getting
clause (b) of \ref{5.3}.
\smallskip

{\em Why does ($\otimes$) hold?}

\noindent Choose $C\subseteq(j(*),\cf(\lambda))$ of order type $\kappa$ such
that $\langle\gamma_i:i<\kappa^+\rangle$, $\langle\gb_i:i<\kappa^+\rangle$ are
well defined and
\[\max\pcf(\{\lambda^+_{\gamma_i}: i\in C\})<\lambda\]
(possible by our being in subcase IIIB, see the definition of $\cP$). Let
$\gd=:\{\lambda^+_{\gamma_i}:i<\kappa^+\}$ and let $\langle\gb_\theta[\gd]:
\theta\in\pcf(\gd)\rangle$ be as in \cite[2.6]{Sh:371}. Let $\theta$ be
minimal such that $\otp(\gb_\theta[\gd])=\kappa$. We can find pairwise
disjoint sets $B_\vare\subseteq C$ (for $\vare<\kappa$) such that
\[\{\lambda^+_\gamma:\gamma\in B_\vare\}\subseteq\gb_\theta[\gd],\quad\quad
\otp(B_\vare)=\kappa.\]
Clearly $\max\pcf(\{\lambda^+_{\gamma_i}:\gamma_i\in B_\vare\})=\theta$, since
$\{\lambda^+_{\gamma_i}:\gamma_i\in B_\vare\}\subseteq \gb_\theta[\gd]$ but
it is not a subset of any finite union of $\gb_{\theta'}[\gc]$, $\theta'<
\theta$. Now letting $\ga^*=:\bigcup\limits_{j\in C}\gb_j$, we find (by
\cite[2.6]{Sh:371}) a subset $\ga$ of $\ga^*$ such that $\theta=\max\pcf(\ga)$
but $\theta\notin\pcf(\ga^*\setminus\ga)$. Now as $\theta\in\pcf(\{
\lambda^+_\gamma:\gamma\in B_\vare\})$, $\lambda^+_\gamma\in\pcf(\gb_\gamma)$
we have (by \cite[1.12]{Sh:345a}) $\theta\in\pcf(\bigcup\limits_{\gamma\in
B_\vare}\gb_\gamma)$. Hence by the previous sentence $\theta\in\pcf(\ga\cap
\bigcup\limits_{\gamma\in B_\vare}\gb_\gamma)$. Let
\[\gc_\vare=: \ga\cap\bigcup\limits_{j\in B_\vare}\gb_j, \quad\quad
\lambda'=\theta.\]
We can apply \ref{3.2} and get that there is an entangled linear order of
cardinality $\lambda'$ (which is more than required) and, of course,
\[\lambda_{j(*)}<\lambda'\in\lambda\cap\pcf(\{\lambda_j: j<\cf(\lambda)\}).\]
The assumptions of \ref{3.2} hold as the $\gc_\vare$ are pairwise
disjoint (by (i) above),
\[\theta\in\pcf(\{\lambda^+_{\gamma_i}:\gamma_i\in B_\vare\}), \quad\quad
\pcf(\bigcup\limits_{j\in B_\vare}\gb_j)=\pcf(\gc_\vare)\qquad\mbox{
and}\]
\[\theta_1\in\ga\ \ \Rightarrow\ \ \max\pcf(\ga\cap\theta_1)<\theta_1,\]
as $\theta_1$ is $(\chi_0,\kappa^+,\aleph_0)$-inaccessible and
\[\theta=\lambda'\geq\sup\{\lambda^+_{\gamma_i}: i\in C\}>\lambda_{j(*)}>
\chi_0.\]
So clause $(\alpha)$ of ($\otimes$) holds and clause $(\beta)$ was done along
the way. Thus we finish subcase IIIb and hence case III.
\QED$_{\ref{5.3}}$

\begin{conclusion}
\label{5.4}
For $\lambda$ as in \ref{5.3} there is a Boolean algebra $B$ of cardinality
$\lambda^+$ satisfying $(\oplus)^B_{\lambda^+}$, $(\otimes)^B_{\sup\limits_{
i<\delta}\lambda^*_i}$ from \ref{4.1} (and also there is an entangled linear
order in $\lambda^+$).
\end{conclusion}

\Proof If (a) of \ref{5.3} holds, apply \ref{4.2}. If (b) of \ref{5.3} holds
use \ref{4.4}. \QED$_{\ref{5.4}}$

\begin{definition}
\label{5.6}
\begin{enumerate}
\item $\pcf^{ex}_\kappa(\ga)=\{\lambda:\mbox{ if }\gb\subseteq\ga,\|\gb|<
\kappa$ then $\lambda\in \pcf(\ga\setminus\gb)\}$ (equivalently: $\lambda
\in\pcf(\ga)$ and $\gb\in J_{<\lambda}[\ga]\ \ \Rightarrow\ \ |\gb_\lambda[
\ga]\setminus\gb| \geq \kappa$).
\item  $J^{ex,\kappa}_{<\theta}[\ga]=:\{\gb:\gb\subseteq\ga$ and for some
$\gc\subseteq\ga$ we have: $|\gc|<\kappa$ and $\gb\setminus\gc\in J_{<\theta}
[\ga]\}$.
\end{enumerate}
\end{definition}

\begin{proposition}
\label{5.5}
Assume that $\mu$ is a singular cardinal which is a fix point (i.e.~$\mu=
\aleph_\mu$) and $\mu^*<\mu$.
\begin{enumerate}
\smallskip
\item For some successor cardinal $\lambda^+\in (\mu, \pp^+(\mu))$ there is an
entangled linear $\cI$ order of cardinality $\lambda^+$ and density
$\in(\mu^*,\mu]$.
\item If $\mu\leq\chi_0$, $\chi^{+\mu^{+4}}_0\leq \pp^+(\mu)$ then we can find
an entangled linear order $\cI$, $|\cI|=\lambda^+\in (\chi_0,\chi_0^{+
\mu^{+4}})$ of density $\in (\mu^*,\mu]$.
\item In both parts we get also a Boolean Algebra $B$ satisfying
$(\oplus)^B_{\lambda^+}$, $(\otimes)^B_\mu$ of \ref{4.2}.
\item In both parts 1), 2),  if $\cJ$ is an interval of $\cI$ or $\cJ\in
[\cI]^{\lambda^+}$ then $\dens(\cJ)=\dens(\cI)$. This applies to
\ref{5.5A}, \ref{5.5B}, \ref{5.5C}, too.
\end{enumerate}
\end{proposition}

\Proof Let $\langle\mu_i:1\leq i<\cf(\mu)\rangle$ be a strictly increasing
continuous sequence with limit $\mu$. Wlog $\mu_1>\mu_0>\mu^* +\cf(\mu)$,
$\aleph_0\leq\cf(\mu_i)<\max\{\cf(\mu),\aleph_1\}$.

\noindent 1)\ \ \ We try to choose by induction on $i<\cf(\mu)$ regular
cardinals $\lambda_i$ such that
\[\mu_i<\lambda_i<\mu,\quad \max\pcf(\{\lambda_j:j<i\})<\lambda_i,\]
and there is an entangled sequence of linear orders each of cardinality
$\lambda_i$ of length $\max\pcf(\{\lambda_j:j<i\})$ (i.e., $\ens(\lambda_i,
\max\pcf(\{\lambda_j: j<i\}))$). For some $\alpha$, $\lambda_i$ is
defined if and only if $i<\alpha$. Clearly, $\alpha$ is a limit ordinal
$\le\cf(\mu)$, and $\lambda=:\max\pcf(\{\lambda_j: j<\alpha\})$ is $>\mu$ [as
otherwise $\lambda<\mu$ (as $\lambda$ is regular by \cite[1.x]{Sh:345a}), so
there is $\lambda_\alpha$ as required among $\{(\lambda+\mu_\alpha)^{+\gamma}:
\gamma<(\lambda+\mu_\alpha)^{+4}\}$]. So clearly $\mu<\lambda=\cf(\lambda)<
\pp^+(\mu)$ and by \ref{3.3}(2) there is an entangled linear order of
cardinality $\lambda$ and density $\leq\sum\limits_{i<\alpha}\max\pcf(\{
\lambda_j: j<i\})\leq\mu$. If $\lambda$ is a successor cardinal then we are
done. Otherwise, clearly $\mu^{+\mu^{+4}}\leq \pp^+(\mu)$, and hence we can
apply part (2).

\noindent 2) It follows by the claims below, each has the conclusion of
\ref{5.5}(2) from assumptions which are not necessarily implied by the
assumption of \ref{5.5}(2), but always at least one applies.

\begin{claim}
\label{5.5A}
\begin{enumerate}
\item Assume $\cf(\mu)\leq\kappa\leq\mu^*<\mu\leq 2^\kappa$, $\mu\leq\chi$
and $\chi^{+\kappa^{+4}}\leq\pp_\kappa(\mu)$. Then there is $\lambda^+\in
[\chi,\chi^{+\kappa^{+4}}]$ in which there is an entangled linear order $\cI$
of density $\leq\mu$ but $\ge\mu^*$.
\item If in addition $(\forall\alpha<\mu)(|\alpha|^\kappa<\mu)$ then we can
add ``$\cI$ is $\cf(\mu)$-entangled.''
\item There is $\gamma<\kappa^{+4}$ and a set $\gb\subseteq\Reg\cap\mu
\setminus\mu^*$ of $(\mu^*,\kappa^+,2)$-inaccessible cardinals, $|\gb|\leq
\kappa$, $\gb$ is the disjoint union of $\gb_\vare$ (for $\vare<\kappa$),
$\sup(\gb_\vare)$ is the same for $\vare<\kappa$ and $J^{\bd}_{\gb}\subseteq
J_{<\chi^{+\gamma+1}}[\gb]$ and $\chi^{+(\gamma+1)}\in\pcf(\gb_\vare)$.
\end{enumerate}
\end{claim}

\noindent{\em Proof of the claim:}\hspace{0.2in} 1)\ \ \ Of course we can
decrease $\mu$ as long as $\mu^*<\mu$, $\cf(\mu)\leq\kappa$, $\chi^{+
\kappa^{+4}}\leq\pp_\kappa(\mu)$. By \cite[2.3]{Sh:355}, without loss of
generality we have:
\[\ga\subseteq (\mu^*,\mu)\cap\Reg\ \&\ |\ga|\leq\kappa\quad
\Rightarrow\quad\max\pcf(\ga)<\pp^+(\mu).\]
We choose by induction on $i<\kappa$, $\gb_i$ and $\gamma_i$ such that:
\begin{description}
\item[(i)]   $\gb_i\subseteq{\rm Reg}\cap\mu\setminus\bigcup\limits_{j<i}
\gb_j\setminus\mu^*$,
\item[(ii)]  $\gamma_i<\kappa^{+4}$ is a successor ordinal,
\item[(iii)] $\chi^{+\gamma_i}\in\pcf(\gb_i)$,
\item[(iv)]  $\gamma_i$ is the first successor ordinal for which $\chi^{+
\gamma_i}\notin\pcf(\bigcup\limits_{j<i}\gb_j)$,
\item[(v)]   all members of $\gb_i$ are $(\mu_0,\kappa^+,2)$-inaccessible
\[(\mbox{i.e.}\quad\theta\in\gb_i\ \&\ \ga\subseteq (\mu_0,\theta)\ \&\ |\ga|
\leq\kappa\quad\Rightarrow\quad\max\pcf(\ga)<\theta),\]
\item[(vi)]  $\gb_i$ has cardinality $\leq\kappa$.
\end{description}
Note that this is possible, since if $|\gb|\leq\kappa$ then $\pcf(\gb)$ cannot
contain the interval $[\chi,\chi^{+\kappa^{+4}}]\cap\Reg$ (see \cite[\S
3]{Sh:410}). Let $\gd=:\{\chi^{+\gamma_i}:i<\kappa\}$, let $\langle \gb_\theta
[\gd]:\theta\in\pcf(\gd)\rangle$ be as in \cite[2.6]{Sh:371}. Note that we
know $\pcf(\gd)\subseteq [\chi,\chi^{+\kappa^{+4}}]$ (by \cite[4.2]{Sh:400}).
Let $\theta\in\pcf(\gd)$ be minimal such that $\otp(\gb_\theta [\gd])\ge
\kappa$, so necessarily $\theta$ is a successor cardinal. Let $\langle
\gd_\alpha: \alpha<\kappa\rangle$ be a partition of $\gb_\theta[\gd]$ to
pairwise disjoint subsets of order type $\geq\kappa$. Let $\gb'_\alpha=
\bigcup\left\{\gb_i:\chi^{+\gamma_i}\in \gd_\alpha\right\}$ (for $\alpha<
\kappa$) and $\ga=\bigcup\limits_{\alpha<\kappa}\gb'_\alpha$. Now we can
finish by \ref{3.2}(1).

\noindent 2)\ \ \ In this case we can in the beginning increase $\mu^*$
(still $\mu^*<\mu$) such that the ``wlog'' in the second sentence of the
proof of \ref{5.5A}(1) holds. Necessarily $\sup(\gb_i)=\mu$ for each
$i<\kappa$.

\noindent 3)\ \ \ Included in the proof above. \hfill$\square_{\ref{5.5A}}$

\begin{claim}
\label{5.5C}
If $\cf(\mu)<\kappa=\cf(\kappa)<\mu\leq\chi$, $\chi^{+\kappa^{+4}}\leq
\pp(\mu)$, $\mu$ is $(*,\kappa^+,2)$-inaccessible then
\begin{enumerate}
\item We can find an increasing sequence $\langle\lambda_i:i<\cf(\mu)\rangle$
of regular cardinals with limit $\mu$ such that for each $i$, $\ens(\lambda_i,
\kappa)$ and $\pcf(\{\lambda_i: i<\cf(\mu)\})$ has a member in $(\chi,
\chi^{+\kappa^{+4}})$ and $\lambda_i>\max(\pcf(\{\lambda_j: j<i\}))$ and
$\lambda_i>\kappa$.
\item In addition $\lambda_i\in\pcf^{ex}_\kappa(\ga_i)$ for some sets
$\ga_i\subseteq\Reg\cap\lambda_i\setminus\bigcup\limits_{j<i}\lambda_j$
of $(\mu,\kappa^+,2)$-inaccessible cardinals of cardinality $\kappa$. If
$\cf(\mu)>\aleph_0$ then $\prod\limits_{i<\cf(\mu)}\lambda_i/J^{\bd}_{
\cf(\mu)}$ has true cofinality.
\end{enumerate}
\end{claim}

\noindent{\em Proof of the claim:}\hspace{0.2in} Choose $\mu^*\in (\kappa,
\mu)$ such that
\[\mu'\in (\mu^*,\mu)\ \&\ \cf(\mu')\leq\kappa\quad\Rightarrow\quad\pp_\kappa
(\mu')<\mu\]
(exists, as $\mu$ is $(*,\kappa^+,2)$-inaccessible; see \cite[2.3]{Sh:355}).
Let $\gb_i,\gamma_i$ (for $i<\kappa$) be as in the proof of \ref{5.5A}; so
$\min(\gb_i)>\mu^*$, $\gd=:\{\chi^{+\gamma_i}:i<\kappa\}$, $\langle\gb_\theta
[\gd]:\theta\in\pcf(\gd)\rangle$ be as in \cite[2.6]{Sh:371}. Let $\theta\in
\pcf(\gd)$ be minimal such that $\otp(\gb_\theta[\gd])=\kappa$. Without loss
of generality $\gb_\theta[\gd]=\gd$, so $\theta=\max\pcf(\gd)$. Note that
$\theta\in\pcf(\gd)\subseteq (\chi,\chi^{+\kappa^{+4}})$ is a successor
cardinal. Let $\langle \mu_i: i<\cf(\mu)\rangle$ be strictly increasing
continuous with limit $\mu$ with $\mu_0>\mu^*$. Let $\ga=:\bigcup\limits_{i<
\kappa}\gb_i$ and let $\langle\gb_\sigma[\ga]:\sigma\in\pcf(\ga)\rangle$ be
as in \cite[2.6]{Sh:371}. For each $\vare<\cf(\mu)$, we can find finite
$\goe_\vare\subseteq\pcf^{ex}_\kappa(\ga\cap\mu_\vare)$ and $\gc_\vare
\subseteq\ga\cap\mu_\vare$, $|\gc_\vare|<\kappa$ such that $\ga\cap\mu_\vare
\subseteq\bigcup\limits_{\sigma\in\goe_\vare}\gb_\sigma[\ga]\cup\gc_\vare$.
As $\cf(\mu)<\kappa=\cf(\kappa)$, we can find $\zeta<\kappa$ such that:
\begin{description}
\item[$(*)$] \qquad $(\forall\vare<\cf(\mu))(\gc_\vare\subseteq\bigcup
\limits_{i<\zeta}\gb_i)$.
\end{description}
So by renaming $\zeta=0$ (so $\gc_\vare=\emptyset$). By \ref{5.7}(4) below
each $\sigma\in\goe_\vare$ satisfies $\ens(\sigma,\kappa)$. As in the proof
of \cite[1.5]{Sh:371} and \cite[6.5]{Sh:430} one of the following holds:
\begin{description}
\item[$(*)_1$]  $\cf(\mu)>\aleph_0$ and for some $S\subseteq\cf(\mu)$
unbounded, $\otp(S)=\cf(\mu)$, and $\sigma_\vare\in\goe_\vare$ for $\vare\in 
S$, $\theta=\tcf(\prod\limits_{\vare\in S}\sigma_\vare/J^{\bd}_S)$,
\item[$(*)_2$]  $\cf(\mu)=\aleph_0$ and for some increasing sequence
$\langle\sigma_\zeta: \zeta<\cf(\mu)\rangle$ of regular cardinals from
$\bigcup\limits_{\vare<\cf(\mu)}\goe_\vare$ with limit $\mu$, and an ideal $I$
on $\omega$ extending $J^{\bd}_{cf(\mu)}$, $\theta=\tcf(\prod\limits_{\zeta<
\cf(\mu)}\sigma_\zeta/I)$.
\end{description}
If $(*)_1$ holds then (by the choice of $\mu_0$) without loss of generality
$\vare\in S\ \Rightarrow\ \sigma_\vare>\max\pcf(\{\sigma_\zeta:\zeta<\vare
\})$, so we can apply \ref{3.1}. If $(*)_2$ holds, necessarily $\max\pcf(\{
\sigma_\zeta:\zeta<\vare\})$ is $\sigma_{\vare-1}$ so we can apply
\ref{3.1}. In both cases we get the desired conclusion.
\hfill$\square_{\ref{5.5C}}$

\begin{claim}
\label{5.5B}
Assume that $\mu$ is a singular strong limit cardinal, or at least that
it is $(*,\kappa,2)$-inaccessible for every $\kappa<\mu$, $\cf(\mu)<\mu\leq
\chi$ and $\chi^{+\mu^{+4}}\le\pp(\mu)$. Then
\begin{enumerate}
\item We can find $\lambda^+\in (\chi,\chi^{+\mu^{+4}})$ in which there is an
entangled linear order with density $\mu$.
\item Moreover we can find a strictly increasing sequence $\langle\lambda_i:
i<\cf(\mu)\rangle$ with limit $\mu$, $\max\pcf(\{\lambda_j: j<i\})<\lambda_i$,
$\lambda^+=\max\pcf(\{\lambda_i:i<\cf(\mu)\})$. Letting $\kappa_i=(
\sum\limits_{j<i}\lambda_j+\mu^*)^+$ we can also find a set $\ga_i\subseteq
\Reg\cap\lambda_i\setminus\bigcup\limits_{j<i}\lambda_j$ of
$(\kappa^+_{j_i},\kappa^+_{j_i},2)$-inaccessible cardinals of cardinality
$\kappa_{j_i}$ such that $\lambda_i\in\pcf^{ex}_{\kappa_i}(\ga_i)$ and $j_i
\leq i$. Also $\cf(\mu)>\aleph_0$ implies $\lambda=\tcf(\prod\limits_{i<
\cf(\mu)}\lambda_i/J^{\bd}_{\cf(\mu)})$ and $(\forall i<\cf(\mu))(j_i=i)$.
\end{enumerate}
\end{claim}

\noindent{\em Proof of the claim:}\hspace{0.2in} Choose $\langle\mu_i: i<
\cf(\mu)\rangle$ be strictly increasing continuous with limit $\mu$. By
induction on $\vare<\cf(\mu)$ we choose $\theta^\vare\in (\chi,\chi^{+
\mu^{+4}})$ and $\langle\lambda^\vare_\zeta:\zeta<\cf(\mu)\rangle$ as follows:
arriving to $\vare$ we apply the proof of \ref{5.5C} to $\kappa_\vare=
\mu^+_\vare$ and\footnote{We could have asked $\chi_\vare=\max\{\max\pcf(\{
\theta_\zeta:\zeta<\vare\}),\chi\}$ and thus later ``the $\{\lambda^\vare_\zeta
:\zeta<\cf(\mu)\}$ are pairwise disjoint'' (while omitting ``few''
$\lambda^\vare_\zeta$ for each $\vare$)} $\chi_\vare=\sup(\{\theta_\zeta:
\zeta<\vare\}\cup\{\chi\})$ and get $\langle\lambda^\vare_\zeta: \zeta<\cf(
\mu)\rangle$ as there. So there is a successor $\theta_\vare\in\pcf(\{
\lambda^\vare_\zeta: \zeta<\cf(\mu)\})\cap [\chi_\vare,\chi_\vare^{
\kappa_\vare^{+4}})$ such that $\theta_\vare=\max\pcf(\{\lambda^\vare_\zeta:
\zeta<\cf(\mu)\})$. Hence $\theta_\vare>\bigcup\limits_{\zeta<\vare}
\theta_\zeta$ and $\chi<\theta_\vare< \chi^{+\mu^{+4}}$, and without loss of
generality $\mu^{+4}_i<\mu_{i+1}$. Let $\chi^{+\gamma_{\vare}}$ be
$\theta_\vare$, $\gb_\vare=:\{\lambda^\vare_\zeta:\zeta<\cf(\mu)\}$ (for
$\vare<\cf(\mu)$), $\theta=\max\pcf(\{\theta_\vare:\vare<\cf(\mu)\})$, $\ga=:
\bigcup\{\gb_\vare:\vare<\cf(\mu)\}$ and without loss of generality
\[J_{<\theta}[\{\theta_\vare: \vare<\cf(\mu)\}]\subseteq \{\{\theta_\vare:
\vare\in a\}: a\subseteq\cf(\mu), |a|<\cf(\mu)\}.\]
Note:
\[\begin{array}{ll}
\pcf(\{\theta_\vare:\vare<\cf(\mu)\})\subseteq & \pcf(\{\chi^{+\gamma+1}:
\gamma<\mu^{+4}\}) \subseteq\\
\ &\Reg\cap[\chi,\chi^{+\mu^{+4}}]\cap(\chi,\pp^+(\mu)),
  \end{array}\]
so each member of $\pcf(\{\theta_\vare: \vare<\cf(\mu)\})$ is a successor
cardinal. Let $\langle\gb_\sigma[\ga]:\sigma\in\pcf(\ga)\rangle$ be as in
\cite[2.6]{Sh:371}.

First assume $\cf(\mu)>\aleph_0$. For every limit $\vare<\cf(\mu)$ let
$\goe_\vare$ be a finite subset of $\pcf(\ga\cap\mu_\vare)$ such that
\begin{description}
\item[($\alpha$)] for some $\zeta_\vare<\vare$, $(\ga\cap\mu_\vare)\setminus
\bigcup\limits_{\sigma\in\goe_\vare} b_\sigma[\ga]\subseteq\{\lambda^\xi_i:
i<\cf(\mu)$ and $\lambda^\xi_i<\mu_{\zeta_\vare}\}$,
\item[($\beta$)]  for every $\zeta<\vare$ and for every $\sigma\in\goe_\vare$,
we have
\[\sigma\in\pcf(\{\lambda^\xi_i: i<\cf(\mu)\ \&\ \mu_\zeta<\lambda^\xi_i<
\mu_\vare\}).\]
(Exist by \cite[6.x]{Sh:430}.)
\end{description}
So for every $\sigma\in\goe_\vare$, by \ref{3.3}(2), there is an entangled
linear order of cardinality $\sigma$. Also, by \cite[6.7]{Sh:430} for some
unbounded set $S\subseteq\cf(\mu)$ and $\sigma_\vare\in\goe_\vare$ we have
$\theta=\tcf\prod\limits_{\vare\in S} \sigma_\vare/J^{\bd}_S$. Note that
without loss of generality $\sigma_\vare>\prod\limits_{\zeta<\vare}
\sigma_\zeta$ (when $\mu$ is strong limit!) or at least $\sigma_\vare>
\max\pcf(\{\sigma_\zeta:\zeta<\vare\})$, so by \ref{3.3}(2) we can get the
desired conclusion.

Now assume $\cf(\mu)=\aleph_0$. Use \cite[6.7]{Sh:430} to find finite
$\goe_\vare\subseteq\pcf(\ga\cap\mu_\vare)$ for $\vare<\omega$ such that
\[\vare<\xi<\omega\quad\Rightarrow\quad\max(\goe_\vare\cup\{\aleph_0\})<\min
(\goe_\zeta\cup\{\mu\})\]
(the $\{\aleph_0\}$, $\{\mu\}$ are for empty $\goe_n$'s), and so $\otp(\goe)
=\omega$, $\sup(\goe)=\mu$ where $\goe=\bigcup\limits_{\vare} e_\vare$ and
$\ga\cap\mu_\vare=\cup\{\gb_\sigma[\ga]:\ \sigma\in\bigcup\limits_{n\leq
\vare}\goe_n\}$; hence $\sigma\in\pcf(\goe)$. Define $h:\pcf(\ga)
\longrightarrow\omega$ by $h(\sigma)=\max\{n<\omega:\ens(\sigma,\mu_n)$ or
$n =0\}$. By \ref{3.4}(2) it suffices to prove, for each $n<\omega$, that
$\{\sigma\in\goe:h(\sigma)\leq n\}\in J_{<\sigma}[\goe]$. This can be
easily checked. \hfill$\square_{\ref{5.5B}}$
\medskip

\noindent 3), 4)\ \ \ Left to the reader. \QED$_{\ref{5.5}}$

\begin{remark}
\label{5.5D}
{\em
\begin{enumerate}
\item Under the assumptions of \ref{5.5C} we can get
\begin{description}
\item[($\circledast$)] there are a successor cardinal $\lambda^+\in (\chi,
\chi^{+\mu^{+4}})$ and an increasing sequence $\bar{\lambda}=\langle
\lambda_i: i<\delta\rangle$ with limit $\mu$ such that

$\delta<\kappa^+$,  $\lambda^+=\tcf(\prod\limits_{i<\delta}\lambda_i/J)$,
$J$ an ideal on $\delta$ extending $J^{\bd}_\delta$,

$\lambda_i$ is $(\mu^*,\kappa^+,2)$-inaccessible (where $\mu^*=
\min\{\mu': \mu^*$ is $(\mu',\kappa^+,2)$-inaccessible$\}$) and

there is an entangled linear order of cardinality $\lambda^+$.
\end{description}
\item In the proof of \ref{5.5B} we can have $\ga=\bigcup\limits_{i<\cf(\mu)}
\ga_i$, $\mu=\sup(\ga_i)$ for each $i<\cf(\mu)$, and let $\ga=\{
\lambda^\vare_\zeta:\vare,\zeta<\cf(\mu)\}$,
\[\ga_i=\{\lambda^\vare_\zeta\in \ga:\vare\in\ga_i\mbox{ but for no }\xi<
\vare\mbox{ does }\lambda^\vare_\zeta\in\{\lambda^\xi_j:j<\cf(\mu)\}\}\]
and apply \ref{3.2}(1).
\end{enumerate}
}
\end{remark}

\begin{proposition}
\label{5.7}
Let $\langle\gb_\theta[\ga]:\theta\in\pcf(\ga)\rangle$ be as in
\cite[2.6]{Sh:371}.
\begin{enumerate}
\item If $|\ga|\geq\kappa$ (and $|\ga|<\min(\ga)$) then $\pcf^{ex}_\kappa(
\ga)$ has a last element.
\item Assume that $\theta=\max\pcf(\ga)$, $\gc\subseteq\theta\cap\pcf(\ga)$,
$|\gc|<\min(\gc)$ and
\[\gb\in J_{<\theta}[\ga]\quad \Rightarrow\quad\gc\cap\pcf(\ga\setminus\gb)
\neq\emptyset.\]
Then $\theta\in\pcf(\gc)$.
\item If $\ga=\bigcup\limits_{i<\sigma}\ga_i$, $\sigma<\cf(\kappa)$, $\theta
\in\pcf^{ex}_\kappa(\ga)$ then we can find finite $\goe_i\subseteq
\pcf^{ex}_\kappa\ga_i$ for $i<\sigma$ such that $|\ga_i\setminus
\bigcup\limits_{\lambda\in\goe_i}\gb_\lambda[\ga_i]|<\kappa$, and $\theta
\in\pcf(\bigcup\limits_{i<\sigma}\goe_i)$.

[And if $\cf(\sigma)>\aleph_0$, $\ga_i$ increasing with $i$ and $S\subseteq
\sigma=\sup S$ then $\theta\in\pcf(\bigcup\limits_{i\in S}\goe_i)$.]
\item  Assume $\chi<\min(\ga)$, and each $\mu\in\ga$ is $(\chi,\kappa,
2)$-inaccessible and $\kappa>\aleph_0$. If $\theta\in\pcf^{ex}_\kappa(\ga)$
then $\ens(\theta,2^\kappa)$ holds exemplified by linear orders of density
$>\chi$.
\end{enumerate}
\end{proposition}

\Proof 1)\ \ \ Among the $\gc\subseteq\ga$ of cardinality $<\kappa$ choose
one with $\max\pcf(\ga\setminus\gc)$ minimal. So $\max\pcf(\ga\setminus\gc)=
\max\pcf^{ex}_\kappa(\ga)$.

\noindent 2)\ \ \ By \cite[1.16]{Sh:345a}.

\noindent 3)\ \ \ Easy, as in \cite[\S 1]{Sh:371}.

\noindent 4)\ \ \ Without loss of generality $\theta=\max\pcf(\ga)$, $\ga$
has no last element and 
\[\mu\in\ga\quad\Rightarrow\quad\theta\notin\pcf^{ex}_\kappa(\ga\cap\mu),\]
so $J^{\bd}_{\ga}\subseteq J^{ex,\kappa}_{<\theta}[\ga]$ (see Definition
\ref{5.6}). We are going to prove the statement by induction on $\theta$.
\smallskip

\noindent If $\ga$ ($=\gb_\theta [\ga])$ can be divided to $\kappa$ sets, no
one of which is in $J_{<\theta}[\ga]+J^{\bd}_{\ga}$, this should be clear
(use e.g.~\ref{3.1}(1): there are such $A_i$ by \cite{EK}, see
e.g.~\cite[Appendix]{Sh:g}).
\smallskip

\noindent Also, if there are $\gc\subseteq\theta\cap\pcf^{ex}_{\kappa}(\ga)$,
(such that $|\gc|<\min(\gc)$ and ) an ideal $I$ on $\gc$, satisfying $\theta
=\tcf\prod\gc/I$, then we can find $\gd\subseteq\theta\cap\pcf(\gc)$ such
that $(\forall\sigma\in\gd)(\max\pcf(\gd\cap\sigma)<\sigma)$, $\theta=\max
\pcf(\gd)$ and then use the induction hypothesis and \ref{3.4}(1). So by
part (1) (of \ref{5.7}) the remaining case is:
\begin{description}
\item[($\boxtimes$)] if $\gc\subseteq\theta\cap\pcf^{ex}_\kappa(\ga)$,
$|\gc|<\min(\gc)$ then $\theta\notin\pcf(\gc)$.
\end{description}
Without loss of generality
\begin{description}
\item[($\boxplus$)] $\ga'\subseteq\ga\ \&\ \theta\in\pcf^{ex}_\kappa(\ga')
\quad\Rightarrow\quad\sup(\theta\cap\pcf^{ex}_\kappa(\ga'))=\sup(\theta\cap
\pcf^{ex}_\kappa(\ga))$.
\end{description}
We can try to choose by induction on $i$, $\theta_i\in\pcf^{ex}_\kappa(\ga)$
such that $\theta>\theta_i>\max\pcf(\{\theta_j: j<i\})$. By localization
(see \cite[\S 3]{Sh:371}) we cannot have $\langle\theta_i: i<|\ga|^+\rangle$
(as $\max\pcf(\{\theta_i: i<|\ga|^+\})\in\pcf(\{\theta_i:i<\alpha\})$ for some
$\alpha<|\ga|^+$, hence $\theta_\alpha<\theta_{\alpha+1}<\max\pcf(\{\theta_i:
i<|\ga|^+\}\leq\max\pcf(\{\theta_i:i<\alpha\})=\theta_\alpha$, a
contradiction). So for some $\alpha<|\ga|^+$, $\theta_i$ is defined if
and only if $i<\alpha$. If $\max\pcf(\{\theta_i:i<\alpha\})<\theta$ we can get
a contradiction to $(\boxplus)$ (by \ref{5.7}(1)). If the equality holds, we
get contradiction to $(\boxtimes)$. \QED$_{\ref{5.7}}$
\medskip

We want to state explicitly the pcf theorems behind \ref{5.5}.

\begin{proposition}
\label{5.7A}
Assume $\mu$ is a singular cardinal which is a fix point and $\mu_0<\mu$.
\begin{enumerate}
\item There are $\lambda$, $\langle(\lambda_i,\bar{\lambda}^i):i<\delta
\rangle$ such that:
\begin{description}
\item[(a)] $\langle\lambda_i:i<\delta\rangle$ is a strictly increasing
sequence of regular cardinals, $\delta$ is a limit ordinal $\leq\cf(\mu)$,
$\lambda_i\in\mbox{\rm Reg}\cap\mu\setminus\mu_0$, $\lambda_i^*=\max\pcf(
\{\lambda_j:j<i\})<\lambda_i$, $\lambda\in\mbox{\rm Reg}\cap (\mu,\pp^+(
\mu))$, $\lambda=\tcf(\prod\limits_{i<\delta}\lambda_i/J^{\bd}_\delta)$,
\item[(b)] $\bar{\lambda}^i=\langle\lambda^i_j: j<(\lambda^*_i)^+\rangle$ is
strictly increasing, $\lambda^i_j\in\mbox{\rm Reg}\cap\lambda_i\setminus
\lambda^*_i$, $\lambda_i=\tcf(\prod\limits_j \lambda^i_j/J^{\bd}_{(
\lambda^*_j)^+})$, $\lambda^i_j>\max\pcf(\{\lambda^i_\zeta:\zeta<j\})$
\end{description}
\item Assume in addition $\mu\leq\chi_0$ and $\chi_0^{+\mu^{+4}}\leq
\pp^+(\mu)$. Then for some $\gamma<\mu^{+4}$, letting $\lambda=\chi^{+\gamma
+1}$, we can find a strictly increasing sequence $\langle\lambda_i:i<\delta
\rangle$ of regular cardinals of the length $\delta=\cf(\delta)<\cf(\mu)$,
$\lambda_i\in\mbox{\rm Reg}\cap\mu\setminus\mu^*$, $\lambda_i>\max\pcf(\{
\lambda_j:j<i\})$, $\lambda=\max\pcf(\{\lambda_i:i<\delta\})$ and letting
$\kappa_i=(\sum\limits_{j<i}\lambda_j+\mu^*)^+$ we can also find sets
$\ga_i\subseteq\mbox{\rm Reg}\cap\lambda_i\setminus\bigcup_{j<i}\lambda_j
\setminus\mu^*$ of $(\kappa_{j_i},\kappa_{j_i}^+,2)$-inaccessible cardinals
of cardinality $\kappa_{j_i}$ such that $\lambda_i\in\pcf^{ex}_{\kappa_i}
(\ga_i)$, $j_i\leq i$. Also if $\cf(\mu)>\aleph_0$ then $\lambda=\tcf(
\prod\limits_{i<\cf(\mu)}\lambda_i/J^{\bd}_{\cf(\mu)})$.
\item If part 2) does not apply then in 1) $\lambda$ is a successor.
\end{enumerate}
\end{proposition}

\Proof Let $\langle\mu_i:1\leq i<\cf(\mu)\rangle$ be increasing continuous
with limit $\mu$, $\mu_1>\mu_0$ and wlog if $\mu$ is $(*,\cf(\mu),
2)$-inaccessible then
\[\mu'\in [\mu_0,\mu)\ \&\ \cf(\mu')\leq\cf(\mu)\quad\Rightarrow\quad
\pp(\mu')<\mu,\]
and hence
\[\ga\subseteq [\mu_0,\mu)\cap\mbox{\rm Reg}\ \&\ |\ga|\leq\cf(\mu)\ \&
\ \sup(\ga)<\mu\quad \Rightarrow\quad \max\pcf(\ga)<\mu.\]

\noindent 1)\ \ \ Try to choose by induction on $i<\cf(\mu)$, $\lambda_i$ and
$\ga_i$, $\lambda_i$, $\bar{\lambda}^i$, $\lambda_i^+$, $j_i$ such that
\begin{description}
\item[$(\alpha)$] $\{\lambda^i_j: j\}\cup\{\lambda_i,\lambda_i^*\}\subseteq
(\mu_{j_i},\lambda)\cap\mbox{\rm Reg}$,
\item[$(\beta)$]  $\mu>\mu_{j_i}$,
\item[$(\gamma)$] $\lambda_i^*$, $\bar{\lambda}^i$, $\lambda_i$ are as
required in 1).
\end{description}
So for some $\alpha$, $(j_i,\lambda^*_i,\bar{\lambda}^i,\lambda_i)$ is defined
if and only if $i<\alpha$, in fact $\alpha$ is limit and $\pcf(\{\lambda_i:i
<\alpha\})\not\subseteq\mu$ (as in the proof of \ref{5.5A}). For some
unbounded set $A\subseteq\delta$ we have $\max\pcf(\{\lambda_i:i\in A\})>\mu$
but it is minimal under this restriction. Now restricting ourselves to $A$
(and renaming) we finish.
\medskip

\noindent 2)\ \ \ It follows from \ref{5.5A}, \ref{5.5B}.

\noindent 3)\ \ \ Should be clear. \QED$_{\ref{5.7A}}$
\medskip

Let us finish this section with stating some results which will be
developed and presented with all details in a continuation of the
present paper.

\begin{proposition}
\label{5.8}
Assume:
\begin{description}
\item[(A)]
\begin{description}
\item[(a)] $T\subseteq\bigcup\limits_{j\leq\delta}\prod\limits_{i<j}
\lambda_i$,
\item[(b)] $T$ closed under initial segments,
\item[(c)] $T_j=:T\cap\prod\limits_{i<j}\lambda_i\neq\emptyset$ for
$j\le\delta$,
\item[(d)] for $j<\delta$, $\eta\in T_j$ we have $(\exists^{\lambda_j}
\alpha<\lambda_j)(\eta {}^\frown\!\langle\alpha\rangle\in T)$,
\item[(e)] $|T_\delta|=\kappa\geq\mu=\cf(\mu)>\sum\limits_{i<\delta}|T_i|
>\sigma$,
\item[(f)] for $\eta\in T_i$, $I_\eta$ is a $\sigma$-complete ideal on
$\lambda_i$,
\item[(g)] $\sigma$ is a regular cardinal, $\cf(\delta)\geq\sigma$,
\end{description}
\item[(B)]
\begin{description}
\item[(a)] $J$ is an ideal on $\delta$ extending $J^{\bd}_\delta$,
\item[(b)] $g_i$ is a function from $T_i $ to $\kappa_i$,
\item[(c)] for $i<\delta$, $\alpha<\kappa_i$, $I^i_\alpha$ is a
$\sigma$-complete ideal on $\lambda_i$,
\item[(d)] if $\langle\eta_{\beta,\vare}:\vare<\vare^*,\beta<\mu\rangle$ are
pairwise distinct members of $T_\delta$, $\vare^*<\sigma$, $i(*)<\delta$, for
each $\vare<\vare^*$, $(\forall\beta<\mu)(\eta_\vare=\eta_{\beta,\vare}\rest
i(*))$ and $\langle\eta_\vare: \vare<\vare^*\rangle$ are pairwise distinct

then for some $A\in J$, for every $i\in\delta\setminus A$ there are
$\nu_\vare\in T_i$ for $\vare<\vare^*$ such that:
\begin{description}
\item[($\alpha$)] $\langle g_i(\nu_\vare): \vare<\vare^*\rangle$ are pairwise
distinct,
\item[($\beta$)]  for every $B_{\nu_{\vare}}\in I^i_{g_i(\nu_\vare)}$
(for $\vare<\vare^*$), for some $\beta<\mu$ we have
\[(\forall\vare<\vare^*)(\nu_\vare\vartriangleleft\eta_{\beta,\vare})\qquad
\quad\mbox{ and }\qquad\quad\eta_{\beta,\vare}(i)\notin B_{\nu_\vare},\]
\end{description}
\end{description}
\item[(C)]  for $i<\delta$, $\langle\cI^i_\zeta: \zeta<\kappa_i\rangle$ is a
sequence of linear orders with universe $\lambda_i$ such that:
\begin{quotation}
\noindent if $\vare(*)<\sigma$, $\langle\zeta_\vare: \vare<\vare(*)\rangle$ is
a sequence of distinct members of $\kappa_1$, $\langle\alpha_{\beta,\vare}:
\beta<\beta(*),\vare<\vare(*)\rangle$ is a sequence of ordinals $<\lambda$
such that
\[(\forall\bar{B}\in\prod\limits_{\vare<\vare(*)}\cI_{\zeta_\vare})(\exists
\beta<\beta^*)(\forall\vare<\vare(*))(\alpha_{\beta,\vare}\notin B_\vare)\]
and $\vare(*)=u\cup v$, $u\cap v=\emptyset$

\noindent then for some $\beta_1,\beta_2<\beta(*)$ we have
\[\vare\in u\quad\Rightarrow\quad\cI^i_{\zeta_\vare}\models\mbox{``}\alpha_{
\beta_1, \vare}<\alpha_{\beta_2,\vare}\mbox{''},\]
\[\vare\in v\quad\Rightarrow\quad\cI^i_{\zeta_\vare}\models\mbox{``}
\alpha_{\beta_1,\vare}>\alpha_{\beta_2,\vare}\mbox{''}.\]
\end{quotation}
\end{description}
Then there is a $(\mu,\sigma)$-entangled linear order of cardinality $\lambda$.
\end{proposition}

\begin{remark}
{\em
\begin{enumerate}
\item  The proof is derived from the proof of \ref{3.3} (and so from
\cite[4.10]{Sh:355}).
\item Are the assumptions reasonable? At least they are not so rare, see
\cite[\S 5]{Sh:430}.
\end{enumerate}
}
\end{remark}

\Proof The desired linear order $\cI$ has the universe $T_\delta$ (which has
cardinality $\lambda)$ with the order:

$\eta <_\cI \nu$\qquad if and only if\qquad for the minimal $i<\delta$ for
which $\eta(i)\neq\nu(i)$ we have
\[\cI_{g(\eta \rest i)} \models\mbox{ `` }\eta(i) < \nu(i)\mbox{ ''}.\]

Details, as said before, will be presented somewhere else, but they
should be clear already. \QED$_{\ref{5.8}}$

\begin{proposition}
\label{5.8A}
We can replace (B)(d) of \ref{5.8} by
\begin{description}
\item[(d1)] if $\langle\eta_{\beta,\vare}:\vare<\vare^*,\beta<\mu\rangle$
are pairwise distinct members of $T_\delta$, $\vare^*<\sigma$,
$i(*)<\delta$ and
\[(\forall\vare<\vare(*))(\forall\beta<\mu)(\eta_\vare=\eta_{\beta,\vare}
\rest i(*))\]
and $\langle\eta_\vare: \vare<\vare^*\rangle$ are pairwise distinct

then for some $\beta<\mu$
\[\hspace{-1cm}\begin{array}{rl}
\big\{i<\delta:\mbox{ for every }\bar{B}\in\prod\limits_{\vare<\vare^*}
I^i_{g_i(\eta_{\beta,\vare}\rest i)}\mbox{ for some }\gamma<\mu\mbox{ we
have }&\ \\
(\forall\vare<\vare^*)(\eta_{\gamma,\vare}\rest i=\eta_{\beta,\vare}\rest i)
\quad\mbox{ and }\quad(\forall\vare<\vare^*)(\eta_{\gamma,\vare}(i)\in
B_\vare)\big\}& = \delta\ \mod\;J,
\end{array}\]
\item[(d2)] if $\eta_\vare\in T_\delta$ ($\vare<\vare^*<\sigma$) are distinct,
then
\[\{i<\delta:\langle g_i(\eta_\vare\rest i):\vare<\vare^*\rangle\mbox{ is with
no repetition }\}\neq\emptyset\ \mod\; J.\]
\end{description}
\end{proposition}

\begin{proposition}
\label{5.8B}
Assume
\begin{description}
\item[(A)] $\lambda=\tcf(\prod\limits_{i<\delta}\lambda_i /J)$, $\lambda_i>
\theta_i=:\max\pcf(\{\lambda_j: j<i\})$,
\item[(B)] $J$ a $\sigma$-complete ideal on the limit ordinal $\delta$,
$\lambda_i=\cf(\lambda_i)>\delta$,
\item[(C)] $g_{i,j}:\theta_i\longrightarrow\kappa_{i,j}$ for $j<j_i$ are such
that for any $w\in [\theta_i]^{<\sigma}$ for some $j<j_i$ the restriction
$g_{i,j}\rest w$ is one to one,
\item[(D)] $\bar{a}=\langle a_i:i<\delta\rangle$, $a_i\subseteq
(\bigcup\limits_{\zeta<i}\{\zeta\}\times j_\zeta)$ are such that for any
$\zeta$ and $j<j_\zeta$ we have
\[\{i:(\zeta,j_\zeta)\in a_i\}\neq\emptyset\ \mod\; J,\]
\item[(E)] $\ens_\sigma(\lambda_i,\prod\limits_{(\zeta,j)\in a_i}
\kappa_{\zeta,j})$.
\end{description}
Then for some $T$ the assumptions of \ref{5.8} hold for $\mu_i=\lambda_i$,
$\kappa_i=:\prod\limits_{(\zeta,j)\in a_i}\kappa_{\zeta,j}$.
\end{proposition}

\begin{proposition}
\label{5.8C}
In \ref{5.8B}, (C) + (D) + (E) holds if
\begin{description}
\item[(C')] $\ens_\sigma(\lambda_i,\kappa)$,
\item[(D')]
\begin{description}
\item[$(\alpha)$] $\kappa^{|\delta|}\geq\sum\limits_{i<\delta}\lambda_i$ but
for $i<\delta$, $\kappa^{|i|}=\kappa$ (so $\delta$ is a regular cardinal) or
\item[$(\beta)$]  $\kappa^{|\delta|}\geq\sum\limits_{i<\delta}\lambda_i$ and
there is a regular ultrafilter $E$ on $\delta$ disjoint from
\[J\cup\{A\subseteq\delta: \otp(A)<\delta\}.\]
\end{description}
\end{description}
\end{proposition}

\begin{proposition}
\label{5.9}
Suppose (A), (B) as in \ref{5.8} and
\begin{description}
\item[(C)] there is a sequence $\langle\cI^i_\zeta: \zeta<\kappa_i\rangle$ of
partial orders with universe $\lambda_i$ such that
\begin{quotation}
\noindent if $\vare(*)<\sigma$, $\langle\zeta_\vare:\vare<\vare(*)\rangle$ a
sequence of ordinals $<\kappa_i$ with no repetitions, $\langle\alpha_{
\beta,\vare}:\beta<\beta(*),\vare<\vare(*)\rangle$ a sequence of ordinals
$<\lambda$ such that
\[(\forall \bar{B}\in\prod\limits_{\vare<\vare(*)}\cI_{\zeta_\vare})(\exists
\beta<\beta^*)(\forall\vare<\vare(*))(\alpha_{\beta,\vare}\notin B_\vare)\]
and $\vare(*)=u\cup v$, $u\cap v=\emptyset$,

\noindent then for some $\beta_1,\beta_2<\beta(*)$ we have:
\[\vare\in u\quad\Rightarrow\quad\cI^i_{\zeta_\vare}\models\mbox{``}
\alpha_{\beta,\vare}<\alpha_{\beta_2,\vare}\mbox{''},\]
\[\vare\in v\quad\Rightarrow\quad\cI^i_{\zeta_\vare}\models\mbox{``}\neg
\alpha_{\beta,\vare}<\alpha_{\beta_2,\vare}\mbox{''}.\]
\end{quotation}
\end{description}
Then there is a Boolean Algebra $B$, $|B|=\lambda$ with neither chain nor
pie of cardinality $\lambda$; moreover for $\vare<\sigma$, $B^\sigma$ has
those properties.
\end{proposition}

\Proof Combine \ref{5.8} and proof of \ref{4.2}. \QED$_{\ref{5.9}}$

\begin{remark}
{\em The parallels of \ref{5.8A}, \ref{5.8B} and \ref{5.8C} hold too.}
\end{remark}

\section{Variants of entangledness in ultraproducts}
In this section we develop results of section 1. The following improves 
\ref{1.8}:

\begin{proposition}
\label{6.1}
Assume that:
\begin{description}
\item[(a)] $D$ is an ultrafilter on $\kappa$, $E$ is an ultrafilter on
$\theta$,
\item[(b)] $g_\vare: \kappa\longrightarrow\theta$ for $\vare<\vare(*)$ are
such that:
\begin{quotation}
\noindent if $\vare_1<\vare_2<\vare(*)$ then $g_{\vare_1}\neq g_{\vare_2}\
\mod\; D$ and

\noindent if $\vare<\vare(*)$, $A\in E$ then $g_\vare^{-1}[A]\in D$,
\end{quotation}
\item[(c)] $\cI$ is a linear order of the cardinality $\lambda\geq\theta$.
\end{description}
Then there exists a sequence $\langle f^\alpha_{\vare}/D:\vare<\vare(*),
\alpha<\lambda\rangle$ of pairwise distinct members of $\cI^\kappa/D$ such
that for each $\alpha<\beta<\lambda$:
\[\mbox{ either }\ (\forall\vare<\vare(*))(f^\alpha_\vare/D<f^\beta_\vare/D)
\qquad\mbox{ or }\ (\forall\vare<\vare(*))(f^\beta_\vare/D<f^\alpha_\vare/D
).\]
In particular, the linear order $\cI^\kappa/D$ is not entangled (here
$\vare(*)=2$ suffices) and the Boolean algebra $(\BA_{\inter}(\cI))^\kappa/D$
is not $\lambda$-narrow.
\end{proposition}

\Proof Choose pairwise distinct $a^\alpha_\zeta\in\cI$ for $\alpha<\lambda$,
$\zeta<\theta$. Let $f^\alpha_\vare:\kappa\longrightarrow\cI$ be given by
$f^\alpha_\vare(i)=a^\alpha_{g_\vare(i)}$. Note that if $\alpha_1\neq
\alpha_2$ then $\{i<\kappa:f^{\alpha_1}_{\vare_1}(i)=f^{\alpha_2}_{\vare_2}
(i)\}=\emptyset$. If $\alpha_1=\alpha_2=\alpha$ but $\vare_1\neq\vare_2$ then
\[\{i\in\kappa:f^\alpha_{\vare_1}(i)=f^\alpha_{\vare_2}(i)\}=\{i\in\kappa:
g_{\vare_1}(i)=g_{\vare_2}(i)\}\notin D\]
(as $g_{\vare_1}\neq g_{\vare_2}\ \mod\; D$). Consequently, if $(\alpha_1,
\vare_1)\neq (\alpha_2,\vare_2)$ then $f^{\alpha_1}_{\vare_1}/D\neq
f^{\alpha_2}_{\vare_2}/D$. Suppose now that $\alpha<\beta<\lambda$. Let
$A=\{\zeta<\theta:a^\alpha_\zeta<a^\beta_\zeta\}$. Assume that $A\in E$ and
let $\vare<\vare(*)$. Then
\[\{i<\kappa: f^\alpha_\vare(i)<f^\beta_\vare(i)\}=\{i<\kappa: a^\alpha_{
g_\vare(i)}<a^\beta_{g_\vare(i)}\}=g^{-1}_\vare[A]\in D\]
and hence $f^\alpha_\vare/D<f^\beta_\vare/D$. Similarly, if $A\notin E$ then
for each $\vare<\vare(*)$, $f^\alpha_\vare/D>f^\beta_\vare/D$.

Now clearly $\cI^\kappa/D$ is not entangled, but what about the narrowness of
$(\BA_{\inter}(\cI))^\kappa/D$? Remember that $\BA_{\inter}(\cI^\kappa/D)$
embeds into $\BA_{\inter}(\cI)^\kappa/D$. So if the cardinality of
$\BA_{\inter}(\cI)$ is regular we can just quote \ref{1.4} (a)\
$\Leftrightarrow$\ (c) (here $\sigma=\aleph_0$, see Definition \ref{1.3}).
Otherwise, just note that for a linear order $\cJ$, if $a_\ell<b_\ell$ (for
$\ell=0,1$) and $[a_\ell,b_\ell)\in \BA_{\inter}(\cJ)$ are comparable and
$\{a_0,b_0,a_1,b_1\}$ is with no repetition then
\[a_0<_{\cJ} a_1\quad \Leftrightarrow\quad\neg(b_0<_{\cJ} b_1);\]
this can be applied by the statement above. \QED$_{\ref{6.1}}$

\begin{remark}
\label{6.2}
{\em  Proposition \ref{6.1} shows that entangledness can be destroyed by
ultraproducts. Of course, to make this complete we have to say how one can
get $D, E, g_\vare$'s satisfying (a)--(b) of \ref{6.1}. But this is easy:
\begin{enumerate}
\item For example, suppose that $E$ is a uniform ultrafilter on $\theta$,
$D=E\times E$ is the product ultrafilter on $\theta\times\theta=\kappa$,
$\vare(*)=2$ and $g_\vare:\theta\times\theta\longrightarrow\theta$ (for
$\vare<2$) are given by $g_0(i,j)=i$, $g_1(i,j)=j$. Then $E$, $D$, $\vare(*)$,
$g_\vare$ satisfy the requirements (a)--(b) of \ref{6.1}.
\item More general, assume that $\theta\leq\kappa$, $E$ is a non-principal
ultrafilter on $\theta$, $\vare(*)\leq 2^\kappa$. Let $g_\vare:\kappa
\longrightarrow\theta$ for $\vare<\vare(*)$ constitute an independent family
of functions. Then the family $\{g^{-1}_\vare[A]:\vare<\vare(*), A\in E\}$
has the finite intersections property so we can complete it to an ultrafilter
$D$. One can easily check that $D,E,g_\vare$'s satisfy (a)--(b) of
\ref{6.1}.
\end{enumerate}
}
\end{remark}

\begin{remark}
\label{6.3}
{\em
\begin{enumerate}
\item In \ref{6.1} we did not use ``$<$ is a linear order''. Thus for any
binary relation $R$ the parallel (with $R$, $\neg R$ instead of $<$, $>$)
holds. In particular we can apply this to Boolean algebras.
\item We can weaken assumption (b) in \ref{6.1} and accordingly the
conclusion. For example we can replace (b) by
\begin{description}
\item[(b)$^*$] $\cP\subseteq \cP(\vare(*))$ and for each $A\in E$
\[\{\vare<\vare(*): g^{-1}_\vare[A]\in D\}\in\cP,\]
\end{description}
and the conclusion by
\[\{\vare<\vare(*): f^\alpha_\vare/D<f^\beta_\vare/D\}\in\cP.\]
\item A kind of entangledness can be preserved by ultraproducts, see
\ref{6.5} below. More entangledness is preserved if we put additional demands
on the ultrafilter, see \ref{6.8}.
\item Let us explain why we introduced ``positive entangledness'' in \ref{6.4}.
The proof of \ref{6.1} excludes not only ``full'' entangledness but many
variants (for the ultrapower). Now the positive $\sigma$-entangledness seems
to be the maximal one not excluded. Rightly so by \ref{6.5}. 
\item So if we can find an linear order $\cI$ which is $\mu^+$-entangled for
some $\mu$ such that $\mu^{\aleph_0}=\mu$ (or at least $\mu^{\aleph_0}<|\cI|$)
then we can answer Monk's problem from the introduction: if $D$ is a non
principal non separative ultrafilter on $\omega$ (see Definition \ref{6.6};
they exist by \ref{6.2}(1)), then $\cI^\omega/D$ is not $\mu^+$-entangled (by
\ref{6.1}). Thus if $B$ is the interval Boolean algebra of $\cI$ then
\[\inc(B)\leq\mu,\quad \mu^{\aleph_0}<|\cI|,\quad\mbox{ but }\ \ \inc(
B^\omega/D)\geq |\cI|^{\aleph_0}\]
(in fact $\inc^+(B^\omega/D)=|\cI^\omega/D|^+=(|\cI|^\omega/D)^+$, $\inc^+(B)
\leq\mu^+$). In fact for any infinite Boolean algebra $B$ and a non principal
ultrafilter $D$ on $\omega$ we have $\inc(B^\omega/D)\geq (\inc(B))^\omega/D$
(as for $\lambda_n$ inaccessible the linear order $\prod\limits_{n<\omega}(
\lambda_n,<)/D$ cannot be $\mu^+$-like (see \cite{MgSh:433}).
\end{enumerate}
}
\end{remark}

\begin{proposition}
\label{6.5}
Suppose that $\kappa<\sigma<\lambda$ are regular cardinals such that
$(\forall\theta<\lambda)(\theta^{<\sigma}<\lambda)$. Assume that $D$ is an
ultrafilter on $\kappa$. Then:
\begin{quotation}
\noindent if $\cI$ is a positively $\sigma$-entangled linear order of the size
$\lambda$

\noindent then $\cI^\kappa/D$ is positively $\sigma$-entangled.
\end{quotation}
\end{proposition}

\Proof Suppose $f^\alpha_\vare/D\in \cI^\kappa/D$ (for $\vare<\vare(*)<\sigma$,
$\alpha<\lambda$) are such that
\[(\forall \alpha<\beta<\lambda)(\forall\vare<\vare(*))(f^\alpha_{\vare}/D
\neq f^\beta_{\vare}/D).\]
Let $u\in\{\emptyset,\vare(*)\}$. For $\alpha<\lambda$ let $A_\alpha=\{
f^\alpha_\vare(i):\vare<\vare(*),i<\kappa\}$, $A_{<\alpha}=\bigcup\limits_{
\beta<\alpha} A_\beta$, $B_\alpha=A_\alpha\cap A_{<\alpha}$. Note that
$|A_\alpha|<\sigma$, $|A_{<\alpha}|<\lambda$ (for $\alpha<\lambda$). If
$\delta<\lambda$, $\cf(\delta)=\sigma$ then $|B_\delta|<\cf(\delta)$ and
consequently for some $h(\delta)<\delta$ we have $B_\delta\subseteq
A_{<h(\delta)}$. By Fodor's lemma we find $\beta_0<\lambda$ and a stationary
set $S\subseteq\lambda$ such that
\[\delta \in S\quad \Rightarrow\quad \cf(\delta)=\sigma\ \& \ B_\delta
\subseteq A_{<\beta_0}.\]
For $\delta\in S$ let $Y_\delta=\{(\vare,i,f^\delta_\vare(i)): \vare<\vare(*),
i<\kappa,f^\delta_\vare(i)\in A_{<\beta_0}\}$. As there are $|A_{<\beta_0
}|^{<\sigma}<\lambda$ possibilities for $Y_\delta$, we can find a stationary
set $S_1\subseteq S$ and $Y$ such that for $\delta\in S_1$ we have
$Y_\delta=Y$. Let $\alpha,\beta\in S_1$, $\alpha<\beta$ and $\vare<\vare(*)$.
If there is $x$ such that $(\vare,i,x)\in Y$ then $f^\alpha_\vare(i)=
f^\beta_\vare(i)$. Thus
\[E_\vare=\{i<\kappa: (\exists x)((\vare,i,x)\in Y)\}=\emptyset\ \mod\; D,\]
as $\{i<\kappa: f^\alpha_\vare(i)=f^\beta_\vare(i)\}=\emptyset\ \mod\; D$ (for
distinct $\alpha,\beta\in S_1$). Clearly if $\alpha,\beta\in S_1$, $\alpha<
\beta$, $\vare<\vare(*)$ and $i\in\kappa\setminus E_\vare$ then
$f^\alpha_\vare(i)\neq f^\beta_\vare(i)$. Hence we may apply the positive
$\sigma$-entangledness of $\cI$ to
\[\{f^\alpha_\vare(i):\alpha\in S_1, \vare< \vare(*), i\in\kappa\setminus
E_\vare\}\ \ \mbox{ and }\ \ u'=\{(\vare,i): \vare\in u, i\in\kappa\setminus
E_\vare\}.\]
Consequently we have $\alpha<\beta$, both in $S_1$ and such that
\[(\forall\vare<\vare(*))(\forall i\in\kappa\setminus E_\vare)(
f^\alpha_\vare(i)<f^\beta_\vare(i)\ \Leftrightarrow\ \vare\in u).\]
Since $\kappa\setminus E_\vare\in D$ we get $(\forall\vare<\vare(*))(
f^\alpha_\vare/D<f^\beta_\vare/D\ \Leftrightarrow\ \vare\in u)$.
\QED$_{\ref{6.5}}$

\begin{definition}
\label{6.6}
An ultrafilter $D$ on $\kappa$ is called {\em separative} if for every
$\bar{\alpha},\bar{\beta}\in{}^\kappa\kappa$ such that $(\forall i<\kappa)
(\alpha_i\neq\beta_i)$ there is $A\in D$ such that
\[\{\alpha_i:i\in A\}\cap\{\beta_i:i\in A\}=\emptyset.\]
\end{definition}

\begin{remark}
\label{6.6A}
{\em  So \ref{6.2}(1) says that $D\times D$ is not separative. }
\end{remark}

\begin{proposition}
\label{6.7}
Suppose that $D$ is a separative ultrafilter on $\kappa$, $n<\omega$ and
$\bar{\alpha}^\ell\in {}^\kappa\kappa$ (for $\ell<n$) are such that
\[(\forall\ell_0<\ell_1<n)(\bar{\alpha}^{\ell_0}/D\neq \bar{\alpha}^{\ell_1}
/D).\]
Then there is $A\in D$ such that the sets $\{\alpha^\ell_i:i\in A\}$ (for
$\ell<n$) are pairwise disjoint.
\end{proposition}

\Proof For $\ell<m<n$, by \ref{6.6}, there is $A_{\ell,m}\in D$ such that
\[\{\alpha^\ell_i: i\in A_{\ell,m}\}\cap\{\alpha^m_i:i\in A_{\ell,m}\}=
\emptyset.\]
Now $A=\bigcap\limits_{\ell<m<n} A_{\ell,m}$ is as required.\QED$_{\ref{6.7}}$

\begin{proposition}
\label{6.8}
Assume $\mu=\cf(\mu)$, $(\forall \alpha<\mu)(|\alpha|^\kappa<\mu)$. Suppose
that $\cI$ is a $(\mu,\kappa^+)$-entangled linear order and $D$ is a
separative ultrafilter on $\kappa$. Then the linear order $\cI^\kappa/D$ is
$(\mu,\aleph_0)$-entangled.
\end{proposition}

\Proof Let $n<\omega$, $u\subseteq n$ and $f^\ell_\alpha/D\in\cI^\kappa/D$ for
$\ell<n$, $\alpha<\mu$ be pairwise distinct. By \ref{6.7} we may assume that
(for each $\alpha<\mu$) the sets $\langle\{f^\ell_\alpha(i): i<\kappa\}:\ell
<n\rangle$ are pairwise disjoint. Applying $\Delta$-lemma we may assume
that $\{\langle f^\ell_\alpha(i): i<\kappa,\ell<n\rangle: \alpha<\mu\}$ forms
a $\Delta$-system of sequences and that the diagram of the equalities does not
depend on $\alpha$.  Let
\[A=\{(i,\ell)\in\kappa\times n:(\forall\alpha<\beta<\mu)(f^\ell_\alpha(i)=
f^\ell_\beta(i))\}\]
(i.e. the heart of the $\Delta$-system). Note that for each $\ell<n$ the set
$\{i\in\kappa:(i,\ell)\in A\}$ is not in $D$ (as $f^\ell_\alpha/D$'s are
pairwise distinct). Consequently we may modify the functions $f^\ell_\alpha$
and we may assume that $A=\emptyset$. Now we have
\[f^{\ell_0}_{\alpha_0}(i_0)=f^{\ell_1}_{\alpha_1}(i_1)\quad\Rightarrow\quad
\alpha_0=\alpha_1\ \&\ \ell_0=\ell_1.\]
It is easy now to apply the $(\mu,\kappa^+)$-entangledness of $\cI$ and find
$\alpha<\beta<\mu$ such that $f^\ell_\alpha/D<f^\ell_\beta/D\ \equiv\ \ell\in
u$. \QED$_{\ref{6.8}}$

\begin{proposition}
\label{6.9}
\begin{enumerate}
\item If $D$ is a selective ultrafilter on $\kappa$ (i.e. for every $f:\kappa
\longrightarrow\kappa$ there is $A\in D$ such that either $f\rest A$ is
constant or $f\rest A$ is one-to-one) then $D$ is separative.
\item If no uniform ultrafilter on $\omega$ is generated by less than
continuum sets (i.e. ${\frak u}=2^{\aleph_0}$) then there exists a separative
ultrafilter on $\omega$.
\end{enumerate}
\end{proposition}

\Proof 1)\ \ \ Suppose that $\bar{\alpha},\bar{\beta}\in{}^\kappa\kappa$ are
such that $(\forall i<\kappa)(\alpha_i\neq\beta_i)$. We find $A\in D$ such
that $\bar{\alpha}\rest A$, $\bar{\beta}\rest A$ are either constant or
one-to-one. If at least one of them is constant then, possibly omitting one
element from $A$, the sets $\{\alpha_i: i\in A\}$, $\{\beta_i: i\in A\}$ are
disjoint, so assume that both sequences are one-to-one. Choose inductively
$m_i\in\{0,1,2\}$ (for $i\in A$) such that 
\[i,j\in A\ \&\ \alpha_i=\beta_j\quad\Rightarrow\quad m_i\neq m_j.\]
There are $m^*$ and $B\subseteq A$, $B\in D$ such that $(\forall i\in
B)(m_i=m^*)$. Then 
\[\{\alpha_i: i\in B\}\cap\{\beta_i: i\in B\}=\emptyset.\]

\noindent 2)\ \ \ Straightforward. \QED$_{\ref{6.9}}$

\begin{proposition}
\label{6.10}
If $D$ is a uniform not separative ultrafilter on $\kappa$, $\cI$ is a linear
ordering of size $\lambda\geq\kappa$ then the linear order $\cI^\kappa/D$ is
not $(\lambda,2)$-entangled.
\end{proposition}

\Proof As $D$ is not separative we have $\bar{\alpha},\bar{\beta}\in
{}^\kappa\kappa$ witnessing it. This means that $\bar{\alpha}\neq\bar{\beta}\
\mod\; D$ and the family
\[\{\{\alpha_i:i\in A\}, \{\beta_i: i\in A\}: A\in D\}\]
has the finite intersection property. Consequently we may apply
\ref{6.1}. \QED$_{\ref{6.10}}$

\begin{conclusion}
\label{6.11}
\begin{enumerate}
\item If $\theta^\kappa/D>2^{2^\theta}$ then $D$ is not separative.
\item If $D$ is regular ultrafilter on $\kappa$, $2^\kappa>\beth_2$ then $D$
is not separative.
\end{enumerate}
\end{conclusion}

\begin{remark}
\label{6.11A}
{\em
If there is no inner model with measurable then every $D$ is regular or close
enough to this to give the result.
}
\end{remark}

\Proof 1)\ \ \ For every $f\in {}^\kappa\theta$ the family $E_f=\{A\subseteq
\theta: f^{-1}[A]\in D\}$ is an ultrafilter on $\theta$. For some $g_0, g_1
\in {}^\kappa\theta$ we have $E_{g_0}=E_{g_1}$ but $g_0/D\neq g_1/D$ and hence
 we are done.

\noindent 2)\ \ \ The regularity implies $\aleph_0^\kappa/D=2^\kappa$
(see \cite{CK}), so by the first part we are done. \QED$_{\ref{6.11}}$

\begin{definition}
\label{6.12}
Let $\kappa$, $\sigma$ be cardinal numbers.
\begin{enumerate}
\item We say that a linear order $\cI$ is {\em strongly $(\mu,
\sigma)$-entangled} if:

($|\cI|,\mu\geq\sigma+\aleph_0$ and) for every $\vare(*)<1+\sigma$,
$t^\zeta_\alpha\in\cI$ (for $\alpha<\mu,\zeta<\vare(*)$) and $u\subseteq
\vare(*)$ such that
\[\alpha<\mu\ \&\ \zeta\in u\ \&\ \xi\in\vare(*)\setminus u\quad\Rightarrow
\quad t^\zeta_\alpha\neq t^\xi_\alpha\]
for some $\alpha<\beta<\mu$ we have:
\begin{description}
\item[(a)] $\vare\in u\quad \Rightarrow\quad t^\vare_\alpha\leq_{\cI}
t^\vare_\beta$,
\item[(b)] $\vare\in\vare(*)\setminus u\quad \Rightarrow\quad t^\vare_\beta
\leq_{\cI} t^\vare_\alpha$.
\end{description}
\item We say that a linear order $\cI$ is {\em strongly positively}
[{\em positively}$^*$] $(\mu,\sigma)$-{\em entangled} if for every
$\vare(*)<1+\sigma$, $t^\zeta_\alpha\in\cI$ (for $\alpha<\mu$, $\zeta<
\vare(*)$) and $u\in\{\emptyset,\vare(*)\}$ for some $\alpha<\beta<\mu$
[for some $\alpha\neq\beta<\mu$] we have
\begin{description}
\item[(a)] $\vare\in u\quad \Rightarrow\quad t^\vare_\alpha\leq_{\cI}
t^\vare_\beta$,
\item[(b)] $\vare\in\vare(*)\setminus u\quad \Rightarrow\quad t^\vare_\beta
\leq_{\cI} t^\vare_\alpha$.
\end{description}
\end{enumerate}
\end{definition}

\begin{remark}
\label{6.12A}
{\em For ``positively$^*$'' it is enough to use $u=\vare(*)$, so only clause
(a) applies. [Why? as we can interchange $\alpha,\beta$.]}
\end{remark}

\begin{proposition}
\label{6.13}
\begin{enumerate}
\item If $\sigma=\aleph_0$ and $\mu=\cf(\mu)\geq\sigma$ then in
Definition \ref{6.12}(1) we can weaken $\alpha<\beta<\mu$ to
$\alpha\neq\beta$ ($<\mu$). [Why? As in \ref{1.2}(6).]
\item For a linear order: ``strongly $(\mu,\sigma)$-entangled'' implies both
``$(\mu,\sigma)$-entangled'' and ``strongly positively $(\mu,\sigma)
$-entangled'' and the last implies ``positively $(\mu,\sigma)$-entangled''.
Lastly ``strongly positively $(\mu,\sigma)$-entangled'' implies ``strongly
positively$^*$ $(\mu,\sigma)$-entangled''.
\item In Definition \ref{6.12} the properties are preserved when
increasing $\mu$ and/or decreasing $\sigma$ and/or decreasing $\cI$.
\item If $\mu=\cf(\mu)$, $(\forall\vare<\sigma)(\forall\theta<\mu)
(\theta^{|\vare|}<\mu)$ then:
\begin{description}
\item[(a)] $\cI$ is $(\mu,\sigma)$-entangled if and only if $\cI$ is
strongly $(\mu,\sigma)$-entangled,
\item[(b)] $\cI$ is positively $(\mu,\sigma)$-entangled if and only if
$\cI$ is strongly positively $(\mu,\sigma)$-entangled.
\end{description}
\item If $\mu=\cf(\mu)$, $(\forall\vare<\sigma)(2^{|\vare|}<\mu)$ then in
Definition \ref{6.12}(1),(2) we can assume
\[(\forall\alpha<\mu)(\forall\zeta<\xi<\vare(*))(t^\zeta_\alpha\neq
t^\xi_\alpha).\]
\item $(\mu,<)$ and $(\mu,>)$ are not strongly $(\mu,2)$-entangled (even if
$\cI$ is strongly $(\mu,2)$-entangled then there is no partial function $f$
from $\cI$ to $\cI$ such that $x<_{\cI} f(x)$, $|\dom(f)|=\mu$ and $f$
preserves $<_{\cI}$).
\item If $(\forall\vare<\sigma)(|\cI|^{|\vare|}<\mu)$ then $\cI$ is strongly
$(\mu,\sigma)$-entangled.
\item Assume $\sigma$ is a limit cardinal, $\cI$ is a linear order. Then
\begin{description}
\item[(a)] $\cI$ is strongly $(\mu,\sigma)$-entangled if and only if for every
$\sigma_1<\sigma$ the order $\cI$ is strongly $(\mu,\sigma_1)$-entangled.
\item[(b)] Similarly for other notions of \ref{1.1}, \ref{6.4},
\ref{6.12}. \QED$_{\ref{6.13}}$
\end{description}
\end{enumerate}
\end{proposition}

\begin{proposition}
\label{6.14}
Assume that $\cI$ is a $(\mu,2)$-entangled linear order and $\theta=\cf(
\theta)<\mu$. Then for some $A\subseteq\cI$, $|A|<\mu$ we have:
\[x<_{\cI}y\ \& \ \{x,y\}\not\subseteq A\quad\ \Rightarrow\quad\ |(x,y)_{\cI}|
\geq\theta.\]
\end{proposition}

\Proof Let $E^\theta_{\cI}$ be the following two place relation on $\cI$:

$x\; E^\theta_{\cI}\;  y$\qquad if and only if

$x=y$ or [$x<_{\cI}y$ \& $|(x,y)_{\cI}|<\theta$] or [$y<_{\cI} x$ \&
$|(y,x)_{\cI}|<\theta$].

\noindent Obviously:
\begin{description}
\item[(a)] $E^\theta_{\cI}$ is an equivalence relation,
\item[(b)] each equivalence class has cardinality $\leq\theta$,
\item[(c)] if an equivalence class has cardinality $\theta$ then there
is a monotonic sequence of length $\theta$ in it.
\end{description}
It is enough to show that the set $A=\{x: |x/E^\theta_{\cI}|>1\}$ has
cardinality less than $\mu$. Suppose that $|A|\geq\mu$. Then there are at
least $\mu$ equivalence classes (as each class is of the size $\leq\theta<
\mu$). Consequently we can find $t^0_\alpha,t^1_\alpha\in\cI$ for $\alpha<\mu$
such that $t^0_\alpha<t^1_\alpha$, $t^0_\alpha\;E^\theta_{\cI}\;t^1_\alpha$
and there is no repetition in $\{t^0_\alpha/E^\theta_{\cI}:\alpha<\mu\}$. For
$\vare(*)=2$ and $u=\{0\}$ we get contradiction to Definition
\ref{6.12}. \QED$_{\ref{6.14}}$ 
\medskip

\noindent{\bf Remark:}\qquad Easily, $\theta=\cf(\theta)$ is redundant.

\begin{proposition}
\label{6.15}
If $\cI$ is strongly $(\mu,\sigma)$-entangled, $\aleph_0\leq\theta<\sigma$,
$\mu\leq |\cI|$ then
\begin{description}
\item[(a)] $2^\theta<\mu$ and
\item[(b)] the cardinal $\chi=:|\{x/E^\theta_{\cI}: x/E^\theta_{\cI}\
\mbox{ not a singleton}\}|$ satisfies $\chi^\theta<\mu$ (where $E^\theta_{\cI}$
is from the proof of \ref{6.14}).
\end{description}
\end{proposition}

\Proof Take $x^0_0>_{\cI} x^1_0<_{\cI} x^0_1<_{\cI} x^1_1$. For $f\in
{}^\theta 2$ let $t^{2\vare+\ell}_f=x^\ell_{f(\vare)}$ (for $\vare<\theta$).
If $2^\theta\geq\mu$ then we may consider $u=\{2\vare:\vare<\theta\}$ and find
(by the strongly $(\mu,\sigma)$-entangledness) functions $f\neq g$ such that
for all $\vare<\theta$:
\[t^{2\vare}_f\leq_{\cI} t^{2\vare}_g,\qquad t^{2\vare+1}_f\geq_{\cI}
t^{2\vare+1}_g.\]
But if $\vare<\theta$ is such that $f(\vare)\neq g(\vare)$ then we get
$x^0_{f(\vare)}<_{\cI} x^0_{g(\vare)}$. Hence $f(\vare)=0$, $g(\vare)=1$ and
consequently $x^1_{f(\vare)}<_{\cI} x^1_{g(\vare)}$. But the last
contradicts to $t^{2\vare+1}_f\geq_{\cI} t^{2\vare+1}_g$. Hence
$2^\theta<\mu$ as required in (a). 

\noindent For (b), let $\langle x^0_i/E^\theta_{\cI}: i<\chi\rangle$ be with
no repetition, $x^1_i\in x^0_i/E^\theta_{\cI}$, $x^0_i<_{\cI} x^1_i$. Let
$\{f_\alpha:\alpha<\chi^\theta\}$ list all functions $f:\theta\longrightarrow
\chi$ (so for each $\alpha<\beta$ there is $\vare<\theta$ such that
$f_\alpha(\vare)\neq f_\beta(\vare)$). Let $t^{2\vare+\ell}_\alpha$ (for
$\alpha<\lambda$, $\vare<\theta$, $\ell<2$) be $x^\ell_{f_\alpha(\vare)}$
and $u=\{2\vare: \vare<\theta\}\subseteq\theta$ (so $\vare(*)=\theta$). If
$\mu\leq\chi^\theta$ then we get $\alpha<\beta<\mu$ by Definition
\ref{6.12}(1) and we get contradiction, so $\mu>\chi^\theta$ as required.
\QED$_{\ref{6.15}}$
\medskip

\noindent{\bf Remark:}\qquad See more in\cite{SaSh:553}.

\begin{proposition}
\label{6.15A}
If $\cI$ is a linear order, $|\cI|\geq\theta\in[\aleph_0,\sigma)$
then $\BA^\sigma_{\inter}(\cI)$ is not $2^\theta$-narrow.
\end{proposition}

\Proof Choose $\cJ\subseteq\cI$, $|\cJ|=\theta$ such that for every $t\in\cJ$
for some $\tau\in\BA^\sigma_{\inter}(\cI)$ we have $\tau\cap\cJ=\{t\}$
($\tau=\bigcap\{[t,s): t<_{\cI}s\in\cJ\}$). Hence for every $\cJ'\subseteq\cJ$
for some $\tau\in\BA^\sigma_{\inter}(\cI)$ we have $\tau\cap\cJ=\cJ'$.
The conclusion now follows. \QED$_{\ref{6.15A}}$

\begin{proposition}
\label{6.16}
In Definition \ref{6.12}(1), if $\cf(\mu)=\mu$ (or less) we can wlog demand
\[(\forall\alpha<\mu)(\forall\zeta<\xi<\vare(*))(t^\zeta_\alpha\neq
t^\xi_\alpha)\]
and for some linear order $<^*$ on $\vare(*)$
\[\zeta<^*\xi, \ \zeta,\xi<\vare(*)\quad\Rightarrow\quad t^\zeta_\alpha<_{\cI}
t^\xi_\alpha.\]
\end{proposition}

\Proof Clearly the new version of the definition implies the old one. So now
assume the old definition and we shall prove the new one. Let $\vare(*)<1+
\sigma$ and $t^\zeta_\alpha$ (for $\alpha<\mu$, $\zeta<\vare(*)$). By
\ref{6.15}(a) we have $\mu>2^{|\vare(*)|}$ so we can replace $\langle
\langle t^\zeta_\alpha: \zeta<\vare(*)\rangle: \alpha<\mu\rangle$ by
$\langle\langle t^\zeta_\alpha:\zeta<\vare(*)\rangle:\alpha\in A\rangle$ for a
suitable $A\in [\mu]^\mu$, and we are done. \QED$_{\ref{6.16}}$

\begin{proposition}
\label{6.17}
\begin{enumerate}
\item Assume $\sigma\geq\aleph_0$ and $\cf(\mu)=\mu$. Then in Definition
\ref{6.12}(1), if we allow first to discard $<\mu$ members of $\cI$ we can add
\begin{description}
\item[(c)] for $\zeta,\xi<\vare(*)$, if $t^\zeta_\alpha<_{\cI} t^\xi_\alpha$
then $\{t^\zeta_\alpha,t^\zeta_\beta\}<_{\cI}\{t^\xi_\alpha, t^\xi_\beta\}$
\end{description}
\noindent (i.e.~we get an equivalent definition).
\item Even without ``if we allow first to discard $<\mu$ members of $\cI$''
part (1) still holds true. It holds also for $(\mu,\sigma)$-entangledness.
\end{enumerate}
\end{proposition}

\Proof 1)\ \ \ The new definition is apparently stronger so we have to prove
that it follows from the old one. Wlog $\mu>\aleph_0$. By \ref{6.14} wlog
\[x<_{\cI} y\quad\Rightarrow\quad |(x,y)_{\cI}|> 2 |\vare(*)|^2.\]
Let $\vare(*)<\sigma$, $t^\zeta_\alpha\in\cI$ be given. By \ref{6.16} we may
assume that
\[(\forall\alpha<\mu)(\forall\zeta<\xi<\vare(*))[t^\zeta_\alpha\neq
t^\xi_\alpha\ \ \&\ \ (t^\zeta_\alpha<t^\xi_\alpha\equiv\zeta<^*\xi)].\]
So for $\alpha<\mu$, $\zeta<^*\xi$ we can choose $s^{\zeta,\xi,\ell}_\alpha
\in\cI$ ($\ell=1,2$) such that for each $\alpha<\mu$:
\[t^\zeta_\alpha<_{\cI} s^{\zeta,\xi,1}_\alpha<_{\cI} s^{\zeta,\xi,2}_\alpha
<_{\cI} t^\xi_\alpha\]
and there are no repetitions in $\{t^\zeta_\alpha,s^{\zeta,\xi,\ell}_\alpha:
\zeta<\vare(*),\zeta<^*\xi\ \mbox{ and } \ell\in\{1,2\}\}$. So for some
$\alpha<\beta$ we have
\[\begin{array}{lcl}
\zeta\in u &\Rightarrow &t^\zeta_\alpha\leq_{\cI} t^\zeta_\beta,\\
\zeta\in\vare(*)\setminus u &\Rightarrow & t^\zeta_\alpha\geq_{\cI}
t^\zeta_\beta,\\
\zeta<^*\xi &\Rightarrow & s^{\zeta,\xi,1}_\alpha\leq_{\cI}
s^{\zeta,\xi,1}_\beta,\\
\zeta<^*\xi &\Rightarrow & s^{\zeta,\xi,2}_\alpha\geq_{\cI}
s^{\zeta,\xi,2}_\beta.
\end{array}\]
Now (c) follows immediately.
\medskip

\noindent 2)\ \ \ Use \ref{6.14} + \ref{6.15}(b).

\noindent For $(\mu,\sigma)$-entangledness -- straightforward.
\QED$_{\ref{6.17}}$

\begin{proposition}
\label{6.18}
Assume $\aleph_0\leq\theta<\sigma$, $\chi=\dens(\cI)$ and:\quad $\cI$ is a
strongly $(\mu,\sigma)$-entangled or $\BA^\sigma_{\inter}(\cI)$ is
$\mu$-narrow. Then $\chi^\theta<\mu$.
\end{proposition}

\Proof First note the following fact:
\begin{claim}
\label{6.18A}
Assume that $\langle\cI_\vare:\vare<\theta\rangle$ is a sequence of pairwise
disjoint convex subsets of $\cI$, $\chi_\vare=\dens(\cI_\vare)\geq\aleph_0$ and
$\chi=\prod\limits_{\vare<\theta}\chi_\vare$. Then $\chi<\mu$ and
$\BA^\sigma_{\inter}(\cI)$ is not $\chi$-narrow.
\end{claim}

\noindent{\em Proof of the claim:}\hspace{0.2in} Choose by induction on
$i<\chi_\vare$, $a^\vare_i<b^\vare_i$ from $\cI_\vare$ such that
$[a^\vare_i,b^\vare_i]_{\cI}$ is disjoint to $\{a^\vare_j, b^\vare_j:
j<i\}$. Let $\{f_\alpha:\alpha<\chi\}$ list $\prod\limits_{\vare<\theta}
\chi_\vare$ (with no repetitions) and let $\vare(*)=\theta$,
$t^{2\vare+\ell}_\alpha$ be:\quad $a^\vare_{f_\alpha(\vare)}$ if $\ell=0$,
$b^\vare_{f_\alpha(\vare)}$ if $\ell=1$, and $u=\{2\vare:\vare<\theta\}$ --
we get a contradiction to ``$\cI$ is strongly $(\chi,\sigma)$-entangled''.
The proof for the Boolean version is similar.
\hfill$\square_{\ref{6.18A}}$
\medskip

\noindent By an argument similar to that of \ref{6.18A} one can show that
$\chi(=\dens(\cI))<\mu$. So the interesting case is when $\chi^\theta>\theta$.
As $2^\theta<\mu$ (see \ref{6.13}) we may assume that $\chi>2^\theta$. Let
$\chi_1=\min\{\lambda:\lambda^\theta\geq\chi\}$, so $\cf(\chi_1)\leq\theta$,
and let $\chi_1=\sum\limits_{\vare<\theta(*)}\chi^+_{1,\vare}$, $\theta(*)=
\cf(\chi_1)\leq\theta$, $\chi_{1,\vare}^\theta<\chi_1$. For each $\vare<
\theta(*)$ we define a two place relation $E^*_\vare$ on $\cI$:
\smallskip

$x\;E^*_\vare\;y$\qquad if and only if

either $x=y$ or $x<y$, $\dens(\cI\rest(x,y))\leq\chi_{1,\vare}$ or

$y<x$, $\dens(\cI\rest(y,x))\leq\chi_{1,\vare}$.
\smallskip

\noindent It is an equivalence relation, each equivalence class has
density $\leq\chi^+_{1,\vare}$, and the number of $E^*_\vare$-equivalence
classes is $\geq\chi$. So by the Erd\"os--Rado theorem we can find a
monotonic sequence $\langle x^\vare_i: i<\theta^+\rangle$ such that $i\neq
j\quad\Rightarrow\quad \neg(x^\vare_i\; E^*_\vare\; x^\vare_j)$.

Without loss of generality, for all $\vare$ the monotonicity is the
same, so wlog $i<j\quad\Rightarrow\quad x^\vare_i<_{\cI} x^\vare_j$.
Throughing away a long initial segment from each $\langle x^\vare_i:i<
\theta^+\rangle$ we may assume that for each $\vare,\zeta<\theta(*)$ either
\[(\forall i<\theta^+)(\exists j<\theta^+)(x^\vare_i<_{\cI} x^\zeta_j\
\&\ x^\zeta_i<_{\cI} x^\vare_j)\]
or
\[\bigcup_{i<j}[x^\vare_i,x^\vare_j]_{\cI},\quad\bigcup_{i<j}[x^\zeta_i,
x^\zeta_j]_{\cI} \quad\quad\mbox{ are disjoint.}\]
Now it is easy to satisfy the assumptions of \ref{6.18A}.
\QED$_{\ref{6.18}}$

\begin{conclusion}
\label{6.19}
Assume that $\mu=\cf(\mu)\geq\sigma=\cf(\sigma)\geq\aleph_0$, $\cI$ is a
linear order of cardinality $\geq\mu$. Then in Definition \ref{6.12}(1) we
can demand $(*)$ of \ref{1.2}(3).
\end{conclusion}

\Proof By \ref{6.18} (and see the proof of \ref{1.2}(3)).
\QED$_{\ref{6.19}}$

\begin{proposition}
\label{6.20}
Assume $\mu=\cf(\mu)>\sigma=\cf(\sigma)\geq\aleph_0$, $\cI$ is a linear order,
$|\cI|\geq\mu$. Then the following conditions are equivalent:
\begin{description}
\item[(a)] $\cI$ is strongly $(\mu,\sigma)$-entangled,
\item[(b)] $\BA^\sigma_{\inter}(\cI)$ is $\mu$-narrow.
\end{description}
\end{proposition}

\Proof (a)$\quad\Rightarrow\quad$ (b)\ \ \ By \ref{6.19} the situation is
similar enough to the one in \ref{1.4} to carry out the proof.
\medskip

\noindent (b)$\quad\Rightarrow\quad$ (a)\ \ \ By \ref{6.18} we can apply the
parallel of \ref{1.2}(3), so the situation is similar enough to the one in
\ref{1.4} to carry out the proof as there. \QED$_{\ref{6.20}}$

\begin{proposition}
\label{6.21}
Assume $\sigma=\theta^+>\aleph_0$, $\mu>\sigma$.
\begin{enumerate}
\item If $\cI$ is not strongly [or strongly positively] [or strongly
positively$^*$] $(\mu_i,\sigma)$-entangled for $i<\theta$ then $\cI$ is not
strongly [or strongly positively] [or strongly positively$^*$]
$(\mu,\sigma)$-entangled for $\mu=\prod_{i<\theta}\mu_i$.
\item If $\mu$ is singular and $\cI$ is strongly $(\mu,\sigma)$-entangled then
for some $\mu'<\mu$, $\cI$ is strongly $(\mu',\sigma)$-entangled (so this
holds for every large enough $\mu'<\mu$).
\item The parallel of (2) holds for ``strongly $(\mu,\sigma)$-entangled'',
``strongly positively $(\mu,\sigma)$-entangled'' and ``strongly positively$^*$
$(\mu,\sigma)$-entangled''.
\end{enumerate}
\end{proposition}

\Proof 1)\ \ \ First we deal with ``strongly $(\mu,\sigma)$-entangled''. We
know $|\cI|\geq\sigma$. Suppose that $\langle t^{i,\zeta}_\alpha:\alpha<
\mu_i, \zeta<\vare_i(*)\rangle$, $u_i$ form a counterexample for $\mu_i$. As
we can extend the sequences and $u_i$ wlog $|u_i|=|\vare_i(*)\setminus u_i|=
\theta$. So by renaming $\vare_i(*)=\theta$, $u_i=\{2\zeta:\zeta<\theta\}$.
Let $f_\beta\in\prod\limits_{i<\theta}\mu_i$ for $\beta<\mu$ be pairwise
distinct. Let us choose $\vare(*)=\theta\cdot\theta$, $t^{\theta i+\zeta}_\beta
=t^{i,\zeta}_{f_\beta(i)}$ for $i<\theta$, $\zeta<\theta$, $\beta<\mu$ and
$u=\{2\zeta: \zeta<\theta{\cdot}\theta\}$. Now check.

For the cases ``strongly positively$^{(*)}$'' $(\mu_i,\sigma)$-entangled first
wlog $\vare_i(*)=\theta$; then as $u_i$ has two values for some $u^*$ we have
\[\prod\{\mu_i: i<\theta,\ u_i=u^*\}=\mu.\]
Thus wlog $u_i=u$ and let $u=\bigcup\limits_{i<\theta}\{i\}\times u^*$. Next
follow as above. 
\medskip

\noindent 2)\ \ \ Assume not. Let $\mu=\sum\limits_{i<\kappa}\mu_i$, $\kappa<
\mu_i<\mu$, $\kappa=\cf(\mu)$, $i<j\ \ \Rightarrow\ \ \mu_i<\mu_j$. So for
each $i<\kappa$ we can find a sequence $\langle t^{i,\zeta}_\alpha:\zeta<
\vare_i(*),\alpha<\mu^+_i\rangle$, $u_i$ exemplifying the failure of ``$\cI$
is $(\mu_i^+,\sigma)$-entangled''. Thus for each $i$ and for every $\alpha<
\mu^+_i$, there are no repetitions in $\{t^{i,\zeta}_\alpha:\zeta<\vare_i(*)
\}$ and $|\cI|\geq \mu$. Now let $\vare(*)=\theta+\theta$, $u=\{2\zeta:
\zeta<\theta+\theta\}$. For $\zeta<\theta$, $\beta\in\mu^+_i\setminus\bigcup
\{\mu^+_j: 0<j<i\}$ we put $t^\zeta_\beta= t^{i,\zeta}_\beta$, $t^{\theta+
\zeta}_\beta=t^{0,\zeta}_i$. Wlog for every $\alpha<\mu$ there are no
repetitions in $\{t^\zeta_\alpha:\zeta<\theta+\theta\}$. Now check.
\medskip

3)\ \ \ Let $\kappa$, $\mu_i$ (for $i<\kappa$), $t^{i,\zeta}_\alpha$ (for
$i<\kappa$, $\zeta<\vare_i(*)$, $\alpha<\mu^*_i$, $u_i$ be as in the proof of
(2) (for the appropriate notion). Again wlog $\vare_i(*)=\theta$, $u_i=u^*$,
but the choice of $\alpha_i$ does not transfer. But by \ref{6.21}(1) we have:
\quad $\mu_i^\theta<\mu$ and hence wlog $(\sum\limits_{j<i}\mu^+_j)^\theta\leq
\mu_i=\mu_i^\theta$, so for some $v_i\subseteq\vare_i(*)$
\[(\zeta\in v_i\quad\Rightarrow\quad [t^{i,\zeta}_\alpha=
t^{i,\zeta}_\beta\ \Leftrightarrow\ \alpha=\beta])\quad\mbox{ and }\quad
[\zeta\in\vare_i(*)\setminus v_i\quad\Rightarrow\quad t^{i,\zeta}_\alpha=
t^{i,\zeta}_\beta]\]
and\quad $t^{i,\zeta}_\alpha= t^{i,\xi}_\alpha\ \Leftrightarrow\
t^{i,\zeta}_\beta=t^{i,\xi}_\beta$. We can omit $\vare_i(*)\setminus v_i$ etc,
so $\langle t^{i,\zeta}_\alpha:\zeta<\vare_i(*),\alpha<\mu_i^+\rangle$ is with
no repetition and proceed as there. \QED$_{\ref{6.21}}$

\begin{remark}
{\em In \ref{6.21}(2) we cannot replace ``strongly entangled'' by
``entangled''; see \cite{SaSh:553}.}
\end{remark}

\begin{conclusion}
\label{6.22}
If $|\cI|\geq\mu>\sigma=\theta^+>\aleph_0$ then:

$\cI$ is $(\mu,\sigma)$-entangled\quad if and only if\quad $\cI$ is strongly
$(\mu,\sigma)$-entangled. \QED
\end{conclusion}

\begin{conclusion}
\label{6.23}
If $\mu\geq\sigma=\theta^+\geq\aleph_0$, $\cI$ is a strongly $(\mu,
\sigma)$-entangled linear order then for some regular $\mu^*$ we have:
\begin{quotation}
\noindent $\mu^*\leq\mu$, $(\forall\alpha<\mu^*)(|\alpha|^\theta<\mu^*)$ and
$\cI$ is (strongly) $(\mu^*,\sigma)$-entangled.
\end{quotation}
Consequently in \ref{6.20} regularity is not needed. \QED
\end{conclusion}
\nocite{M2}

\bibliographystyle{literal-plain}
\bibliography{listb,lista,listc,listx}
\end{document}